\documentclass[reqno]{amsart}

\usepackage{comment}
\usepackage{a4wide}
\usepackage{comment}
\usepackage{color,cite}
\usepackage[normalem]{ulem}
\usepackage{hyperref}                   
\hypersetup{colorlinks,
    linkcolor=blue,%
    citecolor=blue}

\usepackage{amssymb, longtable,mathrsfs, comment, extarrows}
\numberwithin{equation}{section}
\usepackage{amsmath,amsfonts,amssymb,amsthm,amscd,latexsym}

\usepackage{color}
\usepackage{cite}
\usepackage{graphicx}

\usepackage{enumitem}
\usepackage{ulem}
\hypersetup{colorlinks,
	linkcolor=blue,%
	citecolor=blue}	
\DeclareMathOperator{\Tr}{Tr}
\newcommand{\M}{\mathcal M}
\newcommand{\A}{\mathcal A}

\newcommand{\C}{\mathbb C}
\newcommand{\norm}[1]{\left\|#1\right\|}

\newtheorem{theorem}{Theorem}[section]
\newtheorem{lemma}[theorem]{Lemma}

\newtheorem{example}[theorem]{Example}
\newtheorem{prop}[theorem]{Proposition}
\newtheorem{corollary}[theorem]{Corollary}

\newtheorem{remark}[theorem]{Remark}

\newcommand{\cM}{{\mathcal M}}
\newcommand{\cN}{{\mathcal N}}
\newcommand{\cA}{{\mathcal A}}
\newcommand{\cB}{{\mathcal B}}

\newcommand{\cP}{{\mathcal P}}

\newcommand{\cF}{{\mathcal F}}
\newcommand{\cE}{{\mathcal E}}

\newcommand{\cU}{{\mathcal U}}

\begin{document}

\title[Grothendieck noncommutative Marcinkiewicz spaces]{The Grothendieck property of noncommutative Marcinkiewicz spaces}

\author[Y. Gao]{Yi Gao} 
     \address{Institute for  Advanced Study in  Mathematics, Harbin Institute of Technology, Harbin, 150001, China}  \email{{\color{blue}yi.gao@stu.hit.edu.cn}}
\author[J. Huang]{Jinghao Huang}\thanks{Jinghao Huang was supported by the NNSF of China  (No.12031004, 12301160, 12471134 and 12671159).}
    \address{Institute for  Advanced Study in  Mathematics, Harbin Institute of Technology, Harbin, 150001, China}
\email{{\color{blue}jinghao.huang@hit.edu.cn}}   
    \author[R. Xu]{Ran Xu}
     \address{Institute for  Advanced Study in  Mathematics, Harbin Institute of Technology, Harbin, 150001, China}   \email{{\color{blue}xxurann@stu.hit.edu.cn}}


\begin{abstract}
Let $\cM$ be a semifinite von Neumann algebra equipped with a faithful normal semifinite trace. 
We study sufficient  and necessary conditions  of  a noncommutative Marcinkiewicz space affiliated with $\cM $ to have the   Grothendieck property.
 In particular, we show that every noncommutative weak $L_p$-space,
$1<p<\infty$, is Grothendieck. 
We also prove that the intersection   of two Grothendieck noncommutative $\cM$-bimodules is also Grothendieck. 
These answer   several open problems in the literature.

\end{abstract}

\maketitle
\section{Introduction} \label{s1}

A Banach space $X$ is said
to have the {\it Grothendieck property} (or is called  a {\it Grothendieck space}) if each weak $^*$-convergent sequence in the Banach dual $X^*$  is weakly convergent  \cite{PS20}. 
Every reflexive Banach space has the Grothendieck property. 
Grothendieck proved that $C(K)$ has the Grothendieck property when $K$ is an extremally disconnected compact Hausdorff
space \cite[Theorem~9]{G53}. 
In particular, this gives the
Grothendieck property for commutative von Neumann algebras.
Pfitzner 
extended the result to arbitrary von Neumann
algebras 
\cite[Corollary~7]{P94} (see also \cite{GK,FP09,FP}).

Every separable Grothendieck space is reflexive
\cite[Proposition~5.3.10]{M91}. 
In particular, neither $c_0$ nor
$\ell_1$ has the Grothendieck property.  
A classical $L_1$ space has
the Grothendieck property if and only if it is finite-dimensional
\cite[discussion following Proposition~2.3]{Wnuk}. 

Unless otherwise stated, $\Phi:[0,\gamma)\to[0,\infty)$, $\gamma\in(0,\infty]$, 
is a nonzero nondecreasing concave function with $\Phi(0)=0$.
	The classical \emph{Marcinkiewicz space} $M_\Phi(0,\gamma)$  associated with $\Phi$  
consists of the  real- or complex-valued measurable functions $f$
on $(0,\gamma)$ for which
\[
\norm{f}_{M_\Phi(0,\gamma)}
:=
\sup_{0<t<\gamma}
\frac{1}{\Phi(t)}\int_0^t\mu(s;f)\,ds
<\infty,
\]
where $\mu(f)$ is the decreasing rearrangement of $|f|$.  
	If $1<p<\infty$ and $\Phi(t)=t^{1-1/p}$ for $0<t<\gamma$,
	then the corresponding Marcinkiewicz space is the \emph{weak-$L_p$
		space $L_{p,\infty}(0,\gamma)$}.
	If $\Phi(t)=t$ for $0<t<\gamma$, then  $M_\Phi(0,\gamma)=L_\infty(0,\gamma)$ isometrically.

Lotz proved that $L_{p,\infty}(0,\gamma)$ is Grothendieck for
$1<p<\infty$ by applying his sufficient criterion for the
Grothendieck property of Banach lattices
\cite[Theorems~1 and~2]{L10}. 
 Applying the same criterion to
Marcinkiewicz spaces, 
 de Pagter and
	Sukochev \cite[Theorem~5.6]{PS20}  showed  that
 $M_\Phi(0,\gamma)$ has
the Grothendieck property if and only if either  condition (A)
 or (B)  of Theorem~\ref{thm:main} below
holds. 

A natural question in noncommutative integration is which
Banach space properties pass from a symmetric function space
to its associated operator space.  
It is known that the separable part of the weak $L_p$-function space $L_{p,\infty}(0,\infty )$ is $M$-embedded \cite{HWW93}.  
As shown in \cite{HLS20}, the separable part of noncommutative weak $L_p$-space $L_{p,\infty}(\cM,\tau)$, $1<p<\infty$, is $M$-embedded. 
However, \begin{quote}
    It is still unknown whether the noncommutative counterpart $L_{p,\infty }(\cM,\tau)$, $1 < p < \infty$, is a Grothendieck space or not\cite[Remark 4.5]{HLS20} (see also \cite[Remark 3.4]{HL}).
\end{quote}
De Pagter and Sukochev asked a more general question:
\begin{quote}
    It seems to be an interesting open question
whether the results presented in this article also hold for noncommutative Marcinkiewicz
spaces \cite[Introduction]{PS20}.
\end{quote} 
Below, we answer this question by showing  that 
 their sufficient condition extends to the setting of    Marcinkiewicz spaces over
arbitrary semifinite von Neumann algebras, see Theorem~\ref{thm:main} below.
If, in addition,   the underlying  von Neumann algebras are non-atomic, then their necessary condition remains valid,  see Corollary \ref{cor:characterization} below.


 Let $\M$ be a semifinite von Neumann algebra on a Hilbert space
$H$, with identity ${\bf 1}$ and a fixed faithful normal
semifinite trace $\tau$. 
	The \emph{noncommutative Marcinkiewicz space} associated with
	$\Phi$ is defined by
	\begin{align*}
		M_\Phi(\M,\tau)
		&=\{x\in S(\M,\tau): \mu(x)\in  M_\Phi(0, \tau({\bf 1}))\},
	\end{align*} and 
    $\norm{x}_{M_\Phi(\M, \tau)} =\norm{\mu(x)}_{M_\Phi(0, \tau({\bf 1}))}.$




\begin{theorem}\label{thm:main}
Let  $(\cM, \tau)$ be a nonzero semifinite  von Neumann algebra. 
Suppose that one of the following conditions holds: 
\begin{enumerate}[label=(\Alph*)]
\item $\tau(\mathbf1)=\infty$ and 
\begin{align} \label{condition 1}
 \mathop{\liminf}\limits_{t\downarrow0}\frac{\Phi(2t)}{\Phi(t)}>1,
 \qquad
 \mathop{\liminf}\limits_{t\to\infty}\frac{\Phi(2t)}{\Phi(t)}>1;
\end{align}
\item  $\tau(\mathbf1)<\infty$ and 
\begin{align} \label{condition 2}
 \mathop{\liminf}\limits_{t\downarrow0}\frac{\Phi(2t)}{\Phi(t)}>1.
\end{align}
\end{enumerate}
Then, $M_\Phi(\M,\tau)$ has the Grothendieck property.
\end{theorem}

Our proof of Theorem~\ref{thm:main} uses a different approach from that of
de Pagter and Sukochev~\cite{PS20} and does not rely
on the Banach lattice structure.
We construct an isomorphic embedding of the Lorentz
predual $\Lambda_\Phi(\cM,\tau)$ into a Banach space $Z$
whose dual $Z^*$ has the Grothendieck property.
The adjoint map is then a bounded linear surjection
from $Z^*$ onto $M_\Phi(\cM,\tau)$.
Since the Grothendieck property is preserved under
bounded linear surjections
\cite[Proposition~3.1.4]{GK}, it follows that 
$M_\Phi(\cM,\tau)$ is a Grothendieck space.
One of the ingredients in the proof of Theorem~\ref{thm:main} is the Grothendieck property of the direct sum of  varying summands established in  Section~\ref{S:3} (see Corollary~\ref{cor:uc-product} below). 
Kucher  \cite[Conjecture 1]{K99}  (see also \cite[Problem 29]{GK}) asked 
 {\it whether  $\bigoplus_{k\in \mathbb{N}}^{\ell_\infty} E$ is Grothendieck whenever $E$ is a super-reflexive Banach space.}
One of the consequences of results in Section \ref{S:3} (Corollary~\ref{cor:uc-product} below) is  an affirmative answer to Kucher's question. 
During the preparation of the current manuscript, 
Li, Luo, Wu and
Zhou answered this question independently in a more general setting of Pe{\l}czy\'nski’s property (V)
 \cite[Theorem 1.1]{LLWZ26}.

Combined with the commutative
characterization of de Pagter and Sukochev
\cite[Theorem~5.6]{PS20}, Theorem~\ref{thm:main} above shows that the
Grothendieck property of $M_\Phi(0,\tau(\mathbf{1}))$ implies
that of $M_\Phi(\M,\tau)$ for every semifinite von Neumann algebra
$(\M,\tau)$. 
The converse may fail, 
as seen in Example~\ref{atomic} below.


For $1<p<\infty$, the noncommutative weak $L_p$-space
is defined by
\[
L_{p,\infty}(\cM,\tau)
:=
\left\{
x\in S(\cM,\tau): \norm{x}_{L_{p,\infty}(\cM,\tau)}
<\infty
\right\},
\] 
where $\norm{x}_{L_{p,\infty}(\cM,\tau)}
:= 
\mathop{\sup}\limits_{t>0}t^{1/p-1}\int_0^t\mu(s;x)\,ds$
 (see, e.g., \cite[Remark 6.7.9]{DPS}). 
Thus,  $L_{p,\infty}(\cM,\tau)=M_\Phi(\cM,\tau)$
isometrically when $\Phi(t)=t^{1-1/p}$.

\begin{corollary}\label{weakLp}
   Let $(\cM, \tau)$ be a semifinite von Neumann algebra,
and let $1<p<\infty$.  
    Then, weak $L_p$-space $L_{p,\infty}\left(\cM,\tau\right)$ is  Grothendieck.
\end{corollary}
Corollary~\ref{weakLp} extends \cite[Theorem 3]{L10} and \cite[Corollary~6.4]{PS20}
to the noncommutative setting and gives an affirmative
answer to the question raised in \cite{HLS20}
and   \cite[Remark~3.4]{HL}.

 For non-atomic $\M$, Theorem~\ref{thm:main} yields the following full characterization of noncommutative Marcinkiewicz spaces having the Grothendieck property. 
\begin{corollary}\label{cor:characterization}
Let $(\cM, \tau)$ be a nonzero semifinite  non-atomic von Neumann algebra.
Then, $M_\Phi(\M,\tau)$ has the Grothendieck property
if and only if either  condition (A)
or (B) in Theorem~\ref{thm:main} holds.
\end{corollary}


The study of the Grothendieck property of the intersection between Banach spaces plays an important role in the proof of Theorem \ref{thm:main}. 
Recently, it was shown in \cite[Theorem 2.13]{HL}   that
$E\cap L_\infty(0,\infty)$ is Grothendieck whenever $E$ is a reflexive Banach function space on $(0,\infty)$ and satisfies a lower $q$-estimate for some $1 \le q <\infty$, which leads to the following  
  questions \cite[Questions]{HL}:
\begin{quote}
 \begin{enumerate}
     \item Assume that $E(0,\infty)$ is a reflexive (not necessarily symmetric) Banach function space. 
Must $E\cap L_\infty(0,\infty)$ be a Grothendieck space? 
More generally, what about $E\cap M_\Psi(0,\infty)$ for any Marcinkiewicz space $M_\Psi(0,\infty)$ with weakly sequentially complete dual? 
     \item Assume that $E(0,\infty)$ and $F(0,\infty)$ are Banach function spaces that are both Grothendieck spaces. 
     Must $E\cap F$ be a Grothendieck space?
 \end{enumerate}   
\end{quote}
In Section \ref{s6},
we show  that  a Banach $\M$-bimodule $Y$ is Grothendieck if and only if it is positively
Grothendieck (see Theorem 
 \ref{criterion} below), which establishes a noncommutative version of \cite[Proposition 1]{Kuhn80}. 
 As a consequence, 
we prove that the intersection $E\cap F$ of two Grothendieck noncommutative $\cM$-bimodules is also Grothendieck (see Theorem \ref{main}  below).
This, in particular,  gives   affirmative answers to  Questions (1) and (2) (see Corollary \ref{affirmative} below).

\section{Preliminaries}
\label{s2}
\subsection{Basic facts and notation}	
 Let $\mathbb R$ and $\mathbb C$ denote the fields of
real and complex numbers, respectively. 
Set $\mathbb N:=\{0,1,2,\ldots\}$.

  A positive measure on a measurable space $(\Omega,\Sigma)$
is a countably additive map $\lambda:\Sigma\to[0,\infty]$
with $\lambda(\varnothing)=0$. It is finite if
$\lambda(\Omega)<\infty$, and it is a probability measure if
$\lambda(\Omega)=1$ \cite[Section~1.1]{D19}.

For a compact Hausdorff space $K$, let $\mathcal B(K)$ denote
the Borel $\sigma$-algebra generated by the open subsets of $K$.
A complex Borel measure on $K$ is a countably additive map
$\nu:\mathcal B(K)\to\mathbb C$ with $\nu(\varnothing)=0$.
Its total variation is given by
\[
|\nu|(B):=\sup\left\{
  \sum_{j=1}^{m}|\nu(B_j)|:
  B=\bigsqcup_{j=1}^{m}B_j,\quad B_j\in\mathcal B(K)
\right\},
\quad B\in\mathcal B(K),
\] 
where $\bigsqcup$ denotes a disjoint union. 
We call $\nu$ {\it regular} if, for every $B\in\mathcal B(K)$,
\[
|\nu|(B)
=\sup_{\substack{F\subseteq B\\ F\text{ compact}}}|\nu|(F)
=\inf_{\substack{B\subseteq U\\ U\text{ open}}}|\nu|(U).
\]

For a Banach space $(X,\norm{\cdot}_X)$, we denote
its closed unit ball, unit sphere, and Banach dual
by $B_X$, $S_X$, and $X^*$, respectively. 
The weak topology on $X$ and the weak$^*$ topology on $X^*$
are denoted by $\sigma(X,X^*)$ and $\sigma(X^*,X)$, respectively.

Throughout this paper, a subspace means a norm-closed
linear subspace. 
A  subspace $Y$ of  $X$
is called {\it complemented} if there exists a bounded
linear operator $p:X\to X$ such that
$p^2=p$ and $p(X)=Y$. 
	For a bounded linear operator $T:X\to Y$ between Banach spaces,
	we denote its adjoint by $T^*:Y^*\to X^*$, defined by
    \begin{align*}
	T^*(f):=f\circ T,\qquad f\in Y^*.
	\end{align*}

For a family of Banach spaces $(X_k)_{k\in\mathbb N}$
and $1\le p\le\infty$, define
\[
\mathop{\bigoplus}\limits_{k\in\mathbb N}^{\ell_p}X_k
:=
\left\{
x=(x_k): x_k\in X_k,\ 
\norm{x}:=
\left\|(\norm{x_k}_{X_k})_{k\in\mathbb N}\right\|_{\ell_p}
<\infty
\right\}.
\]

Unless otherwise stated, whenever $X$ and $Y$
are Banach spaces continuously embedded in a common
Hausdorff topological vector space, their intersection $X\cap Y$, equipped with the  intersection norm 
\[
x
\mapsto
\max\{\norm{x}_X,\norm{x}_Y\},
\]
is a Banach space
\cite[Lemma~2.3.1]{BL12}.
Consequently, the diagonal map
\begin{align} \label{diagonal map}
X\cap Y\to X\bigoplus^{\ell_\infty} Y,
\qquad x\mapsto(x,x),
\end{align}
is an isometric isomorphism onto the (closed) 
subspace $\{(x,x):x\in X\cap Y\}$.

For a normed space $X$
and a (closed) 
subspace $Y\subseteq X$,
the quotient $X/Y$ is equipped with the quotient norm
\begin{align}\label{defquotientnorm}
\norm{x+Y}_{X/Y}:=
\inf_{y\in Y}\norm{x-y}_X,
\qquad x\in X.
\end{align}
If $X$ is Banach, then so is $X/Y$;
	see \cite[Theorem 1.41 and  p.~32]{R91}.


	\subsection{Symmetric function spaces}
	
	Let \((I, \Sigma, m)\) denote the measure space  \( I = (0, a) \) for some \( a \in (0, \infty] \), endowed with the Lebesgue measure \( m \).
	Denote by \( L(I, m) \) the space of   all (equivalence classes of) $\Sigma$-measurable complex-valued 
	functions on \( I \).  
	Denote by  $L_p
	$
	the space of all $p$-integrable real-valued (or complex-valued) 
	$\Sigma$-measurable functions on $I
	$ ($1\leq p<\infty$)
	and by  $L_{\infty}
	$ 
	the space
	of all $m$-essentially bounded real-valued  (or complex-valued)  $\Sigma$-measurable functions on $I
	$, 
	equipped with the usual norms
	$$\norm{f}_{L_p}=\left( \int_{I} |f|^{p} dm \right)^{\frac{1}{p}}, \quad 1 \le p<\infty, \quad   \norm{f}_{L_{\infty}}= \operatorname*{ess\,sup}_{t\in I}|f(t)|,   
	$$
respectively.
	


	The algebra \( L_\infty (I) \) can be seen as an abelian von Neumann algebra acting via multiplication on the Hilbert space \(  L_2 (I) \), 
	with the trace given by integration with respect to Lebesgue measure \( m \).  
    For $f\in L(I)$, its distribution function is defined by
\[
d(s;f):=m\bigl(\{t\in I:|f(t)|>s\}\bigr),
\qquad s>0.
\] 
The right-continuous inverse (called the decreasing rearrangement of $f$) of $d(|f|)$, $f\in L(I)$, is given by
$$\mu(t;f)=\inf\{s\ge0:d\left(s;|f\right|)\le t\}, \qquad t\ge0.$$
The algebra of all $\tau$-measurable operators affiliated
with $L_\infty(I)$ can be identified with
\[
S(I):=\{f\in L(I,m):d(s;f)<\infty
       \text{ for some }s>0\};
\]
see~\cite{KPS,BS,LT79}. 


        A Banach space $\left(E,\norm{\cdot}_E\right)\subseteq S\left(I\right)$ is said to be a {\it Banach function space} if, from $f\in E, ~g\in S\left(I\right)$ and $|g|\le |f|$, it follows that $g\in E$ and $\norm{g}_E\le \norm{f}_E$ \cite{KPS,LT79}. 
	A Banach  function space $(E,\norm{\cdot}_E) \subseteq S(I)$ is said to be 
	{\it symmetric} if 
	$x\in E$,  $y\in S(I)$  and  $\mu(y)\leq\mu(x)$ imply that  $y\in E$ and $\left\|y\right\|_E\leq \left\|x\right\|_E.$ If, in addition,   
	\(\| y \|_E \leq \| x \|_E\) whenever  \( x, y \in E \) satisfy  $y\prec \prec x$ (that is, $\int_0^t\mu(s;y)ds\le \int_0^t \mu(s;x)ds$ for all $t>0$),  	then 
	$(E, \norm{\cdot}_{E})$ is said to be 	{\it  strongly symmetric}. 
	A symmetric  function space $(E, \norm{\cdot}_{E})$  is said to be 	{\it  fully symmetric}  if  \(\norm{\cdot}_E\) has the additional property that    for any $x\in E$ and $y\in S(I)$,  $y\prec \prec x$  implies that $y\in E$ and $\norm{y}_E\le\norm{x}_E$.

	The {\it K\"othe dual space} \(E^{\times}\) of a symmetric
	function space \(E\) is the collection of all elements \(y \in S(I)\) such that
	\[
	\|y\|_{E^{\times}} = \sup\left\{\int_I  |xy|\,dm : \|x\|_{E} \le 1\right\} < \infty.
	\]
  \subsection{Construction of symmetric operator spaces}   

Let $\M$ be a semifinite von Neumann algebra with identity
${\bf 1}$, equipped with a fixed faithful normal semifinite
trace $\tau$. 
We denote its operator norm by $\norm{\cdot}_{\M}$ or $\norm{\cdot}_\infty$.

 A closed densely defined linear operator $x$ on $ H$ is
\emph{affiliated} with $\M$ if $ux=xu$ for every unitary operator
$u\in\M'$, where $\M'$ is the commutant of $\M$.
In this case, the spectral projection $e^{|x|}(B)$ of
$|x|=(x^*x)^{1/2}$ belongs to $\M$ for every Borel set
$B\subseteq[0,\infty)$. 
Such an operator
 $x$ is
\emph{$\tau$-measurable} if
$\tau(e^{|x|}(s,\infty))<\infty$ for sufficiently large $s>0$. 
	Denote by 	$S(\mathcal{M},\tau)$  the set of all $\tau$-measurable operators \cite{LSZ,DPS,LSZ21,LMSZ23}. 
For $x\in S(\M,\tau)$, its
    \emph{generalized singular
		value function} is defined by
	\[
	\mu(t;x)
	=\inf\left\{s\ge0:
	\tau\bigl(e^{|x|}(s,\infty)\bigr)\le t\right\},
	\qquad t\ge0.
	\]


For a linear subspace $\cE\subseteq S(\cM,\tau)$, we 
write $\cE_+:=\{x\in\cE:x\ge0\}$ for  positive cone of $\cE$.
The extended trace $\tau : S(\M, \tau)_+ \to [0, \infty]$ is given by \cite[Proposition 3.3.3]{DPS}
	\[
	\tau(x) = \int_0^\infty \mu(t; x)  dt, \quad x \in S(\M, \tau)_+. 
	\]
	
    
    Let $(\cM_1, \tau_1)$ and $(\cM_2, \tau_2)$ be two semifinite von Neumann algebras, respectively. If  $x\in S\left(\cM_1,\tau_1\right)$ and $y\in S\left(\cM_2,\tau_2\right)$, then $x$ is said to be {\it submajorized} by $y$, denoted by $x\prec\prec y$, if $\mu(x)\prec \prec \mu(y)$ \cite[Section 3.3]{LSZ}.

    A Banach operator space \( (\cE, \norm{\cdot}_{\cE})\subseteq S(\cM, \tau) \)  will be  called
	\begin{enumerate}
		\item  {\it symmetric} 
		if \( x \in S(\cM, \tau) \), \( y \in \cE \) and \( \mu(x) \leq \mu(y) \) imply that \( x \in \cE \) and \(\| x \|_{\cE} \leq \| y \|_{\cE}\). 
		\item  {\it strongly symmetric}   if  \(\norm{\cdot}_{\cE}\) has the additional property that \(\| x \|_{\cE}\leq \| y \|_{\cE}\) whenever \( x, y \in \cE \) satisfy \( x \prec \prec y \).
		\item  {\it fully symmetric} if    \(\norm{\cdot}_{\cE}\) has the additional property that \( x \in S(\cM, \tau) \), \( y \in \cE \) and \( x \prec \prec y \) imply that \( x \in \cE \) and \(\| x \|_{\cE} \leq \| y \|_{\cE} \).
	\end{enumerate}

	Let
	$(E, \norm{\cdot}_{E})$ be a symmetric  function   space on $(0, \tau({\bf 1}) )$.
	Then, we obtain a symmetric operator space (operator counterpart) by the following  construction: 
	\begin{align*}
		E(\mathcal{M},\tau):=\Big\{x\in S(\mathcal{M},\tau):\ \mu(x)\in E \Big\}
		, \quad  
		\left\|x\right\|_{E (\mathcal{M},\tau)}:= \left\|\mu(x)\right\|_{E}.
	\end{align*}
	Many properties 
	carry over from   \( E \) to its operator   counterpart \( E(\mathcal{M}, \tau) \) \cite[Theorem 6.1.7]{DPS}. 
	When $E=L_p(0, \tau(\mathbf1) )$ for  $1 \leq p \leq \infty$, the associated space $E(\M, \tau)$ coincides with the noncommutative $L_p$-space  $L_p(\M, \tau)$. In particular,  the von Neumann algebra
	$\M$ is also denoted by 
	$L_\infty (\M, \tau)$ \cite[Remark 3.4.15]{DPS}. 


	A  symmetric function/operator space
	\((\cE, \norm{ \cdot }_{\cE}) \) is  said to have the {\it Fatou property}  if for every upward directed net  \( \{x_\alpha\} \) in \( \cE_{+} \), satisfying \( \sup_\alpha \|x_\alpha\|_{\cE} < \infty \), there exists an element \( x \in \cE_{+} \) such that \( 0\leq x_\alpha \uparrow x \) in \( \cE \) and \( \|x\|_{\cE} = \sup_\alpha \|x_\alpha\|_{\cE} \) \cite[Definition 4.1.21]{DPS}.	
	 Note that $M_\Phi\left(\cM,\tau\right)$ is a fully symmetric space with the Fatou property \cite[Theorem 6.4.9]{DPS}. 
	
Recall that a linear subspace $\left(Y,\norm{\cdot}_Y\right)$ of $S(\cM,\tau)$  is called a {\it normed $\cM$-bimodule} if $x\in Y$ and $u,v\in \cM$ imply that $uxv\in Y$ and 
$$\norm{uxv}_Y\le \norm{u}_\infty\norm{v}_\infty\norm{x}_Y.$$
It is a \emph{Banach $\M$-bimodule} if it is complete in this norm.
    
	Suppose that  \( (X, \norm{\cdot}_{X})  \subseteq S(\M, \tau) \) is a normed   \( \M \)-bimodule and  \( p  \) is a nonzero projection  in $\M$. 
    If $x\in \cM$, then $x_p:=(pxp)|_{p(H)}$.
    In this case, $\cM_p$ is a unital $*$-subalgebra of $B\left(p(H)\right)$ and the mapping $\phi_p:p\cM p\to \cM_p$, defined by $\phi_p(x)=x_p$ for all $x\in p\cM p$, is a surjective unital $*$-isomorphism \cite[p. 34]{DPS}.
    We identify the $*$-algebra $p\cM p$ with the von Neumann algebra $\cM_p$ on $p(H)$.
	The reduced von Neumann algebra \( p\M p \) is equipped with the trace \( \tau_p \) given by \( \tau_p(a_p) = \tau(pap) \),  
	\( a \in \M_{+} \). 
    Define 
    the reduced spaces: 
	\[
	p Xp = \{pxp : x \in  X\}, \quad  X_p = \{x_p \in S(p\M p,  \tau_p) : x \in  X\}. 
	\]
	Then,   
	the mapping  \( x \mapsto x_p \) is a \(*\)-preserving linear isomorphism from \( p Xp \) onto the   \( p\mathcal Mp \)-bimodule \(  X_p \) of \( \tau_p \)-measurable operators   \cite[Lemma 3.7.1]{DPS}. 
	Moreover,  endowing 
	\(  X_p \) 
	with the norm
	\[
	\|x_p\|_{ X_p} := \|pxp\|_ X, \quad x \in  X,   
	\]
	the space 
	\(  X_p \) 
	becomes a normed 
    \( 
    p\mathcal Mp \)-bimodule and 
	the mapping  \( x \mapsto x_p \) is an isometry from \( p Xp \) onto \(  X_p \) \cite[p.~241]{DPS}.



\subsection{Lorentz and Marcinkiewicz spaces}

We recall the noncommutative Lorentz space associated with $\Phi$ and its duality with the Marcinkiewicz space $M_\Phi(\M,\tau)$.

Let $\Phi$ be as in Section \ref{s1},  and $\gamma:=\tau(\mathbf{1})\in(0,\infty]$.   
If $\gamma<\infty$, define
$\Phi(\gamma):=\mathop{\lim}\limits_{t\uparrow\gamma}\Phi(t)$
and set $\Phi(t)=\Phi(\gamma)$ for $t>\gamma$.
We write
\[
    \Phi(0_+):=\lim_{t\downarrow0}\Phi(t),
    \qquad
    \Phi(\infty):=\lim_{t\to\infty}\Phi(t).
\] 
Note that 
$ 
1\leq \frac{\Phi(2t)}{\Phi(t)}\le 2, $
 $ t>0 
$ 
 \cite[p. ~60]{DPSS98}. 

The right-hand derivative $\Phi'_+(t)$ exists for all
$t\in(0,\gamma)$ \cite[Theorem $\langle 1\rangle$]{P02}, 
defined by
\[
    \Phi'_+(t)
    :=\lim_{h\downarrow0}
    \frac{\Phi(t+h)-\Phi(t)}{h},
    \qquad t>0, 
\]  and $\Phi$ may be written as \cite[Section 6.3]{DPS}
\begin{align}\label{intergralphi}
\Phi(t)=\Phi(0_+)+\int_0^t\Phi'_+(s)\,ds,
\qquad t>0. 
\end{align} 
Note that the function $\Phi'_+$ is decreasing \cite[Theorem 1.3.3]{NP06}, 
$\Phi'_+(t)\ge0$ on $(0,\gamma)$ and the limit
$ 
0\leq \mathop{\lim}\limits_{t\to\infty}\Phi_+'(t)<\infty
$.

For a projection $e\in \M$ with $0<\tau(e)<\infty$, 
$\mu(e)=\chi_{[0,\tau(e))}$. 
Since  $t/\Phi(t)$ is
nondecreasing, it follows from \cite[p. ~406]{DPS} that 
\begin{equation}\label{eq:projection-norm}
\norm{e}_{M_\Phi(\cM,\tau)}
=\frac{\tau(e)}{\Phi(\tau(e))},
\qquad 0<\tau(e)<\infty.
\end{equation}

For $f\in S(0,\gamma)$,  define
\begin{equation}\label{deffLambdaPhi}
\norm{f}_{\Lambda_\Phi}
=
\Phi(0_+)\norm{f}_\infty
+\int_0^\gamma\mu(t;f)\Phi'_+(t)\,dt.
\end{equation}
The \emph{Lorentz function  space}
$\Lambda_\Phi(0,\gamma)$
is now defined by setting
\[
\Lambda_\Phi
=
\left\{
f\in S(0,\gamma):
\norm{f}_{\Lambda_\Phi}<\infty
\right\}.
\]
The {\it noncommutative Lorentz space}  associated with $\Phi$  is defined by
\[
\Lambda_\Phi(\M, \tau)
=
\left\{
x\in S(\M, \tau):\mu(x)\in\Lambda_\Phi(0,\gamma)
\right\}  \]
and   
$\norm{x}_{\Lambda_\Phi(\M, \tau)}
=
\norm{\mu(x)}_{\Lambda_\Phi}$ for all  $
 x\in\Lambda_\Phi(\M, \tau).
$
The space $\Lambda_\Phi(\M, \tau)$  is a fully symmetric space with the Fatou property \cite[Theorem 6.3.8]{DPS}.

The K\"othe dual of $\Lambda_\Phi\left(\cM,\tau\right)$ is defined by setting \cite[Definition 4.3.1]{DPS}
$$\Lambda_\Phi\left(\cM,\tau\right)^{\times}:=\left\{y\in S\left(\cM,\tau\right): \norm{y}_{\Lambda_\Phi\left(\cM,\tau\right)^\times} <\infty\right\},$$
where 
$$\norm{y}_{\Lambda_\Phi\left(\cM,\tau\right)^\times}=\sup\left\{\tau\left(|xy|\right):x\in \Lambda_\Phi\left(\cM,\tau\right),\norm{x}_{\Lambda_\Phi\left(\cM,\tau\right)}\le1\right\}.$$
If, in addition, $\Lambda_\Phi\left(\cM,\tau\right)$ has order continuous norm, then by \cite[Theorem 5.2.9, Proposition 5.3.2 and Theorem 6.4.9]{DPS}, $\Lambda_\Phi\left(\cM,\tau\right)^*=M_\Phi\left(\cM,\tau\right)$.

We next examine the growth of $\Phi$ under the hypotheses of
Theorem~\ref{thm:main}. 
We first note that if $\mathop{\liminf}\limits_{t\downarrow0}\frac{\Phi(2t)}{\Phi(t)}>1$,
then $\Phi(0_+)=0$. 
Indeed, if $\Phi(0_+)>0$, then
\[
\lim_{t\downarrow0}\frac{\Phi(2t)}{\Phi(t)}
=\frac{\Phi(0_+)}{\Phi(0_+)}=1,
\]
a contradiction.


\begin{lemma}\label{order continuous}
Assume that $\tau(\mathbf 1)=\infty$.
 Suppose that   condition \eqref{condition 1} of Theorem \ref{thm:main} holds.
Then, there exists $\alpha\in(1,2]$ such that, with
$\delta:=\log_2\alpha$,
\begin{align*}
\begin{aligned}
    \alpha^{-1}\lambda^\delta\Phi(t)
    &\le\Phi(\lambda t)\le\lambda\Phi(t),
    &&\lambda\ge1,\ t>0.
\end{aligned}
\end{align*}
Moreover, $\Phi$ is strictly increasing and $\Phi(\infty)=\infty$. 
\end{lemma}
\begin{proof}
        If $\lim_{t\to\infty}\Phi(t)<\infty$, then $\lim_{t\to\infty} \frac{\Phi(2t)}{\Phi(t)}
        =1$, which is a contradiction. Hence,   $\lim_{t\to\infty}\Phi(t)=\infty$.
        We claim that $\Phi$ is strictly increasing.
        Indeed, assume by contradiction that 
        $\Phi(a)=\Phi(b)$ with $a< b<\infty$.
      Since $\Phi$ is concave and nondecreasing, it follows that
$\Phi(t)=\Phi(b)<\infty$ for all $t>b$, which contradicts $\lim_{t\to\infty}\Phi(t)=\infty$.
        This proves the claim.

        Observe that 
          there exists $\alpha>1$ such that (see e.g. the proof of \cite[Proposition 5.1]{PS20})
    \begin{align}\label{fracphi2tgeqalpha}
 \frac{\Phi(2t)}{\Phi(t)}\ge \alpha>1, \qquad t\in \left(0,\infty\right).
    \end{align}
For every $\lambda\ge1$ and $t>0$, concavity of $\Phi$ and
$\Phi(0)=0$ give
\begin{align}\label{Phitgeqlambda1}
    \Phi(t)
    =\Phi\left(\frac1\lambda(\lambda t)
      +\left(1-\frac1\lambda\right)0\right)
    \ge\frac1\lambda\Phi(\lambda t).
\end{align}
Since $\Phi(t)>0$ for $t>0$, it follows that 
\[ \alpha\Phi(t)
\overset{\eqref{fracphi2tgeqalpha}}{\le}
\Phi(2t)\overset{ \eqref{Phitgeqlambda1},~ \lambda=2}{\leq }2\Phi(t), 
\]
which implies that  $\alpha\le2$.

Set $\delta=                                      \log_2\alpha\in(0,1]$.
Given $\lambda\ge1$, choose an integer $m\ge0$ such that
$2^m\le\lambda<2^{m+1}$.
By the monotonicity of $\Phi$, we have 
\begin{align} \label{lambda t ge t}
    \Phi(\lambda t)
    \ge\Phi(2^m t)
    \overset{\eqref{fracphi2tgeqalpha}}{\geq } \alpha^m\Phi(t).
\end{align}
Since
$   \lambda^\delta
< 2^{(m+1)\delta}
    \overset{ \delta=\log_2\alpha}{=} \alpha^{m+1},
$
it follows that $ \frac{1}{\alpha}\lambda^\delta < \alpha^{m}$, and thus 
\[
  \lambda \Phi( t)\overset{\eqref{Phitgeqlambda1}}{\geq }  \Phi(\lambda t)
    \stackrel{\eqref{lambda t ge t}}{\ge}  \alpha^m\Phi(t) > \alpha^{-1}\lambda^\delta\Phi(t).\]
\end{proof}


The following extension of $\Phi$ from $(0,\gamma)$ to $(0,\infty)$ is certainly known to experts. For completeness, we include a full proof below. 
\begin{lemma}\label{lem:growth}
Assume that $\tau(\mathbf 1)<\infty$.
Suppose that   condition \eqref{condition 2} of Theorem \ref{thm:main} holds. 
Then,  there exists an increasing concave function
$\Psi:[0,\infty)\to[0,\infty)$ with $\Psi(0)=0$
satisfying
\[\liminf_{t\downarrow0}\frac{\Psi(2t)}{\Psi(t)}>1,
\qquad
\liminf_{t\to\infty}\frac{\Psi(2t)}{\Psi(t)}>1
\]
and
\begin{align*} 
 \Phi(t)\le\Psi(t)\le 2\Phi(t),\quad 0<t<\gamma.
\end{align*} 
\end{lemma}

\begin{proof}
Set $\Phi(\gamma):=\mathop{\lim}\limits_{t\uparrow\gamma}\Phi(t)$,
which is finite (by concavity) and positive, and define
\begin{align}\label{defPsiequiv}
\Psi(t)=\Phi(\min\{t,\gamma\})
+\frac{\Phi(\gamma)}{\gamma}t,
\qquad t\ge0,
\end{align}
We note that  $\Psi$ is strictly increasing, concave
on $[0,\infty)$, and 
   \[
\Phi(t)
\le\Psi(t)
\stackrel{\eqref{defPsiequiv}}{=}
\Phi(t)+\frac{\Phi(\gamma)}{\gamma}t
\stackrel{\text{concavity of $\Phi$}}{\le}\Phi(t)+\Phi(t)
=2\Phi(t),
\qquad 0<t<\gamma.
\] 

Since for $0<t<\gamma/2$,  $
\Psi(2t)=\Phi(2t)+2\frac{\Phi(\gamma)}{\gamma}t, 
$  it follows    that 
\[
\frac{\Psi(2t)}{\Psi(t)}
=
\frac{\Phi(t)}{\Psi(t)}\cdot 
\frac{\Phi(2t)}{\Phi(t)}
+
\frac{\Phi(\gamma)t}{\gamma\Psi(t)} \cdot  2
\stackrel{\tiny \frac{\Phi(t)}{\Psi(t)}
+\frac{\Phi(\gamma)t}{\gamma\Psi(t)}=1}{\ge}
\min\left\{\frac{\Phi(2t)}{\Phi(t)},2\right\}\stackrel{\tiny\Phi\mbox{ is concave}}{=}\frac{\Phi(2t)}{\Phi(t)}.
\]
Thus,  \[
\liminf_{t\downarrow0}\frac{\Psi(2t)}{\Psi(t)}
\ge
\liminf_{t\downarrow0}\frac{\Phi(2t)}{\Phi(t)}
>1. 
\] 
Moreover,
\[
\frac{\Psi(2t)}{\Psi(t)}
=
\frac{\Phi(\gamma)+\dfrac{\Phi(\gamma)}{\gamma}\,2t}
{\Phi(\gamma)+\dfrac{\Phi(\gamma)}{\gamma}\,t}
=
\frac{\gamma+2t}{\gamma+t}
\rightarrow2, 
\qquad t\to\infty.
\]
\end{proof}

\subsection{Stability properties of Grothendieck spaces}
 We shall use the following stability properties
of Grothendieck spaces.  
\begin{prop}\label{surjective image} \cite[Proposition~3.1.4]{GK}
Let $X$ and $Y$ be Banach spaces.
\begin{enumerate}[label=(\arabic*)]
    \item If $X$ and $Y$ are Grothendieck spaces, then so is
    the direct sum $X\oplus Y$.
    \item If $X$ is a Grothendieck space and there exists a
    bounded surjective linear operator $T:X\to Y$, then
    $Y$ is a Grothendieck space.
\end{enumerate}
\end{prop}

Let $p\in\cM$ be a nonzero projection, and let
$\cA\subseteq p\cM p$ be a von Neumann subalgebra with unit $p$.
Assume that $\tau|_{\cA}$ is semifinite. 
The norm on
$E(\cA,\tau|_{\cA})$ is given by
\begin{align}\label{corneranormMphi}
\norm{a}_{E(\cA,\tau|_{\cA})}
=
\norm{a}_{E(\cM,\tau)},
\qquad a\in E(\cA,\tau|_{\cA}),
\end{align}
where $\mu(a)$ denotes the generalized singular value
function of $a$ relative to $\tau|_{\cA}$. 
Let
\[
\cE_{\cA}:
(L_1+L_\infty)(p\cM p,\tau)
\rightarrow
(L_1+L_\infty)(\cA,\tau|_{\cA})
\]
be the trace-preserving conditional expectation
given by \cite[Proposition~7.1.1]{DPS}. 

The next lemma shows 
the Grothendieck property of a fully symmetric space $E(\cM,\tau)$ (in particular, a Marcinkiewicz space) can be inherited by $E(\cA,\tau|_\cA)$.

\begin{lemma}\label{lem:expectation-quotient} 
Let $E(\cM,\tau)$ be a fully symmetric space. 
The linear contractive  map 
\[
T:E(\cM,\tau)\rightarrow
E(\cA,\tau|_{\cA}),
\qquad
T(x):=\cE_{\cA}(pxp), \quad x\in E\left(\cM,\tau\right)
\]
is    surjective. 
Consequently, if $E(\cM,\tau)$ (in particular, $M_\Phi(\cM,\tau)$) has the
Grothendieck property, then
$E\left(\cA,\tau|_{\cA}\right)$ also has the
Grothendieck property.
\end{lemma}

\begin{proof}
For every $x\in E(\cM,\tau)\subseteq (L_1+L_\infty) (\cM,\tau)$,   the element $pxp\in (L_1+L_\infty)(p\cM p,\tau)$
\cite[p.~241 and Lemma~3.7.1]{DPS}. 
Since 
$\mu(pxp)\le\mu(x)$ and
$\cE_{\cA}(pxp)\prec\prec pxp$
\cite[Proposition~7.1.1(7)]{DPS}, it follows that  
\begin{align*} 
\norm{T(x)}_{E(\cA,\tau|_{\cA})}=
\norm{\cE_{\cA}(pxp)}_{E(\cA,\tau|_{\cA})}  
 \le  
\norm{pxp}_{E(p\cM p,\tau)} 
 \le  
\norm{x}_{E(\cM,\tau)},
\end{align*}
i.e.,  $T$ is contractive. 
Consequently, $T$
is a bounded linear operator from $E(\cM,\tau)$
into $E\left(\cA,\tau|_{\cA}\right)$. 
The surjectivity of $T$ follows from the fact that every
$a\in E\left(\cA,\tau|_{\cA}\right)$ belongs to
$E(\cM,\tau)$ and satisfies
\[
T(a)=\cE_{\cA}(pap)
\stackrel{\text{$\cA\subseteq p\cM p$}}{=}
\cE_{\cA}(a)\stackrel{\text{ \cite[Proposition~7.1.1]{DPS}}}{=}
a.
\]
If $E(\cM,\tau)$ has the
Grothendieck property, then 
by Proposition \ref{surjective image}(2),
$E\left(\cA,\tau|_{\cA}\right)$ has the Grothendieck property. 
\end{proof}

\section{The Grothendieck property for products of noncommutative intersection spaces}\label{S:3}

\subsection{The Grothendieck property of
$\ell_\infty$-sums of uniformly convex spaces}

	Given a family of  Banach spaces $\left(E_k,\nu_k(\cdot)\right)$, $k\in \mathbb{N}$, define the Banach space  
\begin{align}\label{defdirectsumG} 
		G:=\bigoplus_{k\in \mathbb{N}}^{\ell_\infty} E_k.
	\end{align} 
		Let $\mathbf{1}_G:G\to G$ denote the identity operator.
	For $a=\left(a_k\right)_{k\in\mathbb N} \in \ell_\infty$ and $x\in G$, set $ax:=\left(a_kx_k\right)_{k\in\mathbb N} \in G$.  
    For each $A\subseteq\mathbb N$, define the coordinate projection
	$P_A:G\to G$ by
	\begin{align}\label{defPA} P_A(x):=\bigl(x_k\chi_A(k)\bigr)_{k\in \mathbb{N}},
		\qquad x=(x_k)_{k\in \mathbb{N}}
        \in G,
	\end{align} 
	where $\chi_A$ denotes the characteristic function of $A$.

	Let us recall that $\beta\mathbb N$ is the Stone--\v{C}ech
	compactification of $\mathbb N$ endowed with the discrete topology,
	i.e., $\beta\mathbb N$ is the unique compact Hausdorff space
	containing $\mathbb N$ as an open dense subspace such that every bounded
	continuous function on $\mathbb N$ extends to a continuous function
	on $\beta\mathbb N$.
	  Consequently, 
	the restriction map
	\begin{align*}
		R:C(\beta\mathbb N)&\notag \to\ell_\infty, \\ 
		f&\mapsto (f(k))_{k\in\mathbb N},
	\end{align*}
	is an isometric isomorphism \cite[p.~79]{AK06}, where   $C(\beta\mathbb N)$ is  the space of continuous
	complex-valued functions on $\beta\mathbb N$,
	equipped with the supremum norm.
	
	For $x=(x_k)_{k\in\mathbb N}\in G$, define
	\begin{align}\label{defrxknukxk}
	r_x(k):=\nu_k(x_k),\qquad k\in\mathbb N.
	\end{align}
	Then, $r_x\in\ell_\infty$ and $\norm{r_x}_{\ell_\infty}=\norm{x}_G$. 
	Let $\widehat r_x:=R^{-1}(r_x)\in C(\beta\mathbb N)$. Thus,
	\begin{align}\label{propertieswidehatrx}
	\widehat r_x|_{\mathbb N}=r_x,
	\qquad
	\norm{\widehat r_x}_{C(\beta\mathbb N)}=\norm{x}_G.
\end{align}

Let $M(\beta\mathbb N)$ denote the Banach space of finite
regular complex Borel measures on $\beta\mathbb N$, equipped
with the total variation norm
\begin{equation}\label{M(beta N)}
\|\mu\|_{M(\beta\mathbb N)}:=|\mu|(\beta\mathbb N),\quad \mu\in M(\beta\mathbb N).
\end{equation}

    By the Riesz representation theorem  \cite[Theorem~6.19]{R87}, we identify
$M(\beta\mathbb N)$ isometrically with $C(\beta\mathbb N)^*$
via the pairing
\begin{align}\label{defdualityMCbetaN}
\langle\mu,f\rangle
=\int_{\beta\mathbb N}f\,d\mu,
\qquad \mu\in M(\beta\mathbb N),\quad
f\in C(\beta\mathbb N).
\end{align}

Let $\mathcal P(\mathbb N)$ be the Boolean
	algebra  of all subsets of $\mathbb N$.  
    A set function $\rho$ on $\mathcal P(\mathbb N)$ is finitely
additive if $\rho(\varnothing)=0$ and
\[
\rho(A\cup B)=\rho(A)+\rho(B)
\qquad\text{for all disjoint sets}  A,B\subseteq\mathbb N.
\]
It is positive if $\rho(A)\ge0$ for every $A\subseteq\mathbb N$.

	The following lemma follows from the duality between
	$M$-projections and $L$-projections
	\cite[Chapter~I, \S1, pp.~1--2]{HWW93}.
	\begin{lemma}\label{lem:additive-measure}
	Let $G$ be defined by \eqref{defdirectsumG}, and let
	$P_A:G\to G$, $A\subseteq\mathbb N$, be the coordinate
	projections defined by \eqref{defPA}.
		For every $\phi\in G^*$, the set function
		$	\nu_\phi(A):=\norm{\phi\circ P_A}_{G^*} $	
		is a positive, finitely additive measure satisfying
		$\nu_\phi(\mathbb N)=\norm{\phi}_{G^*}$.
	\end{lemma}

		For $A\subseteq\mathbb N$, let
		$\overline A:=\overline A^{\,\beta\mathbb N}\subseteq \beta\mathbb{N}$
		denote its closure in the Stone topology.
	By \cite[Proposition~1.1(a)]{B15},
	$\overline A$ is clopen  in $\beta\mathbb N$, and
		$\chi_{\overline A}$ is the continuous extension of $\chi_A$.
	
	\begin{lemma}\label{lem:dominating-measure}
		Let $G$ and $(P_A)_{A\subseteq\mathbb N}$ be as in
		\eqref{defdirectsumG} and \eqref{defPA}, respectively. 
		For every $\phi\in G^*$, there exists a unique positive 
	 regular Borel measure $\widehat\nu_\phi$
		on $\beta\mathbb N$ such that
		\[
		\widehat\nu_\phi(\overline A)
		=\norm{\phi\circ P_A}_{G^*},
		\qquad A\subseteq\mathbb N.
		\]
		Moreover,
	$
		\widehat\nu_\phi(\beta\mathbb N)=\norm{\phi}_{G^*},
		$ 
		and
		\[
		|\phi(ax)|
		\le\int_{\beta\mathbb N}|a|\widehat r_x\,
		d\widehat\nu_\phi,
		\qquad a\in C(\beta\mathbb N),\quad x\in G,
		\]
		where $|a|$ denotes the pointwise absolute value of $a$.
	\end{lemma}
		\begin{proof} 
	Fix $\phi\in G^*$. By Lemma~\ref{lem:additive-measure},   $\nu_\phi$ is a
	nonnegative finitely additive functional on $\mathcal P(\mathbb N)$.
Identifying $\beta\mathbb N$ with the Stone space of
$\mathcal P(\mathbb N)$ via the canonical homeomorphism
\cite[4A2I(b)(i)]{F03}, we obtain from
\cite[Proposition~416Q(b)]{F03} a unique corresponding
Radon measure on $\beta\mathbb N$.
Denote by $\widehat\nu_\phi$ the restriction of this
Radon measure to the Borel $\sigma$-algebra
$\mathcal B(\beta\mathbb N)$. 
Then,
	\begin{align} \label{widehat phi}
	\widehat\nu_\phi(\overline A)
\overset{\text{ \cite[Proposition~416Q(b)]{F03}}
}{=}	
\nu_\phi(A)
	\overset{
		\text{Lemma~\ref{lem:additive-measure}}
	}{=} \norm{\phi\circ P_A}_{G^*},
	\qquad A\subseteq\mathbb N. 
	\end{align}
In particular,
$
	\widehat\nu_\phi(\beta\mathbb N)
	=\nu_\phi(\mathbb N)
	=\norm{\phi}_{G^*}<\infty.
$
Moreover, by \cite[Proposition 416F, Proposition 412W]{F03}, we have that $\widehat \nu_\phi$ is regular.
		
		Given $x=(x_k)_{k\in \mathbb N}\in G$ and $\varepsilon>0$, choose a finite partition $(A_j)_{j=0}^{m}$ of $\mathbb N$ and nonnegative constants  $\{b_j\}_{j=0}^m$ 
		such that \begin{align}\label{nukxkineq1}
		\nu_k(x_k)
		\overset{\eqref{defrxknukxk}}{=}
		r_x(k)
		\le b_j\le r_x(k)+\varepsilon,
		\qquad k\in A_j,\quad j=0,\dots,m.
		\end{align}
		By \eqref{propertieswidehatrx},  the density of $A_j$ in $\overline A_j$, and continuity of $	\widehat r_x$, we obtain 
		\begin{align}\label{s le r+epsilon}
		\widehat r_x(t)\le b_j\le\widehat r_x(t)+\varepsilon,
		\qquad t\in\overline A_j,\quad j=0,\dots,m.
		\end{align}
		Then,
		\begin{eqnarray*}
			|\phi(x)|
		&	\overset{\eqref{defPA}}{=} & 
			\left|\phi\left(\sum_{j=0}^{m}P_{A_j}(x)\right)\right|
			=\left|\sum_{j=0}^{m}\phi(P_{A_j}(x))\right|
				\\
		&\leq& 
\sum_{j=0}^{m}|\phi(P_{A_j}(x))|
		\leq\sum_{j=0}^{m}
		\norm{\phi\circ P_{A_j}}_{G^*}
		\norm{P_{A_j}(x)}_G		\\
		&\overset{\eqref{widehat phi}}{=} & 
			\sum_{j=0}^m \widehat\nu_\phi(\overline{A_j}) \norm{P_{A_j}(x)}_G
			=
			\sum_{j=0}^m \widehat\nu_\phi(\overline{A_j})\sup_{k\in A_j}\nu_k(x_k)
			\\
		&\overset{\eqref{nukxkineq1}}{	\le} & \sum_{j=0}^{m}b_j\widehat \nu_\phi(\overline{A_j})
	=	\sum_{j=0}^{m} \int_{\overline{A_j}}b_j\,d\widehat\nu_\phi
	\\ &	\stackrel{\eqref{s le r+epsilon}}{\le} & \sum_{j=0}^{m} 
	\int_{\overline{A_j}}
		(\widehat r_x+\varepsilon)\,d\widehat\nu_\phi 
=\int_{\beta \mathbb{N}} \widehat r_x\,d\widehat\nu_\phi+ \varepsilon  \widehat\nu_\phi(\beta\mathbb N)
\\ &	=&\int_{\beta \mathbb{N}} \widehat r_x\,d\widehat\nu_\phi
	+\varepsilon \norm{\phi}_{G^*} .
		\end{eqnarray*}
	Letting $\varepsilon\downarrow0$ yields
		\begin{align} \label{le rx}
			|\phi(x)|\le\int_{\beta \mathbb{N}} \widehat r_x\,d\widehat\nu_\phi.
		\end{align}
		For $a\in C(\beta\mathbb N)$, we have
		\[
		r_{ax}(k)
		\overset{\eqref{defrxknukxk}}{=}
		\nu_k(a(k)x_k)
		=|a(k)|r_x(k),
		\qquad k\in\mathbb N.
		\]
	By continuity of $\widehat r_{ax},|a|\widehat r_x$ and the density of $\mathbb N$ in
	$\beta\mathbb N$, it follows that
	$\widehat r_{ax}=|a|\widehat r_x$.
	Applying \eqref{le rx} to $ax\in G$, we have that 
$
	|\phi(ax)|
	\le \int_{\beta\mathbb N}|a|\widehat r_x\,
	d\widehat\nu_\phi.
$		
	\end{proof}
	
Recall from \cite[p. 347]{AK06} that 	a subset $A$ of a normed space $X$ is said to be {\it relatively weakly compact} if the weak closure of $A$ is compact in the weak topology of $X$.
	\begin{lemma}\label{lem:measure-weak-compactness}
	Let $G$ and $(P_A)_{A\subseteq\mathbb N}$ be as in
	\eqref{defdirectsumG} and \eqref{defPA}, respectively.  	Let $(\phi_n)_{n\in\mathbb N}\subseteq G^*$ converge to zero
		in the weak$^*$ topology.
		For each $n\in\mathbb N$, let $\widehat\nu_{\phi_n}$
		be the positive  regular Borel measure on $\beta\mathbb N$
		given by Lemma~\ref{lem:dominating-measure}.
		Then,
		$
		\{\widehat\nu_{\phi_n}:n\in\mathbb N\}
		$
		is relatively weakly compact in $M(\beta\mathbb N)$.
	\end{lemma}
	
	\begin{proof}
		Assume by contradiction that
		$\{\widehat\nu_{\phi_n}:n\in\mathbb N\}$
		is not relatively weakly compact.
		Since $(\phi_n)_{n\in\mathbb N}$ is weak$^*$ null, it follows from the uniform boundedness principle \cite[D.12]{AK06} that
		\begin{align*}
	\sup_{k\in \mathbb{N}}\norm{\widehat \nu_{\phi_k}}_{M(\beta \mathbb{N})}\stackrel{\eqref{M(beta N)}}{=}\sup_{k\in \mathbb{N}}\widehat \nu_{\phi_k}\left(\beta\mathbb{N}\right)
	\overset{
		\text{Lem.~\ref{lem:dominating-measure}}
	}
	{=}
	\sup_{k\in\mathbb{N}}\norm{\phi_k}_{G^*}\stackrel{\text{
			\cite[D.12]{AK06}
	}}{<}\infty.
\end{align*}
	Since 	$\{\widehat\nu_{\phi_n}:n\in\mathbb N\}$ is norm bounded in $M(\beta\mathbb N)$, it follows from  the Dieudonn\'e--Grothendieck criterion
		\cite[Chapter~VII, Theorem~14]{D84} that there exist 
  $\varepsilon>0$, a strictly increasing sequence
		$(n_j)_{j\in\mathbb N}$, and 
		pairwise disjoint open sets $U_j\subseteq\beta\mathbb N$
such that
\[
\widehat\nu_{\phi_{n_j}}(U_j)>2\varepsilon
\qquad (j\in\mathbb N).
\]

By inner regularity of $\widehat\nu_{\phi_n}$, for each $j\in\mathbb N$, there is a
compact set $K_j$ such that
\begin{align} \label{>2epsilon}
K_j\subseteq U_j,
\qquad
\widehat\nu_{\phi_{n_j}}(K_j)>2\varepsilon.
\end{align}
 For each $x\in K_j$, since $U_j$ is open and
$\beta\mathbb N$ has a clopen basis \cite[4A2I(b)(i) and (iii)]{F03}, there is a clopen set
$V_x$ in $\beta\mathbb N$ such that
\[
x\in V_x\subseteq U_j.
\]
By compactness of $K_j$, there exist $x_1,\ldots,x_{m_j}\in K_j$
such that
\[
K_j\subseteq V_j:=\bigcup_{\ell=1}^{m_j}V_{x_\ell}
\subseteq U_j.
\]
Note that the set $V_j$ is clopen as a finite union of clopen sets. 
Set $A_j:=V_j\cap\mathbb N$. Since $V_j$ is clopen and
$\mathbb N$ is dense in $\beta\mathbb N$, we have
\[
\overline{A_j}
=\overline{V_j\cap\mathbb N}
=V_j.
\]
Moreover, for $i\ne j$,
\[
A_i\cap A_j\subseteq V_i\cap V_j\subseteq U_i\cap U_j
=\varnothing.
\]
Now, we have
a strictly increasing
sequence $(n_j)_{j\in\mathbb N}$, and pairwise disjoint subsets
$A_j\subseteq  \mathbb N$ such that 
\[
\norm{\phi_{n_j}\circ P_{A_j}}_{G^*}
\overset{\text{Lem.~\ref{lem:dominating-measure}}}{=}
\widehat\nu_{\phi_{n_j}}(\overline{A_j})
=\widehat\nu_{\phi_{n_j}}(V_j)
\ge \widehat\nu_{\phi_{n_j}}(K_j)
\stackrel{\eqref{>2epsilon}}{>}2\varepsilon.
\]

		For each $j\in\mathbb N$, by the definition of the dual norm $\norm{\cdot}_{G^*}$, we have that $y_j\in G$ with
	$
		\norm{y_j}_G\le1$  and 
	$
		|\phi_{n_j}(P_{A_j}(y_j))|>\varepsilon.
	$
	Then,  $x_j:=P_{A_j}(y_j)$ satisfies
	\begin{align}\label{PAjxjxj}	P_{A_j}(x_j)=x_j,\qquad
	\norm{x_j}_G\le1,\qquad
	|\phi_{n_j}(x_j)|>\varepsilon.
	\end{align}	
Define $\theta:\ell_\infty\to G$ by
\begin{align}\label{defRbsum}
\theta(b):=\sum_{j=0}^{\infty}b_jP_{A_j}(x_j),
\qquad b=(b_j)_{j\in\mathbb N}\in\ell_\infty,
\end{align}
where the sum is taken coordinatewise. 
Then, 
$\norm{\theta(b)}_G\le\norm{b}_{\ell_\infty}$, i.e., 
 $\theta$ is a bounded linear operator with $\norm{\theta}\le1$.

	For each $j\in\mathbb N$, define a bounded complex-valued finitely additive measure
    \begin{align}\label{mujdefphinjRchiB}
	\mu_j:\mathcal P(\mathbb N)\to\mathbb C,
	\qquad
	\mu_j(B):=\phi_{n_j}(\theta(\chi_B)), \quad B\in \cP(\mathbb N),
		\end{align}
	where  $\chi_B=(\chi_B(k))_{k\in\mathbb N}\in\ell_\infty$ denotes the characteristic function of $B$.
    For any $B\in \cP(\mathbb N)$, by weak$^*$ convergence of  $\left(\phi_{n_j}\right)_{j\in\mathbb N}$, we have that
    $\mu_j(B)\to0$ as $j\to \infty$.
	Applying Phillips' lemma \cite[Lemma~3]{A61} to $(\operatorname{Re}\mu_j)$ and
	$(\operatorname{Im}\mu_j)$ separately, we obtain
	$\mu_j(\{j\})\to0$. 
	Thus, 
	\begin{align*}
	\phi_{n_j}(x_j)
	 \overset{\eqref{PAjxjxj}}{=}
	\phi_{n_j}(P_{A_j}( x_j)
	)
	\overset{\eqref{defRbsum}}{=}
	\phi_{n_j}
	(\theta
	(\chi_{\{j\}})
	)
\overset{\eqref{mujdefphinjRchiB}}{=}
	\mu_j(\{j\})
	\rightarrow
	0, \quad j\to \infty, 
	\end{align*}
	which contradicts $|\phi_{n_j}(x_j)|>\varepsilon$
	for every $j\in\mathbb N$ by \eqref{PAjxjxj}. 
	Therefore, 		$
	\{\widehat\nu_{\phi_n}:n\in\mathbb N\}
	$ is  relatively weakly  compact in $M(\beta\mathbb N)$.
	\end{proof}

Let $X$ be a nonzero Banach space.  
For $0<\varepsilon\le2$, 
the \emph{modulus of convexity} of $X$ is defined by
\begin
{align}\label{defUC}
\delta_X(\varepsilon)
:=\inf\left\{
1-\left\|\frac{x+y}{2}\right\|_X:
\|x\|_X\le1,\ \|y\|_X\le1,\
\|x-y\|_X\ge\varepsilon
\right\}.
\end{align} 
The space $X$ is called \emph{uniformly convex} if
$\delta_X(\varepsilon)>0$ for every $\varepsilon\in(0,2]$. 
A family $(E_k,\nu_k)_{k\in\mathbb N}$ of Banach spaces   is said to admit a
{\it common modulus of uniform convexity} if
\[
\inf_{k\in\mathbb N}\delta_{E_k}(\varepsilon)>0,
\qquad 0<\varepsilon\le2,
\]
where $\delta_{E_k}$ is computed with respect to the norm
$\nu_k$ (see, e.g.,  \cite[Section~2, p.~6]{HMW23},  \cite[pp.~504--505]{D41} and \cite[Section 7.14.2]{DPS}).

\begin{lemma}\label{lem:reflexive-completion}
	Let $G$ be as in \eqref{defdirectsumG}, and suppose that
	$(E_k,\nu_k)_{k\in\mathbb N}$ admits a common modulus
	of uniform convexity.
	Let $\lambda$ be a  regular Borel probability measure 
	on $\beta\mathbb N$, and define
	\begin{align}\label{defPGseminorm}
	p:G \to[0,\infty),\qquad 
	x\mapsto \left(
	\int_{\beta\mathbb N}\widehat r_x^{\,2}\,d\lambda
	\right)^{1/2},
\end{align}
 where $\widehat r_x\in C(\beta\mathbb N)$
	is as in \eqref{propertieswidehatrx}. 
	Then, $p$ is a continuous seminorm on $G$.  
Set $\ker p:=\{x\in G:p(x)=0\}$ and equip the quotient
vector space $G/\ker p$ with the norm induced by $p$, given by
	\begin{align}\label{defnormGkerp}
	\norm{x+\ker p}_{G/\ker p}
  = p(x), 
\end{align}
	and let $(\mathcal H,\norm{\cdot}_\mathcal H)$ be its completion.
	Then, $\mathcal H$ is 
	reflexive. 
\end{lemma}

\begin{proof}
	 It is straightforward to verify that $p$ is a seminorm.
	  Since  $\lambda(\beta\mathbb N)=1$ and for all $t\in\beta\mathbb N$, 
	 $
	 0\le\widehat r_x(t)\le
	 \|\widehat r_x\|_{C(\beta\mathbb N)}
	 \overset{\eqref{propertieswidehatrx}}{=}
	 \|x\|_G,
	 $
	 it follows that 
	 \begin{align}\label{pxxnormG}
	 	p(x)\le\norm{x}_G.
	 \end{align}

	 
	 By the Milman--Pettis theorem \cite[p.~131]{D84},
	 every uniformly convex Banach space is reflexive.
	 Since uniform convexity is preserved under completion
	 \cite[22.42]{S97}, it suffices to prove that
	 $(G/\ker p, \norm{\cdot}_{G/\ker p})$ is uniformly convex.
	 
	 By \cite[Theorem 2.2(3)]{D25},
    there exists
	 a net $(\lambda_i)_{i\in I}\subseteq  M(\beta\mathbb N)$
	 of finitely supported
Borel probability measures on $\beta\mathbb N$, indexed by a directed set $I$,
	 of the form
	 \begin{align}\label{deflambdaiw*}
	 	\lambda_i:=\sum_{k\in F_i}a_{i,k}\delta_k
	 	\rightarrow
	 	\lambda
	 	\quad\text{in the  }\sigma(M(\beta\mathbb N), C(\beta\mathbb N))\text{-topology},
	 \end{align}
	 where $F_i\subseteq\mathbb N$ is finite and nonempty,
	 $a_{i,k}>0$, $\sum_{k\in F_i}a_{i,k}=1$, and
	 $\delta_k$ is the Dirac measure at
	 $k\in\mathbb N$.
	 For each $i\in I$, let 
	 \[
	 H_i:=\left\{
	 z=(z_k)_{k\in F_i}:
	 z_k\in E_k,\ 
	 \norm{z}_{H_i}
	 :=\left(
	 \sum_{k\in F_i}
	 a_{i,k}\cdot
	 \left( \nu_k(z_k)
	 \right)
	 ^2
	 \right)^{1/2}
	 \right\}
	 \]
	 and define
	 \begin{align}\label{defQi}
	 	Q_i:G&\notag \to H_i\\ x=(x_k)_{k\in \mathbb N} &\mapsto  Q_i (x):=(x_k)_{k\in F_i}.
	 \end{align}
	 Observe that, for every $i\in I$ and  $x\in G$, 
	 \begin{align*}
	 \norm{Q_i(x)}_{H_i}^{2}
	  \overset{\eqref{defQi}}{=}
	 \sum_{k\in F_i}a_{i,k}
	 \left( \nu_k(x_k)
	 \right)
	 ^2	
	 \overset{
	 	\eqref{defrxknukxk}, 
	 	~\eqref{propertieswidehatrx}
	 }{=}
	 \sum_{k\in F_i}a_{i,k} \left( 
	 \widehat r_x(k)
	 \right)^2
	 \overset{\eqref{deflambdaiw*}}{=}
	 \int_{\beta\mathbb N}\widehat r_x^{2}\,d\lambda_i
	 \overset{\eqref{defdualityMCbetaN}}{=}
	 \langle \lambda_i,
	 \widehat r_x^{2}
	 \rangle.
	 \end{align*}
	 Therefore, 		
	 \begin{align}\label{Qi(x)normtopx}
	 	\norm{Q_i (x)
	 	}_{H_i}^{\,2}
	 	=	\langle \lambda_i,
	 	\widehat r_x^{2}
	 	\rangle 
	 	\stackrel{\eqref{deflambdaiw*}}{
	 		\rightarrow} 
	 	\langle \lambda,
	 	\widehat r_x^{2}
	 	\rangle 
	 	\overset{\eqref{defdualityMCbetaN}}{=}
	 	\int_{\beta\mathbb N}\widehat r_x^{2}\,d\lambda
	 	\stackrel{\eqref{defPGseminorm}}{=}
	 	p(x)^2.
	 \end{align}

	 Since
	 $(E_k,\nu_k)_{k\in\mathbb N}$ admits a common modulus
	 of uniform convexity, it follows from \cite[Theorem~3]{D41} that 
	 the space
$E:=\mathop{\bigoplus}\limits_{k\in\mathbb N}^{\ell_2}E_k$
	 is uniformly convex.
	 Note that every $H_i$ embeds isometrically into $E$ via
	 $
	 (z_k)_{k\in F_i}
	 \mapsto
	 \bigl(\sqrt{a_{i,k}}\,z_k\bigr)_{k\in F_i} 
	 $ with zero coordinates outside $F_i$. 
  Consequently,
	 \begin{align}\label{deltaHideltaE}
	 	\delta_{H_i}(\varepsilon)\ge\delta_E(\varepsilon)
	 	\overset{\text{E is uniformly convex}}{>}0
	 \end{align}
	 for every $0<\varepsilon\le2$, i.e.,  the family $(H_i)_{i\in I}$ admits a common modulus 	of uniform convexity.
	 
	 Fix $0<\varepsilon\le2$ and $x,y\in G$ such that 
	 $p(x),p(y)\le1$ and $p(x-y)\ge\varepsilon>0$.
	 Set
	 \[
	 c_i:=\max\{1,\norm{Q_i (x)}_{H_i},\norm{Q_i (y)}_{H_i}\}.
	 \]
	 Then, 	\[
	 \left\|\frac{Q_i (x) }{c_i}\right\|_{H_i}\le1,
	 \quad
	 \left\|\frac{Q_i (y)}{c_i}\right\|_{H_i}\le1, \quad \text{and }\,  c_i
	 \overset{\eqref{Qi(x)normtopx}}{\to } 
	 1 \text{ as } i\to \infty, 
	 \]   and so  $\frac{1}{c_i}	\left\|
	 Q_i (x-y)
	 \right\|_{H_i}\overset{\eqref{Qi(x)normtopx}}{\to}  p(x-y)$. 
	 This implies that for all sufficiently large $i$,
	 \[
	 \left\| \frac{Q_i (x)}{c_i} -
	 \frac{Q_i (y)}{c_i}
	 \right\|_{H_i}
	 \ge p(x-y)-\frac{\varepsilon}{2} \ge \varepsilon-\frac{\varepsilon}{2}
	 =\frac{\varepsilon}{2}.
	 \]
By the definition of the modulus of convexity and
\eqref{deltaHideltaE}, we have that  for all sufficiently large $i$,  
	 \begin{eqnarray*}
	 	\left\|Q_i\left(\frac{x+y}{2}\right)\right\|_{H_i}
	 	&=&c_i\left\|
	 	\frac12\left(\frac{Q_i (x) }{c_i}+\frac{Q_i (y)}{c_i}\right)
	 	\right\|_{H_i}
	 	\\ & \stackrel{\eqref{defUC}}{\le}&
	 	c_i\bigl(1-\delta_{H_i}(\varepsilon/2)\bigr)
	 	\stackrel{\eqref{deltaHideltaE}}{\leq } c_i\bigl(1-\delta_E(\varepsilon/2)\bigr).
	 \end{eqnarray*}
	 Passing to the limit,  using \eqref{Qi(x)normtopx} and $c_i\to1$,  we obtain
	 \[
	 p\left(\frac{x+y}{2}\right)
	 \le 1-\delta_E (\varepsilon/2).
	 \]
	 Together with the definition of the norm $\norm{\cdot}_{G/\ker p}$ (see  \eqref{defnormGkerp}), we conclude that  $G/\ker p$ 
	 is uniformly convex. 
     Therefore, 
 	 $(\mathcal H, \norm{\cdot}_\mathcal H)$ is reflexive. 
\end{proof}

	Let $(\Omega,\Sigma,\lambda)$ be a finite positive measure
space. A bounded family $\mathcal F\subseteq L_1(\Omega,\lambda)$
is called \emph{uniformly integrable} if
\[
\lim_{R\to\infty}\sup_{f\in\mathcal F}
\int_\Omega (|f|-R)_+\,d\lambda=0,
\]
where $t_+:=\max\{t,0\}$ for $t\in\mathbb R$ \cite[Proposition~2.1]{K14}, see also e.g.,   \cite[Chapter~III, \S2, Definition~10]{DU}.
By the Dunford--Pettis theorem, a bounded family
$\mathcal F\subseteq L_1(\Omega,\lambda)$ is uniformly
integrable  if and only if 
$\mathcal F$ is relatively
weakly compact in $L_1(\Omega,\lambda)$
\cite[Chapter~III, \S2, Theorem~15]{DU}.  	
Let 
$\nu$  be another  finite positive measure on
$(\Omega,\Sigma)$. 
We write $\nu\ll\lambda$ if $\nu$ is absolutely continuous
with respect to $\lambda$, that is,
\[
\lambda(A)=0\implies\nu(A)=0,
\qquad A\in\Sigma.
\]	
If $\nu\ll\lambda$, we denote the Radon--Nikod\'{y}m derivative
of $\nu$ with respect to $\lambda$ by $d\nu/d\lambda\in L_1(\Omega,\lambda)$, satisfying
\begin{align}\label{defnuAderivative}
\nu(A)=\int_A \frac{d\nu}{d\lambda}\,d\lambda,
\qquad A\in\Sigma; 
\end{align}
 see \cite[p.~59]{DU}.

Recall Grothendieck's criterion
\cite[Theorem~3.44]{AOB06}:
 a subset $D$ of a Banach space
$Y$ is relatively weakly compact if and only if, for every
$\varepsilon>0$, there exists a weakly compact subset
$W\subseteq Y$ such that
\[
D\subseteq W+\varepsilon B_Y,
\]
where
$B_{Y}:=\{x\in Y:\|x\|_{Y}\le1\}$
is the closed unit ball of $Y$.

	\begin{lemma}\label{lem:controlled-families}
			Let $G$  be as in
		\eqref{defdirectsumG}, and suppose that
		$(E_k,\nu_k)_{k\in\mathbb N}$ admits a common modulus
		of uniform convexity.
	Let $\lambda$ be a
	regular Borel probability measure on $\beta\mathbb N$, and let
		$\mathcal F\subseteq L_1(\beta\mathbb N,\lambda)_+$
		be bounded and uniformly integrable. 
        Then,  the family of all $\eta\in G^*$ for which there exists
an $f\in\mathcal F$ satisfying
\begin{equation}\label{the family}
  |\eta(x)|\le \int_{\beta\mathbb N} f\widehat r_x\,d\lambda
  \qquad \text{for every } x\in G,
\end{equation}
is relatively weakly compact in $G^*$,  
	 where $\widehat r_x\in C(\beta\mathbb N)$
		is as in \eqref{propertieswidehatrx}. 
	\end{lemma}
	
	\begin{proof}

Let $p$ and $\mathcal H$ be constructed from $\lambda$ as in
Lemma~\ref{lem:reflexive-completion}. 	
Fix $\eta\in G^*$ and $f\in\mathcal F$ such that
\begin{align}\label{Tetaxinterleq}
	|\eta(x)|\le\int_{\beta\mathbb N}  f \widehat r_x \,d\lambda 
\end{align} for every $x\in G$. 
Define the linear operator 
\[
T_\eta:C(\beta\mathbb N)\to G^*,
\qquad
(T_\eta (g) )(x):=\eta(gx),
\quad g\in C(\beta\mathbb N)=\ell_\infty,\ x\in G.
\]
Since $\widehat r_{gx}=|g|\widehat r_x$ (see the proof of Lemma \ref{lem:dominating-measure}), we have
\begin{eqnarray} \label{exten before}
|(T_\eta (g))(x)|
 &=&|\eta(gx)| 
\stackrel{\eqref{Tetaxinterleq}}{\le}
\int_{\beta\mathbb N}f\widehat r_{gx}\,d\lambda
  =
\int_{\beta\mathbb N}|g|f\widehat r_x\,d\lambda\\
&\stackrel{\widehat r_x\le \norm{x}_G}{\le}&
\|x\|_G\int_{\beta\mathbb N}|g|f\,d\lambda \notag.
\end{eqnarray}
Consequently,
\begin{align} \label{contractive}
\|T_\eta (g)\|_{G^*}
\le\int_{\beta\mathbb N}|g|f\,d\lambda  
= 
\int_{\beta \mathbb{N}}|g|d(f\lambda)
=
\norm{g}_{L_1\left(\beta\mathbb{N}, f\lambda\right)}
<\infty,
\end{align}
where $f\lambda$ is the finite positive regular Borel measure
defined by
\begin{align}\label{defflambdaaction}
(f\lambda)(A):=\int_A f\,d\lambda,
\qquad A\in\mathcal B(\beta\mathbb N).
\end{align}
By the density of $C(\beta\mathbb N)$ in
$L_1(\beta\mathbb N,f\lambda)$ \cite[Theorem~3.14]{R87}, $T_\eta$ extends uniquely
to a linear contraction
\[
T_\eta:L_1(\beta\mathbb N,f\lambda)\to G^*,
\]
where the same notation is used for this extension.
By continuity of $T_\eta$ and \eqref{exten before}, we have
\begin{align}\label{Tetagxineq}
|(T_\eta (g))(x)|
\le\int_{\beta\mathbb N}|g|f\widehat r_x\,d\lambda,
\qquad
g\in L_1(\beta\mathbb N,f\lambda),\ x\in G.
\end{align}

Let $\bf 1$ denote the constant
function with value $1$ on $\beta\mathbb N$.
For $R>0$, define the function  $h_R:\beta\mathbb N\to[0,1]$ by
\[
h_R(t):=
\begin{cases}
	\min\{1,R/f(t)\},& f(t)>0,\\
	0,& f(t)=0.
\end{cases}
\]
Then, 
\begin{align}\label{defhRfmulti}
h_R\cdot f=\min\{f,R\mathbf 1\},\qquad
(\mathbf 1-h_R) \cdot f=(f-R\mathbf 1)_+.
\end{align}
Since $f\in L_1(\beta\mathbb N,\lambda)_+$, we have
\[
\begin{aligned}
	\int_{\beta\mathbb N}|h_R|\,d(f\lambda)
	=\int_{\beta\mathbb N}h_R\cdot f\,d\lambda
	=\int_{\beta\mathbb N}\min\{f,R\mathbf 1\}\,d\lambda
	\le\int_{\beta\mathbb N}f\,d\lambda<\infty.
\end{aligned}
\]
Thus, $h_R\in L_1(\beta\mathbb N,f\lambda)_+$.

Using $T_\eta(\bf 1)=\eta$ and  $\norm{T_\eta}\stackrel{\eqref{contractive}}{\leq} 1$,
we obtain
\begin{align} \label{uniformly integr}
	\norm{\eta-T_\eta(h_R)}_{G^*}
	&
	=\norm{T_\eta(
		{\bf 1} -h_R)
		}_{G^*}
	\stackrel{\eqref{contractive}}{\le}\norm{	{\bf 1} -h_R}_{L_1(\beta\mathbb N,f\lambda)}\notag\\
	&
	=\int_{\beta\mathbb N}(	{\bf 1} -h_R)f\,d\lambda
\overset{\eqref{defhRfmulti}}{=}
	\int_{\beta\mathbb N}(f-R\mathbf 1)_+\,d\lambda.
\end{align}
For every $x\in G$, 
\begin{eqnarray*}
	|(T_\eta(h_R))(x)|
	&\stackrel{\eqref{Tetagxineq}}{\le}&
	\int_{\beta\mathbb N}h_Rf\widehat r_x\,d\lambda \\ 
	& \stackrel{\eqref{defhRfmulti}}{=} & 
	\int_{\beta\mathbb N}\min\{f,R\mathbf 1\}\widehat r_x\,d\lambda
	\\ 
	&  \le &
	R\int_{\beta\mathbb N}\widehat r_x\,d\lambda\\
	&\stackrel{\text{Cauchy--Schwarz ineq.}}{\le}&
	R\left(
	\int_{\beta\mathbb N}\widehat r_x^{\,2}\,d\lambda
	\right)^{1/2}
	\left(
	\int_{\beta\mathbb N}\mathbf1^2\,d\lambda
	\right)^{1/2}\\
	&\stackrel{\lambda(\beta\mathbb N)=1}{=}&
	R\left(
	\int_{\beta\mathbb N}\widehat r_x^{\,2}\,d\lambda
	\right)^{1/2}
	\stackrel{\eqref{defPGseminorm}}{=}
	Rp(x)
	\\ 
	&  \stackrel{\eqref{defnormGkerp}}{=}	 & 
	R\norm{x+\ker p}_{G/\ker p}.
\end{eqnarray*}
Consequently, $T_\eta(h_R)$ vanishes on $\ker p$,
and the map
\begin{align}\label{xidefker}
\xi:	G/\ker p  \to \mathbb{C}, \quad  \xi ( x+\ker p):= (T_\eta(h_R))(x)
\end{align}
is a well-defined bounded linear functional  of norm at most $R$.
By density of $G/\ker p$ in $\mathcal H$, the functional on
$G/\ker p$ extends uniquely to a bounded linear
functional on $\mathcal H$ of norm at most $R$. Denote this extension again by $\xi$.

Define
the canonical map
\begin{align}\label{defJGH}
J:G \to \mathcal H,\qquad J(x):=x+\ker p.
\end{align}
Then,  $\norm{J(x)}_\mathcal H\stackrel{\eqref{defnormGkerp}}{=}p(x)
\overset{\eqref{pxxnormG}}{\leq}
\norm{x}_G$. 
Its adjoint is given by 
 \begin{align}\label{defJ*H*G*}
J^*:\mathcal H^*\to G^*,\qquad J^*(\theta):=\theta\circ J, \quad \theta\in \mathcal H^*.
\end{align}
In particular, \(\|J^*\|\le1\).  
Since  for  all 
$ x\in G,
$
$$
(T_\eta(h_R))(x)
\stackrel{\eqref{xidefker}}
{=} 
\xi ( x+\ker p)
\stackrel{\eqref{defJGH}}
{=}
\xi(J(x))
\stackrel{\eqref{defJ*H*G*}}
{=}
(J^*(\xi))(x),$$ 
 it follows that 
\begin{align} \label{Teta}
T_\eta(h_R)
\in J^*\bigl(\{\zeta\in \mathcal H^*:
\norm{\zeta}_{\mathcal{H}^*}\le R\}\bigr).
\end{align}
Since $\mathcal H$ is reflexive (see Lemma \ref{lem:reflexive-completion}), it follows from  \cite[G.2 and p.348]{AK06}
that  the closed ball
$\{\zeta\in \mathcal H^*:\norm{\zeta}_{\mathcal H^*}\le R\}$ is weakly compact.
Since $J^*$ is bounded and linear, it follows that 
$
J^*\bigl(\{\zeta\in \mathcal H^*:\norm{\zeta}_{\mathcal H^*}\le R\}\bigr)
$
is weakly compact in $G^*$.

Let $\varepsilon>0$. 
By uniform integrability of
$\mathcal F$, choose $R>0$ such that
\begin{align} \label{eta <}
\norm{\eta-T_\eta(h_R)}_{G^*} \stackrel{\eqref{uniformly integr}}{\le} \sup_{f\in\mathcal F}
\int_{\beta\mathbb N}(f-R\mathbf 1)_+\,d\lambda<\varepsilon.
\end{align}
Consider the following decomposition
\[
\eta=T_\eta(h_R)+\bigl(\eta-T_\eta(h_R)\bigr),
\]
where $
J^*\bigl(\{\zeta\in \mathcal H^*:\norm{\zeta}_{\mathcal H^*}\le R\}\bigr)
$
is weakly compact in $G^*$ and
\[
T_\eta(h_R)\stackrel{\eqref{Teta}}{\in}
J^*\bigl(\{\zeta\in \mathcal H^*:\|\zeta\|_{\mathcal H^*}\le R\}\bigr),
\qquad
\|\eta-T_\eta(h_R)\|_{G^*}\stackrel{\eqref{eta <}}{<}\varepsilon.
\]
Consequently,  the family defined by \eqref{the family} is contained in
\[
J^*\bigl(\{\zeta\in \mathcal H^*:\norm{\zeta}_{\mathcal H^*}\le R\}\bigr)
+\varepsilon B_{G^*}.
\]
By Grothendieck's criterion  \cite[Theorem~3.44]{AOB06},  
the family defined by \eqref{the family} is relatively weakly compact in $G^*$.   
	\end{proof}

The following corollary is a consequence of 
 \cite[Theorem~1.1]{LLWZ26}, which answers \cite[Conjecture~1]{K99} in the affirmative, i.e., $\bigoplus^{\ell_\infty}_{k\in \mathbb{N}}(E)$ has the Grothendieck property whenever $E$ is a super-reflexive Banach space (recall that every super-reflexive space admits an equivalent uniformly convex norm \cite[Theorem~9.18]{FHH01}). 
We include an alternative proof below.

	\begin{corollary}\label{cor:uc-product}
		Let $(E_k,\nu_k)_{k\in\mathbb N}$ be a family of Banach
		spaces admitting a common modulus of uniform convexity.
	Then, $G:=\bigoplus_{k\in \mathbb{N}}^{\ell_\infty} E_k$ 
    has the Grothendieck property. 
	\end{corollary}
	
\begin{proof}
	Let $(\phi_n)_{n\in\mathbb N}$ be a weak$^*$ null sequence
	in  $	G^* $.  
    For each $n\in\mathbb N$, let $\widehat\nu_{\phi_n}$ be
	the unique  regular Borel measure on
	$\beta\mathbb N$ given by
	Lemma~\ref{lem:dominating-measure}. 
    Since $(\phi_n)_{n\in\mathbb N}$ is weak$^*$ null, it follows from the uniform boundedness
principle \cite[D.12]{AK06} that 
	\begin{align}\label{boundedfnL1lambda}
	\sup_{n\in\mathbb N}\norm{\phi_n}_{G^*}<\infty.
	\end{align}
Define a finite positive regular Borel measure $\lambda_0$
on $\beta\mathbb N$ by
	\[
	\lambda_0:=
	\sum_{n\in\mathbb N}
	\frac{2^{-n}\widehat\nu_{\phi_n}}{1+\norm{\phi_n}_{G^*}}.
	\]
	If $\lambda_0=0$, then
	$\norm{\phi_n}_{G^*}
	\overset{\text{	Lemma~\ref{lem:dominating-measure}}}{=}
	\widehat\nu_{\phi_n}(\beta\mathbb N)=0$
	for every $n\in\mathbb N$, and hence  $(\phi_n)_{n\in\mathbb N}$ converges weakly to zero.
	If $\lambda_0\neq 0$, define
	\[
	\lambda:=
	\frac{\lambda_0}{\lambda_0(\beta\mathbb N)}.
	\] 
	Then, $\lambda$ is a regular Borel probability measure on
$\beta\mathbb N$.
    Since $\widehat\nu_{\phi_n}\ll\lambda$, we may define
	\begin{equation}\label{deffnderivativelambda}
		f_n:=\frac{d\widehat\nu_{\phi_n}}{d\lambda}
		\in L_1(\beta\mathbb N,\lambda)_+,
		\qquad n\in\mathbb N.
	\end{equation} 
    Then, $\widehat\nu_{\phi_n}\stackrel{\eqref{defnuAderivative}, \eqref{defflambdaaction}}{=}f_n\lambda$ and
	\[
	\int_{\beta\mathbb N}f_n\,d\lambda
	\stackrel{\eqref{defnuAderivative}}{=}\widehat\nu_{\phi_n}(\beta\mathbb N)
 \overset{\text{	Lemma~\ref{lem:dominating-measure}}}{=} \norm{\phi_n}_{G^*}<\infty.
	\]
	Consequently,
\begin{equation}\label{boundedfnL1lambdafn}
\sup_{n\in\mathbb N}\norm{f_n}_{L_1(\beta\mathbb N,\lambda)}
=\sup_{n\in\mathbb N}\norm{\phi_n}_{G^*}
\overset{\eqref{boundedfnL1lambda}}{<} \infty.
\end{equation}

	By Lemma~\ref{lem:measure-weak-compactness},
	the family $\{  \widehat\nu_{\phi_n}: n\in\mathbb N\} $
	is relatively weakly compact in $M(\beta\mathbb N)$. 
For every $f\in L_1(\beta\mathbb N,\lambda)$, let $f\lambda$
denote the finite regular complex Borel measure defined by
\eqref{defflambdaaction}. The map $f\mapsto f\lambda$ is an
injective linear isometry from $L_1(\beta\mathbb N,\lambda)$
into $M(\beta\mathbb N)$, since 
	\[
	\norm{f\lambda}_{M(\beta\mathbb N)}
\stackrel{\eqref{M(beta N)},  \eqref{defflambdaaction}}
{=}	
\int_{\beta\mathbb N}|f|\,d\lambda
	=\norm{f}_{L_1(\beta\mathbb N,\lambda)}.
	\]
Since  $ \{  \widehat\nu_{\phi_n}: n\in\mathbb N\}  \stackrel{\eqref{deffnderivativelambda}}{=} 
\{ f_n\lambda: n\in\mathbb N\}$ is relatively
weakly compact in $M(\beta\mathbb N)$, the weak-to-weak
continuity of the inverse of $f\mapsto f\lambda$ on
its closed range implies that
$\{f_n:n\in\mathbb N\}$ is relatively weakly compact
in $L_1(\beta\mathbb N,\lambda)$ \cite[Theorem 3.10]{B11}. 
	The Dunford--Pettis criterion \cite[Chapter~III, \S2, Theorem~15]{DU}  therefore
	implies that $(f_n)_{n\in\mathbb N}$ is uniformly
	integrable.

By Lemma~\ref{lem:dominating-measure} with $a={\bf 1}\in  C(\beta\mathbb N)$, we obtain 
\[
|\phi_n(x)|
\le\int_{\beta\mathbb N}\widehat r_x\,
d\widehat\nu_{\phi_n}
\stackrel{\eqref{deffnderivativelambda}}{=}
\int_{\beta\mathbb N}f_n\widehat r_x\,d\lambda, \quad   
x\in G, 
\quad n\in\mathbb N.
\]
Since $\{ f_n: n\in\mathbb N\} \subseteq L_1(\beta\mathbb N,\lambda)_+$ is bounded (see \eqref{boundedfnL1lambdafn})  and   uniformly
integrable,
 it follows from   Lemma~\ref{lem:controlled-families}  that
	$\{  \phi_n: n\in\mathbb N\} $ is relatively weakly compact 	in $	G^* $.
  Let $(\phi_{n_k})_{k\in\mathbb N}$ be an arbitrary subsequence
of $(\phi_n)_{n\in\mathbb N}$. By relative weak compactness of $\phi_n$ and
the Eberlein--\v{S}mulian theorem
\cite[Theorem~13.1, Chapter~V]{C19},
$(\phi_{n_k})_{k\in\mathbb N}$ has a further subsequence
converging weakly in $G^*$. Since weak convergence implies
weak$^*$ convergence and $\phi_n\to0$ in $\sigma\left(G^*,G\right)$, the weak
limit of this further subsequence is $0$.
Thus,  every subsequence of $(\phi_n)_{n\in\mathbb N}$ has a
further subsequence converging weakly to $0$. It follows
from \cite[Theorem~2.3.3]{D19} that $(\phi_n)_{n\in\mathbb N}$
is weakly null. Therefore, $G$ has the Grothendieck property. 
\end{proof}

\subsection{The Grothendieck property
of $\ell_\infty$-sums of subspaces}

We first record an auxiliary result on weak$^*$
convergent subsequences, which is certainly known to experts. 

\begin{lemma}\label{lem:weak-star-cluster}
Let $G$ be a Banach space, and let
$(\Psi_n)_{n\in\mathbb N}$ be a norm-bounded sequence
in $G^*$. Set
\[
K:=
\bigcap_{N\in\mathbb N}
\overline{\{\Psi_n:n\ge N\}}^{\,\sigma(G^*,G)}.
\]
If $K$ is weakly compact in $G^*$, then
$(\Psi_n)_{n\in\mathbb N}$ has a subsequence $(\Psi_{n_j})_{j\in \mathbb N}$ 
converging in $\sigma(G^*,G)$ to an element of $K$.
\end{lemma}
\begin{proof}

We claim that
\begin{align*}
\overline{\{\Psi_n:n\in\mathbb N\}}^{\,\sigma(G^*,G)}
=K\cup\{\Psi_n:n\in\mathbb N\}.
\end{align*}
In particular,   $\overline{\{\Psi_n:n\in\mathbb N\}}^{\,\sigma(G^*,G)}
\setminus K$ is at most countable. 
Indeed, let
$
\Psi\in
\overline{\{\Psi_n:n\in\mathbb N\}}^{\,\sigma(G^*,G)}$ 
and $\Psi\notin\{\Psi_n:n\in\mathbb N\}.
$
For every $N\in\mathbb N$, the finite set
$\{\Psi_n:n<N\}$ is weak$^*$ closed and hence,
\[
\begin{aligned}
	\overline{\{\Psi_n:n\in\mathbb N\}}^{\,\sigma(G^*,G)}
	=\{\Psi_n:n<N\}
	\cup
	\overline{\{\Psi_n:n\ge N\}}^{\,\sigma(G^*,G)}.
\end{aligned}
\]
Since $\Psi\notin\{\Psi_n:n\in\mathbb N\}$, it follows that
$
\Psi\in
\overline{\{\Psi_n:n\ge N\}}^{\,\sigma(G^*,G)}$ 
for every $N\in\mathbb N$. 
Thus,
\[
\Psi\in
\bigcap_{N\in\mathbb N}
\overline{\{\Psi_n:n\ge N\}}^{\,\sigma(G^*,G)}
=K.
\]
This proves one inclusion, while the reverse inclusion
follows directly from the definition of $K$,
which proves the claim.

Moreover, every point of
$\overline{\{\Psi_n:n\in\mathbb N\}}^{\,\sigma(G^*,G)}
\setminus K$
is isolated in the relative weak$^*$ topology.
Indeed, let $\Psi\in \overline{\{\Psi_n:n\in\mathbb N\}}^{\,\sigma(G^*,G)}
\setminus K$.
By the definition of $K$, there exists $N\in\mathbb N$
such that
$\Psi\notin
\overline{\{\Psi_n:n\ge N\}}^{\,\sigma(G^*,G)}$. 
Note that the set
$G^*\setminus
\left[  
\overline{\{\Psi_n:n\ge N\}}^{\,\sigma(G^*,G)}
\cup(\{\Psi_n:n<N\}\setminus\{\Psi\})
\right]$
is a weak$^*$ open neighbourhood of $\Psi$.
Since
\[
\overline{\{\Psi_n:n\in\mathbb N\}}^{\,\sigma(G^*,G)}
=
\{\Psi_n:n<N\}
\cup
\overline{\{\Psi_n:n\ge N\}}^{\,\sigma(G^*,G)},
\] 
it follows that   $\Psi$ is the only point of the intersection of
$\overline{\{\Psi_n:n\in\mathbb N\}}^{\,\sigma(G^*,G)}$
and
$G^*\setminus\left[
\overline{\{\Psi_n:n\ge N\}}^{\,\sigma(G^*,G)}
\cup(\{\Psi_n:n<N\}\setminus\{\Psi\})
\right]$.
Thus, $\Psi$ is isolated in	$\overline{\{\Psi_n:n\in\mathbb N\}}^{\,\sigma(G^*,G)}$  in the relative weak$^*$ topology.

Since $\Psi_n$ are norm-bounded, it follows from \cite[F.5]{AK06} that
the set
$\overline{\{\Psi_n:n\in\mathbb N\}}^{\,\sigma(G^*,G)}$
is weak$^*$ compact in $G^*$.
Then, $K$ is a weak$^*$ closed subset of $\overline{\{\Psi_n:n\in\mathbb N\}}^{\,\sigma(G^*,G)}$. 
 Together with the fact that  every point of $\overline{\{\Psi_n:n\in\mathbb N\}}^{\,\sigma(G^*,G)}
  \setminus K$ 
  is isolated in the relative weak$^*$ topology, 
  it follows from \cite[Definition~2.1]{MR} 
  that 
$\overline{\{\Psi_n:n\in\mathbb N\}}^{\,\sigma(G^*,G)}$
is a countably discrete extension of $K$.
Since $K$ is weakly compact in $G^*$, it follows from \cite[p.~487]{MiR77} that $K$  is an
Eberlein compactum when equipped with the topology
induced by $\sigma(G^*,G^{**})$. 
By \cite[p.~494]{MiR77}, we have that every Eberlein compactum is Corson compact.

  Note that  the weak and weak$^*$ topologies coincide on $K$ (since $K$ is weakly compact in $G^*$).
Therefore, every nonempty closed subset of $K$,
with closedness taken in the weak$^*$ topology on $K$,
is Corson compact and has a $G_\delta$-point
\cite[Theorem~14.41]{FHH}. 
By \cite[Lemma~2.2]{MR},
$\overline{\{\Psi_n:n\in\mathbb N\}}^{\,\sigma(G^*,G)}$
is sequentially compact.  
Hence, there exist a strictly increasing sequence
$(n_j)_{j\in\mathbb N}$ in $\mathbb N$ and $\Psi\in G^*$
such that
\[
\Psi_{n_j}\rightarrow\Psi
\quad\text{in }\sigma(G^*,G).
\]
For each fixed  $N\in\mathbb N$, we have $n_j\ge N$
for all sufficiently large $j$, and therefore,
$\Psi_{n_j}\in\{\Psi_n:n\ge N\}$.
Hence, we have 
\[
\Psi\in
\overline{\{\Psi_n:n\ge N\}}^{\,\sigma(G^*,G)}.
\]
Thus,  
$$\Psi\in \bigcap_{N\in\mathbb N}
\overline{\{\Psi_n:n\ge N\}}^{\,\sigma(G^*,G)}=K.$$
The proof is complete. 
\end{proof}

	Let $(E_k,\nu_k)_{k\in\mathbb N}$  and  $G$   be as in
	\eqref{defdirectsumG}, and let $Y_k$ be a  subspace of
$E_k$ for each $k\in\mathbb N$. Set
\[
Y:=\mathop{\bigoplus}\limits_{k\in\mathbb N}^{\ell_\infty}Y_k,
\]
equipped with the norm inherited from $G$. 
 For each $k\in\mathbb N$, let $Q_k:=E_k/Y_k$ equipped   with the  quotient norm
\begin{align}\label{qunormEkYk}
\norm{x+Y_k}_{E_k/Y_k}
\stackrel{\eqref{defquotientnorm}}{:=}\inf_{y\in Y_k}\nu_k(x-y),
\qquad x\in E_k.
\end{align}

\begin{lemma}\label{lem:controlled-extensions} 
Suppose that the quotients $Q_k$ admit norms $\norm{\cdot}_{Q_k}$
uniformly equivalent to their quotient norms
and having a common modulus of uniform convexity.
Then, every weak$^*$ null sequence
$(\phi_n)_{n\in\mathbb N}$ in $Y^*$ admits a
norm-bounded sequence of extensions
$(\Psi_n)_{n\in\mathbb N}$ in $G^*$   with
$\Psi_n|_Y=\phi_n$ for every $n\in\mathbb N$, such that
\[
K:=
\bigcap_{N\in\mathbb N}
\overline{\{\Psi_n:n\ge N\}}^{\,\sigma(G^*,G)}
\]
is weakly compact in $G^*$.
\end{lemma}
\begin{proof}

Since $\norm{\cdot}_{Q_k}$ is equivalent to
$\norm{\cdot}_{E_k/Y_k}$ with equivalence constants
independent of $k$, it follows that there exists a constant $C\ge1$ 
such that
\begin{align}\label{uniformequiv}
C^{-1}\norm{u}_{E_k/Y_k}
\le\norm{u}_{Q_k}
\le C\norm{u}_{E_k/Y_k},
\qquad u\in Q_k,\quad \forall k\in\mathbb N.
\end{align}
	Set
	\begin{align*}
	Q:=\bigoplus_{k\in\mathbb N}^{\ell_\infty}Q_k,
	\qquad
	\norm{z}_Q:=\sup_{k\in\mathbb N}\norm{z_k}_{Q_k},
	\quad z=(z_k)_{k\in\mathbb N}\in Q.
	\end{align*}
	Define the linear operator
	\begin{align*}
	q:G &\notag  \to Q,
	\\ 
x=(x_k)_{k\in\mathbb N} &\mapsto 	q(x):=(x_k+Y_k)_{k\in\mathbb N}.
	\end{align*}
Note that   $\ker q=Y$ and $	\norm{q}\leq C$, since for all $x=(x_k)_{k\in \mathbb N}\in G$, 
	\begin{eqnarray*}
	\norm{q(x)}_Q
	&=&\sup_{k\in\mathbb N}\norm{x_k+Y_k}_{Q_k} 
	\stackrel{\eqref{uniformequiv}}{\le}
	 C\sup_{k\in\mathbb N}\|x_k+Y_k\|_{E_k/Y_k}
	\\ & \stackrel{\eqref{qunormEkYk}}{\le}& C\sup_{k\in\mathbb N}\nu_k(x_k)
	 \overset{	\eqref{defdirectsumG}}{=}	C\|x\|_G.
	\end{eqnarray*}
	For every $z=(z_k)_{k\in\mathbb N}\in Q$ and
	$\varepsilon>0$, choose representatives $x_k\in E_k$
	such that
	\begin{align}\label{decompzkxy}
	x_k+Y_k=z_k,
	\quad
	 \nu_k(x_k)
\stackrel{\eqref{qunormEkYk}}{\le}
\norm{z_k}_{E_k/Y_k}+\varepsilon
		\stackrel{\eqref{uniformequiv}}{\le} C\norm{z_k}_{Q_k}+\varepsilon,
	\quad k\in\mathbb N.
	\end{align}
	Then, 
	\begin{align}\label{qx=z}
	x=(x_k)_{k\in\mathbb N}\in G\quad \text{and} \quad   q(x)=z.
	\end{align}
	Thus, $q$ is surjective.
	In particular, \begin{align}\label{1.4}
	r_x(k)
	\stackrel{\eqref{defrxknukxk}}{=}
	\nu_k(x_k) 	\stackrel{\eqref{decompzkxy}}{\le} C\norm{z_k}_{Q_k}+\varepsilon,
		\quad k\in\mathbb N.
	\end{align}

Let $(\phi_n)_{n\in\mathbb N}\subseteq Y^*$ be
weak$^*$ null, and assume that not all $\phi_n$ are zero.
For each $n\in\mathbb N$, let $\widehat\nu_{\phi_n}$
be the measure given by
Lemma~\ref{lem:dominating-measure}, applied to $Y$. 
Following the measure construction and compactness
argument in the proof of Corollary~\ref{cor:uc-product},
we obtain a  regular Borel probability measure
\[
\lambda:=
\frac{\displaystyle
	\sum_{n\in\mathbb N}2^{-n}\widehat\nu_{\phi_n}}
{\displaystyle
	\sum_{n\in\mathbb N}2^{-n}\norm{\phi_n}_{Y^*}}
\]
such that $\widehat\nu_{\phi_n}\ll\lambda$ for every
$n\in\mathbb N$ and the family
 	$\left\{ 
 C\frac{d\widehat\nu_{\phi_n}}{d\lambda}: n\in\mathbb N 
 \right\}$
 is relatively weakly compact
in
$L_1(\beta\mathbb N,\lambda)_+$. By the Dunford--Pettis criterion \cite[Chapter~III, \S2, Theorem~15]{DU},  the weakly compact closure of $\left\{ 
C\frac{d\widehat\nu_{\phi_n}}{d\lambda}: n\in\mathbb N 
\right\}$ is
bounded and uniformly integrable.

For each $n\in\mathbb N$, the map
\[
x\mapsto
\int_{\beta\mathbb N}\widehat r_x\,
d\widehat\nu_{\phi_n},
\qquad x\in G,
\]
is a continuous seminorm on $G$.
By Lemma~\ref{lem:dominating-measure}  with $a={\bf 1}\in  C(\beta\mathbb N)$  and the
seminorm version of the Hahn--Banach theorem
\cite[Theorem~3.3]{R91}, for each $n\in\mathbb N$,
there exists $\Psi_n\in G^*$ such that
\begin{equation}\label{eq:HBcontrol}
	\Psi_n|_Y=\phi_n,
	\qquad
	|\Psi_n(x)|
	\le\int_{\beta\mathbb N}\widehat r_x\,
	d\widehat\nu_{\phi_n},
	\quad x\in G, 
\end{equation} 
where $\widehat r_x$
is   as in
\eqref{propertieswidehatrx}. 
In particular,
\[
\sup_{n\in\mathbb N}\norm{\Psi_n}_{G^*}
\stackrel{\widehat r_x\le \norm{x}_G}{=}
\sup_{n\in\mathbb N}\norm{\phi_n}_{Y^*}
\stackrel{\text{\cite[D.12]{AK06}}}{<}\infty.
\]
By the Banach--Alaoglu theorem
\cite[F.5]{AK06},
$\{\Psi_n:n\in\mathbb N\}$ has a weak$^*$ compact
closure in $G^*$. 
Hence, by the finite intersection property of $\{ \overline{\{\Psi_n:n\in\mathbb N,\ n\ge N\}}
^{\,\sigma(G^*,G)} \}_{N\in \mathbb N}$ \cite[Definition, p.~169]{M20} and  \cite[Theorems~26.2 and~26.9]{M20},
  the set
\[
K:=
\bigcap_{N\in\mathbb N}
\overline{\{\Psi_n:n\in\mathbb N,\ n\ge N\}}
^{\,\sigma(G^*,G)} 
\]
is nonempty and weak$^*$ compact in $G^*$.

	We next show that $K$ is weakly compact in $G^*$. 
	Fix $\Psi\in K\subseteq G^*$. 
By the definition of $K$, 
$\Psi\in\overline{\{\Psi_n:n\in \mathbb N, n\ge N\}}^{\sigma(G^*,G)}$ for all $N\in \mathbb N$.
Hence, 	there exists a subnet
	$(\Psi_{n_i})_{i\in I}$, indexed by a directed set $I$,
	such that 
	\begin{align}\label{Psini}
	\Psi_{n_i}\rightarrow\Psi
	\quad\text{in the  }\sigma(G^*,G)\text{-topology}.
\end{align}
By the weak compactness of the weak closure of
	$\left\{ 
C\frac{d\widehat\nu_{\phi_{n_i}}}{d\lambda}: i\in I 
\right\}$ in $L_1(\beta\mathbb N,\lambda)_+$,
there exists a function
$f\in L_1(\beta\mathbb N,\lambda)_+$
belonging to the weak closure of $
\left\{
C\frac{d\widehat\nu_{\phi_{n_i}}}{d\lambda}:
i\in I
\right\}
$ \cite[Theorem~4.29]{F99}
such that, after passing to
a further subnet and retaining the same notation,
\begin{align}\label{1.1}
C\frac{d\widehat\nu_{\phi_{n_i}}}{d\lambda}
\rightarrow f
\quad\text{weakly in }L_1(\beta\mathbb N,\lambda).
\end{align}

	Since $\Psi_{n_i}|_Y
	\stackrel{\eqref{eq:HBcontrol}}{=}
	\phi_{n_i}$ and
	$(\phi_n)_{n\in\mathbb N}\subseteq Y^*$ is weak$^*$ null,
it follows from \eqref{Psini} that  \begin{align}\label{PsiY=0}
	 \Psi|_Y=0.
\end{align}
For every $x\in G$, we have
$\widehat r_x\in C(\beta\mathbb N)
\subseteq L_\infty(\beta\mathbb N,\lambda)$.
Thus, 
\begin{eqnarray}\label{Psixineq}
	|\Psi(x)|
	& \notag \stackrel{\eqref{Psini}}{=}&
	\lim_i|\Psi_{n_i}(x)|
	\stackrel{\eqref{eq:HBcontrol}}{\le}
	\lim_i\int_{\beta\mathbb N}
	\widehat r_x\,d\widehat\nu_{\phi_{n_i}}\\
	& =&C^{-1}\lim_i\int_{\beta\mathbb N}
	\left(C\frac{d\widehat\nu_{\phi_{n_i}}}{d\lambda}\right)
	\widehat r_x\,d\lambda
    \stackrel{\eqref{1.1}}{=}
	C^{-1}\int_{\beta\mathbb N}
	f\widehat r_x\,d\lambda.
\end{eqnarray}
For each $z\in Q$, choose $x\in G$ such that $q(x)=z$, and define 
\begin{align} \label{the def of chi}
\chi(z):=\Psi(x).
\end{align}
Since $\Psi|_Y=0$ by \eqref{PsiY=0}  and $\ker q=Y$, it follows that the map $\chi$
is well-defined and linear.
For every $z\in Q$ and $\varepsilon>0$, choose
$x\in G$ satisfying \eqref{decompzkxy}. 
Then,
\[
|\chi(z)|=|\Psi(x)|
\le \norm{\Psi}_{G^*}\norm{x}_G
\le \norm{\Psi}_{G^*}
   \bigl(C\norm{z}_Q+\varepsilon\bigr).
\]
Letting $\varepsilon\downarrow0$, we obtain
\[
|\chi(z)|\le C\norm{\Psi}_{G^*}\norm{z}_Q,
\]
and hence $\chi\in Q^*$.
Moreover, for every $x\in G$,
\begin{align} \label{=Psi(x)}
\left(q^*(\chi)\right)(x)=\chi(q(x))\stackrel{q(x)=z}{=}\Psi(x).
\end{align}
Therefore, $q^*(\chi)=\Psi$.

	For $z=(z_k)_{k\in\mathbb N}\in Q$ and
	$\varepsilon>0$, choose $x=(x_k)_{k\in\mathbb N}\in G$
	satisfying \eqref{decompzkxy}.  Define
$
	r_z(k):=\norm{z_k}_{Q_k}$ for all
$ k\in\mathbb N
$
and 	let $\widehat r_z\in C(\beta\mathbb N)$
	denote the  continuous extension of $r_z$, as in
	\eqref{propertieswidehatrx}.
Since	 $\mathbb N$ is dense in
	$\beta\mathbb N$, 
	it follows from  \eqref{1.4} that 
	\begin{align} \label{rx le rz}
	\widehat r_x(t)
	\le C\widehat r_z(t)+\varepsilon,
	\qquad t\in\beta\mathbb N.
	\end{align}
	Consequently,
	\[ 
	\begin{aligned}
		|\chi(z)|
		&\stackrel{\eqref{qx=z}}{=}|\chi(q(x))|
		\stackrel{\eqref{=Psi(x)}}{=}|\Psi(x)|\\
		&\overset{\eqref{Psixineq}}{\le}
		 C^{-1}
     		\int_{\beta\mathbb N}f\widehat r_x\,d\lambda
     		\stackrel{\eqref{rx le rz}}{\le}
		\int_{\beta\mathbb N}f\widehat r_z\,d\lambda
		+C^{-1}\varepsilon
		\int_{\beta\mathbb N}f\,d\lambda.
	\end{aligned}
	\]
	Letting $\varepsilon\downarrow0$, we obtain
	\[
	|\chi(z)|
	\le\int_{\beta\mathbb N}f\widehat r_z\,d\lambda,
	\qquad z\in Q.
	\]
    This implies that $\chi$ satisfies \eqref{the family}.
	
	Applying Lemma~\ref{lem:controlled-families} to
	$(Q_k,\norm{\cdot}_{Q_k})_{k\in\mathbb N}$ and the
	weak closure of 
		$\left\{ 
	C\frac{d\widehat\nu_{\phi_n}}{d\lambda}: n\in\mathbb N 
	\right\} $, we obtain that
$\{\chi\}\subseteq  Q^*$ defined by \eqref{the def of chi} is  relatively weakly compact. 
	Since $\Psi=q^*(\chi)$ and $q^*$ is weak-to-weak continuous \cite[Theorem 3.10]{B11}, it follows that 
	$K$ is
	relatively weakly compact in $G^*$.
	Since $K$ is weak$^*$ closed, it follows that  $K$ is  weakly
	closed. 
	Therefore, $K$ is weakly compact in $G^*$.     
\end{proof}

Mart\'inez-Cervantes and Rodr\'iguez
\cite[Corollary~1.4]{MR} proved that a 
subspace
of a Grothendieck space is Grothendieck whenever
the quotient is reflexive.  
However,  reflexivity is not preserved under infinite $\ell_\infty$-sums.

The following lemma establishes the Grothendieck property
for the $\ell_\infty$-sum of the subspaces  when the coordinate
quotients admit uniformly equivalent norms with a common
modulus of uniform convexity.

\begin{lemma}\label{lem:coordinate-kernel}

	Suppose that the Banach space 
	$
	(G, \norm{\cdot}_{G})
$ 
	has the Grothendieck property.	
Assume that each quotient $Q_k$ admits a norm
$\norm{\cdot}_{Q_k}$ equivalent to its quotient norm $\norm{\cdot}_{E_k/Y_k}$, with
equivalence constants independent of $k$, and that
the spaces $(Q_k, \norm{\cdot}_{Q_k})_{k\in\mathbb N}$
admit a common modulus of uniform convexity.
	Then,
$
	Y:=\mathop{\bigoplus}\limits_{k\in\mathbb N}^{\ell_\infty}Y_k,
	$
	equipped with the norm inherited from $G$, has the
	Grothendieck property.
\end{lemma}

\begin{proof}

Let $(\phi_n)_{n\in\mathbb N}$ be a weak$^*$ null
sequence in $Y^*$.
By Lemma~\ref{lem:controlled-extensions}, there exists
a norm-bounded sequence $(\Psi_n)_{n\in\mathbb N}$
in $G^*$ with $\Psi_n|_Y=\phi_n$ for every
$n\in\mathbb N$, such that
\[
K:=
\bigcap_{N\in\mathbb N}
\overline{\{\Psi_n:n\ge N\}}^{\,\sigma(G^*,G)}
\]
is weakly compact in $G^*$. 
By Lemma~\ref{lem:weak-star-cluster}, there exist
a strictly increasing sequence $(n_j)_{j\in\mathbb N}$
and $\Psi\in K$ such that
\[
\Psi_{n_j}\rightarrow\Psi
\quad\text{in }\sigma(G^*,G).
\] 
Since $\Psi_{n_j}|_Y=\phi_{n_j}$ and
$(\phi_n)_{n\in\mathbb N}$ is weak$^*$ null,
it follows that  $\Psi|_Y=0$.
Since $G$ has the Grothendieck property, it follows that 
$\Psi_{n_j}\to\Psi$ weakly in $G^*$.
Since the restriction operator
$R:G^*\to Y^*$, $R(\psi)=\psi|_Y$, is bounded and linear,
it follows that 
\begin{align} \label{adjoint R}
R^*:Y^{**}\to G^{**},\qquad R^*(y)(g)=y(Rg), \quad y\in Y^{**},\, g\in G^{*}.
\end{align}
Then, for all $y\in Y^{**}$, by the weak convergence of $\Psi_{n_j}$ to $\Psi$, we have 
\begin{eqnarray} \label{www}
y(\Psi_{n_j}|_Y) &\notag=&y(R(\Psi_{n_j}))  \\ &  \notag \stackrel{\eqref{adjoint R}}{=}&R^*(y)(\Psi_{n_j})
\to
R^*(y)(\Psi) \\ & 
\stackrel{\eqref{adjoint R}}{=}&y(R\Psi)  =y(\Psi|_Y).
\end{eqnarray}
This implies that
\begin{align}\label{subseconverge0weakly} 
\phi_{n_j}
=\Psi_{n_j}|_Y
\xrightarrow{\eqref{www}}\Psi|_Y=
0
\quad\text{weakly in }Y^*.
\end{align}
Since $(\phi_n)_{n\in\mathbb N}$ is an arbitrary
weak$^*$ null sequence in $Y^*$, it follows from \eqref{subseconverge0weakly}
 that every weak$^*$ null sequence in $Y^*$
admits a weakly null subsequence. 
Noting that every subsequence of $(\phi_n)_{n\in\mathbb N}$
is weak$^*$ null, we have that any subsequence admits a further
subsequence converging weakly to zero. 
By \cite[Theorem 2.3.3]{D19}, we have that $\phi_n\to 0$ weakly. 
Since $(\phi_n)_{n\in\mathbb N}$ is arbitrary, it follows that 
$Y$ has the Grothendieck property.
\end{proof}

 \subsection{The Grothendieck property of
$\ell_\infty$-sums of noncommutative $L_\infty\cap L_2$ spaces}
Let $\M$ be a von Neumann algebra acting on a Hilbert space
$H$. 
For each integer $n\ge1$, let $\mathbb M_n(\M)$
denote the von Neumann algebra of all $n\times n$
matrices with entries in $\M$,
equipped with the norm
$$\norm{(x_{ij})}_{\mathbb M_n(\cM)}=\sup_{\sum_{j}\norm{\xi_j}_H^2\le1}\left(\sum_{i=1}^n\norm{\sum_{j=1}^nx_{ij}\xi_j}_H^2\right)^{1/2}, \qquad (x_{ij})\in \mathbb M_n(\cM).$$
Following \cite{BlecherPaulsen1991}  and \cite[p.~525]{JP08},   
let
 $R_n(\M)$ denote  {\it the row  subspace of  $\mathbb M_n(\M)$}, consisting of matrices whose
 entries outside the first row are zero and
 equipped with the operator norm inherited from
 $\mathbb M_n(\M)$.

\begin{lemma}\label{lem:corner-maps}
	Let $(\M, \tau)$ be a semifinite von Neumann algebra. 
	Let $a\in(0,\infty)$ and let  $e\in\M$ be a nonzero  projection
	with $\tau(e)<\infty$. Regard $\cN:=e\M e$ as a Banach space with the norm
	\begin{align}\label{defcNnormmax}
	\norm{x}_{\cN}:=
	\max\left\{
	\norm{x}_{\M},
	a\norm{x}_{L_2(\M,\tau)}
	\right\},
	\qquad x\in\cN.
	\end{align}
	Then, there exist a von Neumann algebra $\cA$
	and two bounded linear operators
	$ 
	T:\left(\cN, \norm{\cdot}_\cN \right) \to \left(\cA,\norm{\cdot}_\infty\right)$ and $
 S:\left(\cA,\norm{\cdot}_\infty\right) \to \left(\cN, \norm{\cdot}_\cN \right)
	$ 
	such that
	$
	S\circ T=\mathbf 1_{\cN}$ and  $
\norm{T},\norm{S}\le C,
	$   
	where $C>0$ is a universal constant independent
	of $\M$, $\tau$, $a$, and $e$.\footnote{It is obvious that the norm $\norm{\cdot}_\cN$ on $\cN$ is equivalent with $\norm{\cdot}_\cM$. 
    However, the equivalence constant here depends on $a,\tau$ and $e$.  
     Taking $\cA=\cN$, the canonical inclusion $T(x)=x$ and
the compression $S(z)=eze$ satisfy
$S\circ T=\mathbf1_{\cN}$, with
$\|T\|\le1$ and
$\|S\|\le\max\{1,a\tau(e)^{1/2}\}$,  see \eqref{equivinftynormN}. 
In this  lemma,  $\cA$, $T$ and $S$ can
be chosen so that both operator norms are bounded by
a universal constant.}
\end{lemma}

\begin{proof} 
It suffices to prove the lemma when $a=1$.
Indeed, the general case follows by replacing $\tau$
with $a^2\tau$, since
\[
\norm{x}_{L_2(\M,a^2\tau)}
=a\norm{x}_{L_2(\M,\tau)},
\qquad x\in\cN.
\]
Moreover, $a^2\tau$ is a faithful normal semifinite trace
and $(a^2\tau)(e)<\infty$.
We therefore assume that $a=1$.


Note that $\cN=e\M e$ is a von Neumann algebra
with identity $e$, equipped with the operator norm
inherited from $\M$.
Since $0<\tau(e)<\infty$, it follows that  the restriction
$\tau|_{\cN}$ is a finite faithful normal trace.
Hence,
$ 
\varphi(\cdot ):=\frac{1}{\tau(e)} \tau(\cdot) 
$ 
defines a faithful normal tracial state on $\cN$.

For every $x\in\cN$, we have 	
	\begin{align}\label{2normL2M2}
	\norm{x}_{L_2(\M,\tau)}^2
	=\tau(x^*x)
	\le\tau(e)\norm{x}_{\M}^2<\infty, 
	\end{align}
which implies that  $\cN\subseteq L_2(\M,\tau)$, and
	\begin{align}\label{equivinftynormN}
	\norm{x}_{\M}
	\stackrel{\eqref{defcNnormmax}}{\le}
	\norm{x}_{\cN}
 \stackrel{\eqref{defcNnormmax}, \eqref{2normL2M2}
 }{\le}	\max\{1,
\sqrt{\tau(e)}\}\norm{x}_{\M},
	\qquad x\in\cN.
\end{align}
	Moreover,   for every $x\in\cN$,
	\begin{align}\label{xL2eMkevarphi2=frac}
		\|x\|_{L_2(\cN,\varphi)}^2
		&\notag =\varphi(x^*x)=\frac{\tau(x^*x)}{\tau(e)}
		=\frac{\|x\|_{L_2(\M,\tau)}^2}{\tau(e)},
		\\ 
		\norm{x}_{L_2(\cN,\varphi)}^2
		& =\varphi(x^*x)
		\le\norm{x}_{\M}^2\varphi(e)
		=\norm{x}_{\M}^2.
	\end{align} 

Define the  conditional expectation
	\begin{align*} 
		\mathcal E: \cN &\notag \to \C e  \\ 
	x&\mapsto
	\varphi(x)e=\frac{\tau(x)}{\tau(e)}e,
	\qquad x\in\cN. 
	\end{align*}
	For every $x\in\cN$, we have 
	\begin{align}\label{cExx*value}
	\norm{\left(\mathcal E (xx^*) \right)^{1/2}
	}_{\M}
 	=	\left\|
		\left(\frac{\tau(xx^*)}{\tau(e)}e\right)^{1/2}
		\right\|_{\M}
		=
		\left(\frac{\tau(x^* x)}{\tau(e)}\right)^{1/2}
		=
		\frac{\norm{x}_{L_2(\M,\tau)}}{\sqrt{\tau(e)}}.
	\end{align}

Now choose
$
n:=\lceil 
\tau(e)\rceil\geq 1, 
$
where $\lceil t\rceil$ denotes the smallest integer
greater than or equal to $t$.
Following \cite{JP08}, 
let $R_{\infty,1}^n(\cN,\cE)$ denote the vector
space $\cN$ equipped with the norm
\begin{align}\label{normParcet}
\norm{x}_{R_{\infty,1}^n(\cN,\cE)}
:=
\max\left\{
\norm{x}_{\M},
\sqrt n\,\norm{\left(\cE(xx^*)^{1/2}\right)}_{\M}
\right\},
\qquad x\in\cN.
\end{align} 
Since
$
\tau(e)\le n\le 
\tau(e)+1,
$
  it follows that for every $x\in\cN$, 
\begin{eqnarray*}
	\norm{x}_{L_2(\M,\tau)}^2 
	&\stackrel{\eqref{xL2eMkevarphi2=frac}}{=}&
	\tau(e) 
	\norm{x}_{L_2(\cN,\varphi)}^2	\leq 
	\frac{n }{\tau(e)} \norm{x}_{L_2(\M,\tau)}^2  
	\\ & \stackrel{\eqref{xL2eMkevarphi2=frac}}{=}&
	n\norm{x}_{L_2(\cN,\varphi)}^2 
	\leq 
 \tau(e) 
	\norm{x}_{L_2(\cN,\varphi)}^2 
	+
	\norm{x}_{L_2(\cN,\varphi)}^2 
	\\
	&\stackrel{\eqref{xL2eMkevarphi2=frac}
    }{\leq}& 
	\norm{x}_{L_2(\M,\tau)}^2
	+\norm{x}_{\M}^2
	\stackrel{\eqref{defcNnormmax}}{\leq }
	2\norm{x}_{\cN}^2
	, 
\end{eqnarray*}
which implies that 
\begin{eqnarray}\label{xNnormsqrt2x}
	\norm{x}_{\cN}
	&\notag =&\max\left\{
	\norm{x}_{\M},
	\norm{x}_{L_2(\M,\tau)}
	\right\}
	\\
	& \notag \stackrel{n\ge \tau(e)}{\le}&
	\max\left\{
	\norm{x}_{\M},
	\sqrt{\frac{n}{\tau(e)}}
	\norm{x}_{L_2(\M,\tau)}
	\right\}\\
	&\stackrel{\eqref{cExx*value}, \eqref{normParcet}}{=}&\norm{x}_{R_{\infty,1}^n(\cN,\cE)}
	\stackrel{\tiny \substack{\eqref{defcNnormmax}\\\tiny \frac{n }{\tau(e)} \norm{x}_{L_2(\M,\tau)}^2\le 2\norm{x}_\cN^2}}{\le}
\sqrt2\,\norm{x}_{\cN}.
\end{eqnarray}
This   yields that $R_{\infty,1}^n\left(\cN,\cE\right)$ coincides as sets with $\cN$ and their norms are equivalent.

By \cite[Lemma~1.3]{JP08}, there exist a von Neumann
algebra $\cB$ with  identity $\mathbf1_{\cB}$, an isomorphic embedding
\[
T^0:R_{\infty,1}^n(\cN,\cE)\to R_n(\cB),
\]
and a bounded linear  operator 
$
P:R_n(\cB)\to T^0(\cN)
$
such that $P^2=P$,  for a universal constant $c_1\ge1$,
\begin{align}\label{c1constant}
c_1^{-1}\norm{x}_{R_{\infty,1}^n(\cN,\cE)}
\le\norm{T^0(x)}_{R_n(\cB)}
\le c_1\norm{x}_{R_{\infty,1}^n(\cN,\cE)},
\qquad x\in\cN,
\end{align}
and
\begin{align}\label{PT0x=T0x}
P(T^0(x))=T^0(x),
\qquad x\in\cN,
\qquad \norm{P}\le c_1.
\end{align}
Since 
\[
\begin{aligned}
	\norm{(T^0)^{-1}(T^0(x))}_{\cN}
=\norm{x}_{\cN}
	\stackrel{\eqref{cExx*value},\eqref{normParcet}}{\le}\norm{x}_{R_{\infty,1}^n(\cN,\cE)}
	 \stackrel{\eqref{c1constant}}{\le} c_1\norm{T^0(x)}_{R_n(\cB)},
	\quad x\in\cN,
\end{aligned}
\]
it follows that 
the inverse operator
$
(T^0)^{-1}:T^0(\cN)\to
(\cN,\norm{\cdot}_{\cN})
$
satisfies $\norm{(T^0)^{-1}}\le c_1$.

Define a bounded linear  projection   by 
	\begin{align*}
	P_{\mathrm{row}}: \mathbb M_n(\cB) &\notag \to R_n(\cB)
	\\ z & \mapsto \operatorname{diag}(\mathbf1_{\cB},0,\ldots,0) \cdot z.
	\end{align*}
	Then, $\norm{P_{\mathrm{row}}}=1$.
	Define
	\begin{align}\label{defTT0}
	T:\cN\to  \mathbb M_n(\cB),
	\qquad T(x):=T^0(x),
	\qquad x\in\cN,
	\end{align}
	and
	\[
	S:=(T^0)^{-1}\circ P\circ P_{\mathrm{row}}
	: \mathbb M_n(\cB) \to (\cN,\norm{\cdot}_{\cN}).
	\]
	For every $x\in\cN$, we have
$
	P_{\mathrm{row}}(T^0(x))=T^0(x)
$,
	and hence
	\begin{align*}
	S(T(x)) & =	S(
	T^0(x)
	)
	=\left( 
	(T^0)^{-1}\circ P\circ P_{\mathrm{row}}
	\right) (T^0(x))
	\\ & 	=\left( 
	(T^0)^{-1}\circ P
	\right) (T^0(x))
\stackrel{\eqref{PT0x=T0x}}{=}	(T^0)^{-1}(T^0(x))
\\ & 	=x=\mathbf 1_{\cN}(x).
	\end{align*}
	Furthermore,
	 \[
		\norm{T(x)}_{\mathbb M_n(\cB) }
	\stackrel{\eqref{defTT0}}{=}
		\norm{T^0(x)}_{R_n(\cB)}
		\stackrel{\eqref{c1constant}}{\le} c_1\norm{x}_{R_{\infty,1}^n(\cN,\cE)}
	\stackrel{\eqref{xNnormsqrt2x}}{	\le}
	\sqrt2\,c_1\norm{x}_{\cN}, \qquad x\in\cN
	\]  and 
	\[
	\norm{S}
	\le
	\norm{(T^0)^{-1}}\norm{P}\norm{P_{\mathrm{row}}}
	\le c_1^2.
	\]
Taking $\cA:=\mathbb M_n(\cB)$ and
$C:=\max\{\sqrt2\,c_1,c_1^2\}$, we obtain
\[
S\circ T=\mathbf 1_{\cN},
\qquad \norm{T},\norm{S}\le C.
\]
Since $c_1$ is universal, we have that $C$ is independent of
$\M$, $\tau$, $a$, and $e$.
This completes the proof.
\end{proof}

\begin{lemma}\label{lem:compression-duality}
	Let $(\M, \tau)$ be a semifinite von Neumann  algebra  
and let  $a\in(0,\infty)$. 
	Equip
$
	F:=\M\cap L_2(\M,\tau)
$  and $	V:=L_1(\M,\tau)+L_2(\M,\tau)
$
	with the Banach  norms
	\begin{align}
	\norm{x}_F & \notag 
	:=
	\max\left\{
	\norm{x}_{\M},
	a\norm{x}_{L_2(\M,\tau)}
	\right\},
	\quad x\in F,
\\ 
	\norm{y}_V & \notag 
	:=
	\inf_{\substack{	y=u+v \\ 
			u\in L_1(\M,\tau),\ v\in L_2(\M,\tau)}}
	\left(
	\norm{u}_{L_1(\M,\tau)}
	+a^{-1}\norm{v}_{L_2(\M,\tau)}
	\right),
	\quad y\in V, 
	\end{align} 
	respectively. 
Defining, for each $x\in F$, the linear functional
$\phi_x:y\mapsto\tau(xy)$ on $V$, then 
the map $x\mapsto\phi_x$ is a linear isometry
from $F$ onto the dual space $V^*$ of $V$. 
	Moreover, let $(e_\alpha)$ be a net of finite-trace
	projections in $\M$ with $e_\alpha\uparrow\mathbf1$. 
    Then, 
\begin{equation*}
		e_\alpha xe_\alpha \rightarrow x
		\quad\text{weak$^*$ in }F,
		\qquad x\in F.
	\end{equation*}
\end{lemma}

\begin{proof}
Since $F$ is the Banach dual (and K\"othe dual) of 
 $V$, it follows from 
\cite[Proposition 6.1.9]{DPS} and \cite[Theorem 2.7.1 and p. 54, (2)]{BL12} that
$x\mapsto\phi_x$ is surjective and satisfies
\[
\norm{\phi_x}_{V^*}
=
\max\left\{
\norm{x}_{\M},
a\norm{x}_{L_2(\M,\tau)}
\right\}
=\norm{x}_F,
\qquad x\in F.
\]

Since ${\bf 1}-e_\alpha \downarrow 0$, it follows from \cite[Theorem 5.5.12]{DPS} that for every $y\in  L_2(\M,\tau)$, 
	\[
	\begin{aligned}
	\norm{(\mathbf1-e_{\alpha})y}_{L_2(\M,\tau)}
		 	\rightarrow0
	\end{aligned}
	\]
	 and 
	\[
	\begin{aligned}
		\norm{y(\mathbf1-e_\alpha)}_{L_2(\M,\tau)} 
		\rightarrow0.
	\end{aligned}
	\]
	Consequently, 	for a fixed  $x\in F\subseteq L_2\left(\cM,\tau\right)$, 
	\begin{align}\label{1.5}
		\norm{e_{\alpha} x
		e_{\alpha}-x}_{L_2(\M,\tau)}&\le\norm{(\mathbf1-e_{\alpha})xe_\alpha}_{L_2(\M,\tau)}
		+\norm{x(\mathbf1-e_{\alpha})}_{L_2(\M,\tau)}\notag\\
		&\le
		\norm{(\mathbf1-e_{\alpha})x}_{L_2(\M,\tau)}
		+\norm{x(\mathbf1-e_{\alpha})}_{L_2(\M,\tau)}\rightarrow0.
	\end{align}

 Since $e_\alpha \uparrow {\bf 1}$, it follows that for all $u\in L_1(\M,\tau)$, 
	\begin{eqnarray} \label{u in L1}
	|\tau((e_\alpha xe_\alpha-x)u)|
  &  \stackrel{\mbox{\tiny \cite[Coro. 3.4.6]{DPS}}}{\le}& \tau\left(|\left(e_\alpha xe_\alpha-e_\alpha x\right)u|\right) +\tau\left(|(e_\alpha x -x)u|\right)\notag\\
&    =& \norm{(e_\alpha xe_\alpha-e_\alpha x)u}_{L_1(\M,\tau)} + \norm{(e_\alpha x - x)u}_{L_1(\M,\tau)}
    \notag\\
    &\stackrel{\mbox{\tiny \cite[Thm. 5.5.12]{DPS}}}{\rightarrow}&0.
	\end{eqnarray}
	For every $v\in L_2(\M,\tau)$, by the
	noncommutative H\"older inequality \cite[Theorem 1]{Sukochev16} (see also \cite[p. 157]{DPS}), we have that
    \begin{align} \label{v in L2}
	|\tau((e_\alpha xe_\alpha
	-x)v)|
	\le
	\norm{e_\alpha xe_\alpha -x}_{L_2(\M,\tau)}
	\norm{v}_{L_2(\M,\tau)}
	\stackrel{\eqref{1.5}}{\rightarrow}
	0.
	\end{align}
	Note that every $y\in V$ admits a decomposition
	$y=u+v$, where $u\in L_1(\M,\tau)$ and
	$v\in L_2(\M,\tau)$.
Then,
\[
\phi_{e_\alpha xe_\alpha}(y)-\phi_x(y)
=
\tau((e_\alpha xe_\alpha-x)u)
+\tau((e_\alpha xe_\alpha-x)v)
\stackrel{\eqref{u in L1}, ~\eqref{v in L2}}{\rightarrow}0.
\]
Thus, $\phi_{e_\alpha xe_\alpha}\to\phi_x$
in $\sigma(V^*,V)$.
Under the isometric identification $F=V^*$,
we therefore obtain
 $e_\alpha xe_\alpha\to x$
in $\sigma(F,V)$.
\end{proof}

The following lemma is a consequence of inheritance of Grothendieck property for complemented subspaces/surjections. 

\begin{lemma}\label{lem:local-maps}
Let $(V_k)_{k\in\mathbb N}$ be a family of Banach spaces,
and set
$
G:=
\mathop{\bigoplus}\limits_{k\in\mathbb N}^{\ell_\infty}V_k^*. 
$ 
	Let $D$ be a nonempty directed set and let
	$(\mathcal B_d)_{d\in D}$ be a family of von Neumann algebras.
	Suppose that
	bounded linear operators
$
	T_d:G\to\mathcal B_d
$  and  $
	R_d:\mathcal B_d\to G
$
	satisfy 
\[
\sup_{d\in D}\max\{\norm{T_d},\norm{R_d}\}<\infty
\]
and, for every $x=(x_k)\in G$ and $k\in\mathbb N$,
the net
$\bigl(
(R_d\circ T_d) (x)
\bigr)_k$,
indexed by $D$,
converges to $x_k$ in the  $\sigma(V_k^*,V_k)$-topology.
Then, $G$ has the Grothendieck property.
\end{lemma}

\begin{proof}
	Define the von Neumann algebra $\mathcal B:=
	\mathop{\bigoplus}\limits_{d\in D}^{\ell_\infty}
	\mathcal B_d$.
	Define
	\[
	T:G\to\mathcal B,
	\qquad
	T(x):=(T_d(x))_{d\in D},
	\qquad x\in G.
	\]
	By the uniform boundedness of $(T_d)_{d\in D}$, we have that
    $T$ is well-defined and linear, with
	\[
	\norm{T}\le\sup_{d\in D}\norm{T_d}<\infty.
	\]

	Let $m\ge1$ and $d_1,\ldots,d_m\in D$.
	There exists
	$d_*\in D$ such that $d_*\ge d_j$ for every
	$1\le j\le m$.
	Consequently, 
$
	d_*\in
	\bigcap_{j=1}^{m}\{d\in D:d\ge d_j\}.
	$ 
	Thus, for every integer $m\ge1$ and every choice of
	$d_1,\ldots,d_m\in D$, we have 
	\[
	\bigcap_{j=1}^{m}\{d\in D:d\ge d_j\}
	\ne\varnothing, 
	\]
i.e.,  the family
	$
	\bigl\{\{d\in D:d\ge d_0\}:d_0\in D\bigr\}
$
	has the finite intersection property \cite[Definition, p.~169]{M20}.
	By the ultrafilter lemma \cite[4A1I(a)]{F03}, there exists an ultrafilter
	$\mathcal U$ on $D$ such that 
$
	\{d\in D:d\ge d_0\}\in\mathcal U$ 
for every $ d_0\in D.
$

Fix $z=(z_d)_{d\in D}\in\mathcal B=	\mathop{\bigoplus}\limits_{d\in D}^{\ell_\infty}
\mathcal B_d$ and
$k\in\mathbb N$.
For every $d\in D$, we have  $R_d(z_d)
\in G=\mathop{\bigoplus}\limits_{k\in\mathbb N}^{\ell_\infty}V_k^*$. 
For every $y\in V_k$ and $d\in D$, we obtain
\begin{align}\label{f z,k}
\left|\bigl(R_d(z_d)\bigr)_k(y)\right|
\le\norm{R_d}\norm{z_d}_{\mathcal B_d}
\norm{y}_{V_k}
\le 
\left(\sup_{d\in D}\norm{R_d}\right)
\norm{z}_{\mathcal B}\norm{y}_{V_k}.
\end{align}
This implies that for each fixed $y\in V_k$, 
the scalar net
$
\left\{ \bigl(R_d(z_d)\bigr)_k(y)
\right\}_{d\in D}
$
takes values in the closed disk
$
\left\{\zeta\in\mathbb C:
|\zeta|\le  \left(\sup_{d\in D}\norm{R_d}\right) \norm{z}_{\mathcal B}\norm{y}_{V_k}
\right\},
$
which is compact and Hausdorff in the usual
topology of $\mathbb C$.
This together with \cite[Theorem~3.1.10]{G22} implies that
the scalar ultralimit
$
\mathop{ \lim}\limits _{d\to\mathcal U}
\bigl(R_d(z_d)\bigr)_k(y)
$
exists uniquely for every $y\in V_k$. 
Taking ultralimits in  the 
estimate \eqref{f z,k} and using  continuity of the absolute value, we obtain  that 
the map
\begin{align}\label{mapfromVk}
V_k\to\mathbb C,
\qquad
y\mapsto
\lim_{d\to\mathcal U}
\bigl(R_d(z_d)\bigr)_k(y)
\end{align}
is a bounded  functional with norm at most
$\left(\sup_{d\in D}\norm{R_d}\right) \norm{z}_{\mathcal B}$. It is straightforward to verify that the map defined in \eqref{mapfromVk}
is linear.

For each $z=(z_d)_{d\in D}\in\mathcal B$ and
$k\in\mathbb N$, define 
$(S(z))_k\in V_k^*$ by
$$
	\bigl((S(z))_k\bigr)(y)
	:=
	\lim_{d\to\mathcal U}
	\bigl(
	R_d(z_d)\bigr)_k
	(y),
	\qquad y\in V_k.
	$$
	Consequently,
$
	S(z):=\big(  (S(z))_k
	\big)_{k\in\mathbb N}
	\in G,
$
	and
	\[
	\norm{S(z)}_G
	\le
	\left(\sup_{d\in D}\norm{R_d}\right)
	\norm{z}_{\mathcal B}.
	\]
Moreover, by the	linearity of each $R_d$, 
	  $S:\mathcal B\to G$ is linear.
	
	Fix $x=(x_k)_{k\in \mathbb N}\in G$, $k\in\mathbb N$, and $y\in V_k$. 
	By $\bigl(
(R_d\circ T_d) (x)
\bigr)_k \to_D x_k$  in
	$\sigma(V_k^*,V_k)$, for every $\varepsilon>0$,
	there exists $d_0\in D$ such that
	\begin{align}\label{convergenceassump}
	\left|
	\bigl(R_d(T_d(x))\bigr)_k(y)-x_k(y)
	\right|<\varepsilon
	\qquad\text{whenever }d\ge d_0.
	\end{align}
	In other words, 
	\[
	\{d\in D:d\ge d_0\}
	\subseteq
	\left\{
	d\in D:
	\left|
	\bigl(R_d(T_d(x))\bigr)_k(y)-x_k(y)
	\right|<\varepsilon
	\right\}.
	\]
	Since $\{d\in D:d\ge d_0\}\in\mathcal U$ for any $ d_0\in D$, it follows from  the upward-closure property of filters
	\cite[Definition~1.1.1(3)]{G22}   that 
 $$	\left\{
 d\in D:
 \left|
 \bigl(R_d(T_d(x))\bigr)_k(y)-x_k(y)
 \right|<\varepsilon
 \right\} \in 
 \mathcal U.$$
Since $\varepsilon>0$ is arbitrary, it follows from the definition
of an ultralimit \cite[Definition~3.1.1]{G22},  \eqref{convergenceassump} and the uniqueness of $\mathop{ \lim}\limits _{d\to\mathcal U}
\bigl(R_d(z_d)\bigr)_k(y)$ for any fixed $z=(z_d)_{d\in  D}\in \cB$  that
\[
\lim_{d\to\mathcal U}
\bigl(R_d(T_d(x))\bigr)_k(y)
=x_k(y).
\]
	By $(T(x))_d=T_d(x)$ and the definition of $S$, we have that
	\[
	\bigl(S(T(x))\bigr)_k(y)
	=
	\lim_{d\to\mathcal U}
	\bigl(R_d(T_d(x))\bigr)_k(y)
	=x_k(y).
	\]
	Therefore, $(S(T(x)))_k=x_k$ for every
	$k\in\mathbb N$ and $x\in G$.
	Thus, $S\circ T=\mathbf 1_G$, and $S:\mathcal B\to G$
	is surjective.
	
Since  $\mathcal B$ has the Grothendieck property
	by \cite[Corollary~7]{P94},
    it follows from Proposition~\ref{surjective image}(2)  that   $G$ has the Grothendieck property.
\end{proof}

    For each $k\in\mathbb N$, let $(\M_k, \tau_k)$ be a semifinite 
von Neumann algebra, and let $a_k>0$.  
	Equip $F_k:=\M_k\cap L_2(\M_k,\tau_k)$ with the complete norm
\[
\norm{x}_{F_k}:=
\max\left\{
\norm{x}_{\M_k},
a_k\norm{x}_{L_2(\M_k,\tau_k)}
\right\}, \quad x\in F_k, 
\]
where $\norm{\cdot}_{\M_k}$ denotes the operator norm
on $\M_k$. 

The following lemma is the main result of this subsection. 

\begin{lemma}\label{lem:square-product}
 The Banach \(^*\)-algebra 
	$G:=
	\mathop{\bigoplus}\limits_{k\in\mathbb N}^{\ell_\infty}
	F_k$    
  has the
	Grothendieck property.
\end{lemma}

\begin{proof}
	Fix $k\in\mathbb N$ and a  projection
	$0\ne e\in\M_k$ with $\tau_k(e)<\infty$.
	For every $x\in e\M_ke$, we have
	\[
	\tau_k(x^*x)\le\tau_k(e)\norm{x}_{\M_k}^2<\infty.
	\]
	Thus, $e\M_ke\subseteq F_k$.
	Equip $e\M_ke$ with the norm $\norm{\cdot}_{F_k}$ inherited from $F_k$.
	Applying Lemma~\ref{lem:corner-maps},
	we obtain a von Neumann algebra $\mathcal B_{k,e}$
	and bounded linear operators
	\[
	T_{k,e}:e\M_ke\to\mathcal B_{k,e},
	\qquad
	S_{k,e}:\mathcal B_{k,e}\to e\M_ke
	\]
such that for all $x\in e\M_k e$ and $z\in\mathcal B_{k,e}$,
\begin{align}\label{eq:corner-maps}
  S_{k,e}\circ T_{k,e}=\mathbf 1_{e\M_k e},
 &\notag  \quad
  \norm{T_{k,e}(x)}_{\mathcal B_{k,e}}\le C\norm{x}_{F_k},
  \\ 
  \norm{S_{k,e}(z)}_{F_k} & \le C\norm{z}_{\mathcal B_{k,e}},
\end{align}
	where $C>0$ is independent of
	$k$, $\M_k$, $\tau_k$, $a_k$ and $e$.

	For each $k\in\mathbb N$, equip
	\[
	V_k:=L_1(\M_k,\tau_k)+L_2(\M_k,\tau_k)
	\]
	with the norm
	\[
	\norm{y}_{V_k}:=
	\inf_{\substack{	y=u+v \\ u\in L_1\left(\cM_k,\tau_k\right),v\in L_2\left(\cM_k,\tau_k\right)}}
	\left(
	\norm{u}_{L_1(\M_k,\tau_k)}
	+a_k^{-1}\norm{v}_{L_2(\M_k,\tau_k)}
	\right),
	\quad y\in V_k.
	\]
    where the infimum is taken over all such  $	u\in L_1(\M_k,\tau_k)$ and $
			v\in L_2(\M_k,\tau_k)$. 
	By Lemma~\ref{lem:compression-duality}, defining
	\[
	\phi_x(y):=\tau_k(xy),
	\qquad x\in F_k,\quad y\in V_k,
	\]
	the map $x\mapsto\phi_x$ is a surjective linear
	isometry from $F_k$ onto $V_k^*$.
	We identify $F_k$ with $V_k^*$ via $x\mapsto\phi_x$. 
Let $(e_\alpha)$ be the net of all finite-trace
projections in $\M_k$, ordered by inclusion of their ranges, that is, $e_\alpha\le e_\beta$ if and only if
$e_\alpha(H)\subseteq e_\beta(H)$.
By semifiniteness of the traces $\tau_k$, $e_\alpha\uparrow\mathbf1_{\M_k}$.
	Applying Lemma~\ref{lem:compression-duality}
	to this net, we obtain
	\begin{equation}\label{eq:compressionk}
		e_\alpha xe_\alpha\rightarrow x
		\quad\text{in }\sigma(F_k,V_k),
		\qquad x\in F_k.
	\end{equation}

	Let $D$ be the index set consisting of
$
	d=(F,(e_k)_{k\in F}),
$
	where $F\subseteq\mathbb N$ is nonempty and finite,
	and each $e_k\in\M_k$, $k\in F$, is a projection
	with $\tau_k(e_k)<\infty$.
	Order $D$ by
	\[
	(F,(e_k)_{k\in F})
	\le(F',(e'_k)_{k\in F'})
	\quad\Longleftrightarrow\quad
	F\subseteq F'
	\ \text{and}\
	e_k\le e'_k,\quad k\in F.
	\]
  It is readily verified that   $D$ is directed.

	For every $d=(F,(e_k)_{k\in F})\in D$ and $k\in F$,
	let $\mathcal B_{k,e_k}$ be the von Neumann algebra
	associated with the corner $e_k\M_ke_k$
	in \eqref{eq:corner-maps}.
	Define
	the von Neumann algebra
	\[
	\mathcal B_d:=
	\mathop{\bigoplus}\limits_{k\in F}^{\ell_\infty}
	\mathcal B_{k,e_k},
	\qquad
	\norm{z}_{\mathcal B_d}:=
	\max_{k\in F}\norm{z_k}_{\mathcal B_{k,e_k}},
	\quad z=(z_k)_{k\in F}\in\mathcal B_d.
	\]
	Define
	\[
	T_d:G\to\mathcal B_d,
	\quad
	T_d(x):=
	\bigl(T_{k,e_k}(e_kx_ke_k)\bigr)_{k\in F},
	\quad x=(x_k)_{k\in\mathbb N}\in G,
	\]
	and define $R_d:\mathcal B_d\to G$ by
	\[
	(R_d(z))_k:=
	\begin{cases}
		S_{k,e_k}(z_k),&k\in F,\\
		0,&k\in \mathbb N\setminus F,
	\end{cases}
	\qquad
	z=(z_j)_{j\in F}\in\mathcal B_d.
	\]
Note that  $R_d$ and $T_d$ are linear.
	By \eqref{eq:corner-maps},
	\[
	\begin{aligned}
		\norm{T_d(x)}_{\mathcal B_d}
		&\le C\max_{k\in F}\norm{e_kx_ke_k}_{F_k}
	\le C\max_{k\in F}\norm{x_k}_{F_k}	
        \leq 
		C\norm{x}_G,
		\quad x\in G,\\
		\norm{R_d(z)}_G
		&\le C\max_{k\in F}\norm{z_k}_{\mathcal B_{k,e_k}}
		=C\norm{z}_{\mathcal B_d},
		\quad z\in\mathcal B_d.
	\end{aligned}
	\]
	Consequently,
	\[
	\sup_{d\in D}\max\{\norm{T_d},\norm{R_d}\}
	\le C<\infty.
	\]
	Furthermore, by \eqref{eq:corner-maps}, we have that
	\begin{align}\label{defkinFRdTd}
	\bigl((R_d\circ T_d)(x)\bigr)_k
	=
	\begin{cases}
		e_kx_ke_k,&k\in F,\\
		0,&k\in \mathbb N\setminus F,
	\end{cases}
	\qquad x\in G.
	\end{align}

Fix $x\in G$, $k\in\mathbb N$, $y\in V_k$,
and $\varepsilon>0$. 
By \eqref{eq:compressionk}, there exists
a $\tau_k$-finite projection   $e\in\M_k$ such that for every $\tau_k$-finite projection $f\in\M_k$
satisfying $f\ge e$, we have 
\begin{align}\label{fgeqeconvergence}
\left|\tau_k((fx_kf-x_k)y)\right|<\varepsilon.
\end{align}
Set $d_0:=(\{k\},(e))\in D$.
For every $d=(F,(e_j)_{j\in F})\ge d_0$,
by the definition of the order on $D$, we have that
$k\in F$ and $e_k\ge e$.
Hence,  for every $d=(F,(e_j)_{j\in F})\ge d_0$,
\[
\begin{aligned}
	\left|
	\tau_k\left(
	\left(\bigl((R_d\circ T_d)(x)\bigr)_k-x_k\right)y
	\right)
	\right|
	\stackrel{\eqref{defkinFRdTd}}{=}
	\left|\tau_k((e_kx_ke_k-x_k)y)\right|
	\stackrel{\eqref{fgeqeconvergence}}{<} 
	\varepsilon.
\end{aligned}
\]
Since $y\in V_k$ and $\varepsilon>0$ are arbitrary,
\[
\bigl((R_d\circ T_d)(x)\bigr)_k
\rightarrow x_k
\quad\text{in }\sigma(F_k,V_k).
\]

Under the isometric identifications
$F_k=V_k^*$ given by
Lemma~\ref{lem:compression-duality},
all assumptions of Lemma~\ref{lem:local-maps}
are now satisfied.
Hence, by Lemma~\ref{lem:local-maps}, we have that $G$ has the Grothendieck property.
\end{proof}

It is noteworthy that the result of Lemma \ref{lem:square-product} 
may fail when $L_2$ is replaced by an   arbitrary reflexive symmetric space, i.e.,   $F_k$ is the intersection of a semifinite von Neumann algebra and an arbitrary reflexive symmetric space.

\begin{remark}
There is a reflexive fully symmetric Banach function space $E$
on $(0,\infty)$ such that the $\ell_\infty$-sum of countably many
copies of $L_\infty(0,\infty)\cap E$ is not Grothendieck.

Take $\mathcal M=L_\infty(0,\infty)$ with the integral trace, and
let $E:=\Gamma_{2,w}(0,\infty)$ be the Lorentz Gamma space with norm
\[
\|u\|_E
=\left(\int_0^\infty [u^{**}(t)]^2 w(t)\,dt\right)^{1/2},
\qquad
u^{**}(t)=\frac1t\int_0^t u^*(s)\,ds, \quad u\in E,
\]
where
\[
w(t)=
\begin{cases}
\dfrac{1}{t(1-\log t)^2}, & 0<t<1,\\[5pt]
\dfrac{t}{(1+\log t)^2}, & t\ge1.
\end{cases}
\]
Let $\ell_1^n:=\mathop{\bigoplus}\limits_{1\le k\le n}^{\ell_1}
	 \mathbb C$.
Note that the space $E$ is a reflexive fully symmetric space by \cite[Theorem~1.1]{KS21}.
Write $\varphi(t)=\|\chi_{(0,t)}\|_E$. 
By
\cite[Lemma~3.1 and formulas~(4.5), (4.9)]{KS21},
\[
\varphi(t)\longrightarrow\infty,
\qquad
\frac{\varphi(nt)}{\varphi(t)}\longrightarrow n,
\quad t\to\infty
\]
for each fixed $n\ge1$.

Equip $F:=\mathcal M\cap E$ with the norm
$\|u\|_F=\max\{\|u\|_\infty,\|u\|_E\}$.
For each $n\ge1$, choose $t_n>0$ such that
\begin{align} \label{varphi(tn)}
\varphi(t_n)\ge1,\qquad
\varphi(nt_n)\ge \tfrac n2\,\varphi(t_n).
\end{align}
For $1\le i\le n$, let $I_{n,i}:=((i-1)t_n,it_n)$. Define
$J_n:\ell_1^n\to F$ and $Q_n:F\to\ell_1^n$ by
\[
J_n b=\frac{1}{\varphi(t_n)}
       \sum_{i=1}^n b_i\chi_{I_{n,i}}, \quad b=(b_i)_{i=1}^n\in \ell_1^n,
\qquad
(Q_nu)_i=\frac{\varphi(t_n)}{t_n}\int_{I_{n,i}}u(s)\,ds, \quad u\in F.
\]
Since $\chi_{I_{n,i}}/\varphi(t_n)$ has $\norm{\cdot}_F$-norm $1$ (by the definition of $\varphi$ and $\varphi(t_n)\ge1$), it follows that
$\|J_n\|\le1$ for any $n\ge1$. 
Note that 
$Q_nJ_n=\bf 1_{\ell_1^n}$.

Set $A_n:=(0,nt_n)$. 
For any $u\in F$, the function
\[
v:=\left(\frac1{nt_n}\int_{A_n}|u(s)|\,ds\right)\chi_{A_n}\stackrel{\tiny\mbox{\cite[Lem. 3.6.2]{LSZ}}}{\prec\prec} u.
\]
Thus, by the full symmetry of $E$, we have
 $\|v\|_E\le\|u\|_E$, and hence
\[
\|Q_nu\|_1
\le \frac{\varphi(t_n)}{t_n}\int_{A_n}|u(s)|\,ds
= \frac{n\varphi(t_n)}{\varphi(nt_n)}\|v\|_E
\stackrel{\eqref{varphi(tn)}}{\le} 2\|u\|_F.
\]

Since $\sup_n\|Q_n\|\le2$ and $\sup_n\|J_n\|\le1$, it follows from
  $Q_nJ_n=\bf 1_{\ell_1^n}$ that
the
coordinatewise maps define a bounded surjection
\[
Q:\ell_\infty(F)\longrightarrow
\mathop{\bigoplus}\limits_{n\ge1}^{\ell_\infty}
	 \ell_1^n,
\qquad Q((u_n)_n)=(Q_nu_n)_n, \quad (u_n)_n\in \ell_\infty(F).
\]
By \cite[p. 303]{Johnson72}, $\mathop{\bigoplus}\limits_{n\ge1}^{\ell_\infty}
	 \ell_1^n$ is not Grothendieck.
By Proposition \ref{surjective image}, we have that 
$\ell_\infty(F)$ is not Grothendieck.
\end{remark}

\subsection{The Grothendieck property of
$\ell_\infty$-sums of  noncommutative
$L_\infty\cap L_r$ spaces for $1<r<2$} 

We begin with an elementary estimate for the modulus of
convexity of a quotient space, obtained by choosing representatives
with norms arbitrarily close to the corresponding quotient norms (see Lemma \ref{quotient} below), which uses the following lemma.

\begin{lemma}\label{lem:quotient-convexity}\cite[Section~4.1, p.~86]{Pisier11}
	Let $X$ be a uniformly convex normed space,
	and let $Z$ be a  subspace of $X$.
	Then,
	\[
	\delta_{X/Z}(\varepsilon)
	\ge\delta_X(\varepsilon/2),
	\qquad 0<\varepsilon\le2,
	\] 
	where $\delta_{X/Z}$ is the modulus of convexity of $\left( X/Z, \norm{\cdot}_{X/Z} \right)$.
\end{lemma}

 	Let $1<r<2$, and let $(\M,\tau)$ be a semifinite
		von Neumann algebra. 	For $i\in\{r,2\}$, equip
$F_i:=\M\cap L_i(\M,\tau)$ with the intersection norm $\norm{\cdot}_{F_i}$.
Set
\begin{align}\label{defYclosed}
A:=F_2\bigoplus^{\ell_\infty}L_r(\M,\tau),
\qquad
Y:=\{(x,x):x\in F_r\}.
\end{align}
    Fix $x\in F_r$.
    Applying \cite[Theorem 7.13.2]{DPS} to $V:\mathbb C\to \cM, ~V(\lambda)=\lambda x$, $\sigma(t)=t, ~t\ge 0$, $p_0=p_1=2$, $q_0=\infty$, $q_1=r$ and $\theta=\frac r2$ (see also \cite[p.~1466, (2.1)]{PX03}), 
we have that
\begin{align}\label{L2Frnorm}
\norm{x}_{L_2(\M,\tau)}
\le
\norm{x}_{\M}^{1-r/2}
\norm{x}_{L_r(\M,\tau)}^{r/2}
\le\norm{x}_{F_r},
\qquad x\in F_r.
\end{align}
Hence, $F_r\subseteq F_2$ and
$\norm{x}_{F_2}\le\norm{x}_{F_r}$ for every $x\in F_r$.
By  \eqref{diagonal map},
the map $x\mapsto(x,x)$ is an
isometric isomorphism from $F_r$ onto  the closed subspace $Y$ of  $A$.

\begin{lemma} \label{quotient}

    The quotient Banach  space $\left( A/Y, \norm{\cdot}_{A/Y}
    \right)$ admits an  equivalent norm $\nu(\cdot)$ satisfying
\begin{equation}\label{eq:quotient-equiv}
2^{-1/2}\nu(u)
\le\norm{u}_{A/Y}
\le2^{1/2}\nu(u),
\qquad u\in A/Y,
\end{equation}
and $$\delta_{\left(A/Y,\nu\right)}(\varepsilon) \ge w_r(\frac \varepsilon 2)>0, \qquad 0<\varepsilon\le2,$$
    where $w_r(\frac \varepsilon 2)$ depends only on $r$ and $\varepsilon$.
\end{lemma}
\begin{proof}

 Let $h$ denote the standard $\ell_2$-sum norm on
$L_2(\cM,\tau)\mathop{\bigoplus}\limits^{\ell_2}L_r(\cM,\tau)$, 
and define
\[ 
\nu((a,b)+Y):=\inf_{x\in F_r}h(a-x,b-x), \qquad (a,b)\in A.
\]
    We claim that $\nu(\cdot)$ and $\norm{\cdot}_{A/Y}$ are equivalent. Indeed, for any $(a,b)\in A$, 
\begin{align} \label{one dire}
\nu((a,b)+Y)
&=\inf_{x\in F_r}
  \bigl(\norm{a-x}_{L_2(\cM,\tau)}^2+\norm{b-x}_{L_r(\cM,\tau)}^2\bigr)^{1/2}\notag\\
&\le\inf_{x\in F_r}
  \bigl(\norm{a-x}_{F_2}^2+\norm{b-x}_{L_r(\cM,\tau)}^2\bigr)^{1/2}\notag\\
&\le\sqrt2\,\inf_{x\in F_r}
  \max\left\{\norm{a-x}_{F_2},\norm{b-x}_{L_r(\cM,\tau)}
  \right\}\notag\\
&=\sqrt2\,\norm{(a,b)+Y}_{A/Y}.
\end{align}
    For the inverse estimate, fix $(a,b)\in A=F_2\mathop{\bigoplus}\limits^{\ell_\infty}L_r(\M,\tau)$, where $ F_2=\M\cap L_2(\M,\tau)$.
    
  We first consider the case that $a\ne0$.
   Set  $C:=\norm{a}_{L_2\left(\cM,\tau\right)}<\infty$.  
  Applying the spectral trace formula
\cite[Remark~3.3(3)]{FK86} (see also \cite[Proposition 3.2.8]{DPS}) to $g(s)=s^r\chi_{\left(C,\infty\right)}(s)$ and $T=a$, we obtain
    \begin{align} \label{a1 le a r}
\norm{ae^{|a|}\left(C,\infty\right) }_{L_r\left(\cM,\tau\right)}^r& \notag =\tau\left(\left|ae^{|a|}\left(C,\infty\right)\right|^r\right)\\  & \notag =\int_0^\infty \mu(t;a)^r\chi_{\{\mu(t;a)>C\}}dt\notag\\
&\notag \le C^{r-2}\int_0^\infty\mu(t;a)^2 dt
\\ & \le C^{r}
\overset{C=\norm{a}_{L_2\left(\cM,\tau\right)}}{<}
\infty.
    \end{align}
 Since $\norm{a}_{\cM}<\infty$, it follows  that $ae^{|a|}\left(C,\infty\right)\in \M\cap L_r(\M,\tau)=F_r$, and so
 $$\left( ae^{|a|}\left(C,\infty\right), ae^{|a|}\left(C,\infty\right)\right)
    \overset{\eqref{defYclosed}}{\in}
    Y.$$
    Then,
    $$(a,b)+Y=
    \left(ae^{|a|}[0,C], b-ae^{|a|}\left(C,\infty\right) \right)+Y.$$
    This implies that for every $(a,b)+Y\in A/Y$ with $a\ne0$,
    \begin{eqnarray} \label{arbitrary}
         \norm{(a,b)+Y}_{A/Y} 
       &\le &
\norm{\left(ae^{|a|}[0,C],
        b-ae^{|a|}\left(C,\infty\right)\right)}_A\notag\\
        &=&\max\left\{
        \norm{ae^{|a|}[0,C]}_{F_2},\norm{b-ae^{|a|}\left(C,\infty\right)}_{L_r\left(\cM,\tau\right)}
        \right\}\notag\\
        &\stackrel{
        \eqref{a1 le a r}
        }{\le}&\max\{C,\, \norm{b}_{L_r\left(\cM,\tau\right)}+C\}\notag\\
        &=&\norm{b}_{L_r\left(\cM,\tau\right)}+C\le\sqrt 2 \left(\norm{b}_{L_r\left(\cM,\tau\right)}^2+C^2\right)^{1/2}\notag\\
        &=&\sqrt 2 \left(\norm{b}_{L_r\left(\cM,\tau\right)}^2+\norm{a}_{L_2\left(\cM,\tau\right)}^2\right)^{1/2}
        =\sqrt {2}   h(a,b) .
    \end{eqnarray}
 If $a=0$, then
\begin{align} \label{a=0}
\norm{(0,b)+Y}_{A/Y} & \overset{\eqref{defquotientnorm}}{=}
\inf_{x\in F_r}
\norm{(0-x, b-x)}_{A}
\overset{0\in F_r}{\le} \norm{(0,b)}_A
\notag\\ & \,=\norm{b}_{L_r(\M,\tau)}
=h(0,b)\le \sqrt 2h(0,b).
\end{align}
    For every $x\in F_r$, we have
\[
\norm{(a,b)+Y}_{A/Y}
\overset{(x,x)\overset{\eqref{defYclosed}}{\in}
    Y}{=}
    \norm{(a-x,b-x)+Y}_{A/Y}
\stackrel{\eqref{arbitrary},\eqref{a=0}}{\le}\sqrt2\,h(a-x,b-x).
\]
Taking the infimum over all $x\in F_r$, we obtain that
\[
\norm{(a,b)+Y}_{A/Y}
\le\sqrt2\,\inf_{x\in F_r}h(a-x,b-x)
=\sqrt2\,\nu((a,b)+Y).
\]
This together with \eqref{one dire} proves \eqref{eq:quotient-equiv}, and so, 
 $\nu$ is a complete norm on $A/Y$. 

To see that $Y$ is closed in $(A,h)$\footnote{Here, we do not claim the completeness of $(A,h)$.},
let $(a,b)\in A$ belong to the $h$-closure of $Y$. For every $\varepsilon>0$,
there exists $x\in F_r$ such that
$
h(a-x,b-x)<\varepsilon.
$ 
Consequently,
\[
0\le\nu((a,b)+Y)
=\inf_{x\in F_r}h(a-x,b-x)<\varepsilon.
\]
Since $\varepsilon>0$ is arbitrary, it follows that $\nu((a,b)+Y)=0$.
Since $\nu$ is a norm, it follows that  
$(a,b)\in Y$. 
Therefore, $Y$ is closed in $(A,h)$.

    Applying \cite[Theorem~7.14.4]{DPS} with $p=2$ and $p=r$,
respectively, we obtain
\[
\min\left\{
\delta_{L_2(\cM,\tau)}(\varepsilon),
\delta_{L_r(\cM,\tau)}(\varepsilon)
\right\}
\stackrel{\text{ $1<r<2$}}{\ge}
\frac{r-1}{r}
\left(\frac{\varepsilon}{2}\right)^{r/(r-1)}
>0,
\quad 0<\varepsilon\leq 2,
\]
where the lower bound is independent of the semifinite
von Neumann algebra $(\cM, \tau)$.
  Since $h$ is the standard $\ell_2$-sum norm on
$L_2(\cM,\tau)\mathop{\bigoplus}\limits^{\ell_2}L_r(\cM,\tau)$, it follows from 
\cite[Theorem~3]{D41} that there is
  a function $w_r:(0,2]\to(0,1]$, depending only
on $r$, such that
\begin{align} \label{modu}
\delta_{(A,h)}\left(\frac{\varepsilon}{2}\right)
\ge\delta_h\left(\frac{\varepsilon}{2}\right)
\ge w_r\left(\frac{\varepsilon}{2}\right)>0,
\qquad 0<\varepsilon\le2,
\end{align}
where $\delta_h$ denotes the modulus of convexity
of $L_2(\cM,\tau)\mathop{\bigoplus}\limits^{\ell_2}L_r(\cM,\tau)$
equipped with the norm $h$
and the first inequality follows from the
inclusion of $(A,h)$ as a normed subspace of
$L_2(\cM,\tau)\mathop{\bigoplus}\limits^{\ell_2}L_r(\cM,\tau)$. 
Since $Y=\{(x,x):x\in F_r\}$, it follows from the definition of $\nu$ that
\[
\nu(z+Y)=\inf_{y\in Y}h(z-y),\qquad z\in A.
\]
Thus, $\nu$ is the quotient norm induced by $h$ on $A/Y$.
Applying Lemma~\ref{lem:quotient-convexity} to $(A,h)$
and its closed subspace $Y$, we obtain
\[
\delta_{(A/Y,\nu)}(\varepsilon)
\ge\delta_{(A,h)}\left(\frac{\varepsilon}{2}\right)
\stackrel{\eqref{modu}}{\ge}
 w_r\left(\frac{\varepsilon}{2}\right)>0,
\qquad 0<\varepsilon\le2.
\]
\end{proof}

\begin{prop} \label{longlong}
    Let $1<r<2$.
  For each $k\in\mathbb N$, let $(\cM_k,\tau_k)$
be a semifinite von Neumann algebra and let $b_k>0$.  The space
    \[
G_r:=
\mathop{\bigoplus}\limits_{k\in\mathbb N}^{\ell_\infty}
\left(
\cM_k\cap L_r(\cM_k,\tau_k),
\max\left\{
\norm{x}_{\cM_k},
b_k\norm{x}_{L_r(\cM_k,\tau_k)}
\right\}
\right)
\]
    is Grothendieck.
\end{prop}
\begin{proof}
    It suffices to consider the case when $b_k=1$
for every $k\in\mathbb N$, since \[
\begin{aligned}
\|x\|_{L_r(\cM_k,b_k^r\tau_k)}^r
=(b_k^r\tau_k)(|x|^r)=b_k^r\|x\|_{L_r(\cM_k,\tau_k)}^r.
\end{aligned}
\]
    For $i=r,2$, let
    $F_{i,k}:=\cM_k\cap L_i\left(\cM_k, \tau_k\right)$
    equipped with the intersection norm.
    Let \[ A_k:=F_{2,k}\mathop{\bigoplus}\limits^{\ell_\infty }L_r\left(\cM_k, \tau_k\right)\quad \text{  and }   \quad
Y_k:=\{(x,x):x\in F_{r,k}\}.
\]
By the same argument as used in the proof of  Lemma~\ref{quotient} and \eqref{L2Frnorm},  $Y_k$ is closed in $A_k$.
Moreover,
the diagonal map $x\mapsto(x,x)$ is an
isometric isomorphism from $F_{r,k}$ onto $Y_k$.

Note that
\begin{align}\label{ellinftydirectsum}
\mathop{\bigoplus}\limits_{k\in\mathbb N}^{\ell_\infty}A_k
\cong
\left(
\mathop{\bigoplus}\limits_{k\in\mathbb N}^{\ell_\infty}
F_{2,k}
\right)
\bigoplus^{\ell_\infty}
\left(
\mathop{\bigoplus}\limits_{k\in\mathbb N}^{\ell_\infty}
L_r(\cM_k,\tau_k)
\right).
\end{align}
    Applying Lemma~\ref{lem:square-product} with
$a_k=1$, we obtain that
$\mathop{\bigoplus}\limits_{k\in\mathbb N}^{\ell_\infty}
F_{2,k} =\mathop{\bigoplus}\limits_{k\in\mathbb N}^{\ell_\infty}
\big( 
\cM_k\cap L_2\left(\cM_k, \tau_k\right)
\big)$
has the Grothendieck property. 
By \cite[Theorem~7.14.4]{DPS}, the spaces
$L_r(\cM_k, \tau_k)$, $k\in\mathbb N$, $1<r<2$, 
admit a common modulus of uniform convexity.
This together with Corollary~\ref{cor:uc-product}  implies that
$
\mathop{\bigoplus}\limits_{k\in\mathbb N}^{\ell_\infty}
L_r(\cM_k,\tau_k)
$
has the Grothendieck property.
   By \eqref{ellinftydirectsum} and Proposition~\ref{surjective image}(1),
$\mathop{\bigoplus}\limits_{k\in\mathbb N}^{\ell_\infty}A_k$
has the Grothendieck property.

  By Lemma~\ref{quotient}, applied to
$(\cM_k,\tau_k)$, each quotient Banach space $\left(A_k/Y_k, \norm{\cdot}_{A_k/Y_k} \right)$ 
admits an equivalent norm $\nu_k$ satisfying
\[
2^{-1/2}\nu_k(u)
\le\norm{u}_{A_k/Y_k}
\le2^{1/2}\nu_k(u),
\qquad u\in A_k/Y_k,\quad k\in\mathbb N,
\]
and
\[
\delta_{(A_k/Y_k,\nu_k)}(\varepsilon)
\ge w_r(\varepsilon/2)>0,
\qquad 0<\varepsilon\le2,\quad k\in\mathbb N.
\]
Since 
$\mathop{\bigoplus}\limits_{k\in\mathbb N}^{\ell_\infty}A_k$
has the Grothendieck property, it follows from 
Lemma~\ref{lem:coordinate-kernel}  that
$ 
Y:=
\mathop{\bigoplus}\limits_{k\in\mathbb N}^{\ell_\infty}Y_k
$ 
has the Grothendieck property.
    Since $F_{r,k}$ and $Y_k$ are isometrically isomorphic, it follows that $G_r=\mathop{\bigoplus}\limits_{k\in\mathbb N}^{\ell_\infty} F_{r,k}$ and $\mathop{\bigoplus}\limits_{k\in\mathbb N}^{\ell_\infty}Y_k$
are isometrically  isomorphic.
    Therefore, by Proposition~\ref{surjective image}(2), $G_r$ is Grothendieck.
    This completes the proof.
\end{proof}


\section{The Grothendieck property of noncommutative
Marcinkiewicz spaces} \label{s5}




\subsection{An embedding of the Lorentz predual} 
Let $\cM$ be a semifinite von Neumann algebra equipped with a faithful normal semifinite trace $\tau$.
If $\gamma =\tau(\mathbf1)<\infty$, then
$\Lambda_\Phi(\cM,\tau)\subseteq L_1(\cM,\tau)$,
and hence
$\Lambda_\Phi(\cM,\tau)\subseteq
(L_1+L_q)(\cM,\tau)$ for every $q>1$.  The following lemma establishes a similar inclusion
when $\gamma=\infty$, under condition \eqref{condition 1} of Theorem \ref{thm:main}.


\begin{lemma} \label{subset}
Assume that $\tau(\mathbf 1)=\infty$.
 Suppose that condition \eqref{condition 1} of Theorem~\ref{thm:main} holds. 
 Let $\alpha\in (1,2]$ be given by Lemma~\ref{order continuous}, and let\[q>\frac{1}{\log_2\alpha}.\]
 Then,
 $
	\Lambda_\Phi(\M,\tau)
	\subseteq (L_1+L_q)(\M,\tau).
$
\end{lemma}
\begin{proof}
    Fix $x\in \Lambda_\Phi\left(\cM,\tau\right)$. 
    It suffices to prove that $\mu(x)\in(L_1+L_q)(0,\infty)$.
   Since
$\Lambda_\Phi(0,\infty)\subseteq
(L_1+L_\infty)(0,\infty)$ \cite[p. 108]{KPS},
it suffices to prove $\mu(x)\chi_{(1,\infty)}\in L_q(0,\infty )$.

By the discussion preceding Lemma~\ref{order continuous} in combination with Condition~\eqref{condition 1}, we have that $\Phi(0_+)=0$.
By 
Lemma~\ref{order continuous}, we have that  for any $t\ge1$, $\Phi(t)\ge
\frac{\Phi(1)}{\alpha}\,t^{\log_2\alpha}$. 
Since 
    \begin{eqnarray*} 
\norm{x}_{\Lambda_\Phi(\cM,\tau)} 
&=&
\norm{\mu(x)}_{\Lambda_\Phi\left(0,\infty\right)}
\stackrel{\eqref{deffLambdaPhi}}{=}
\int_0^\infty\mu(s;x)\Phi_+'(s)ds\\
&\stackrel{\Phi_+'(s)\ge0}{\ge}&\int_0^t\mu(s;x)\Phi_+'(s)ds\nonumber\stackrel{ 0\leq s\leq t}{\ge}
\mu(t;x)\int_0^t\Phi_+'(s)ds\nonumber\\
&\stackrel{ \eqref{intergralphi}}{=}&\mu(t;x)\Phi(t)\nonumber, 
    \end{eqnarray*}
    it follows that  for any $t\ge1$, 
\begin{align}\label{2.3}
\mu(t;x)
\le \frac{\norm{x}_{\Lambda_\Phi(\cM,\tau)}}{\Phi(t)}
\stackrel{\Phi(t)\ge
\frac{\Phi(1)}{\alpha}\,t^{\log_2\alpha}}{\le} 
\frac{\alpha\norm{x}_{\Lambda_\Phi(\cM,\tau)}}{\Phi(1)}
t^{-\log_2\alpha}.
\end{align}
 Since $q\log_2\alpha>1$, it follows that
\begin{eqnarray*}
\int_1^\infty\mu(t;x)^q\,dt
&\stackrel{\eqref{2.3}}{\le}&
\left(
\frac{\alpha\norm{x}_{\Lambda_\Phi(\cM,\tau)}}{\Phi(1)}
\right)^q
\int_1^\infty t^{-q\log_2\alpha}\,dt \\
&=&\left(
\frac{\alpha\norm{x}_{\Lambda_\Phi(\cM,\tau)}}{\Phi(1)}
\right)^q \frac1{q\log_2\alpha-1}
<\infty,
\end{eqnarray*}   
which shows that $\mu(x)\chi_{\left[1,\infty\right)}\in L_q\left(0,\infty\right)$.
This   yields that
$\mu(x)=\mu(x)\chi_{(0,1)}+\mu(x)\chi_{[1,\infty)}\in(L_1+L_q)(0,\infty)$. 
\end{proof}

Our next aim is to construct an isomorphic embedding
of $\Lambda_\Phi(\cM,\tau)$ into a suitable Banach space.
For this purpose, we first recall the definition
of the $K$-functional.
The  {\it $K$-functional  associated with a Banach couple} $(X_1,X_2)$  
is defined by
$$K(t,x; X_1,X_2)=\inf\{\norm{x_1}_{X_1}+t\norm{x_2}_{X_2}:x=x_1+x_2,x_1\in X_1,x_2\in X_2\}$$
for all $x\in X_1+X_2,~ t>0$ (see \cite[Section 7.4]{DPS}). 
When no confusion arises, we write $K(t,x)$ for
$K(t,x;X_1,X_2)$.
For each fixed $t>0$, $K(t,\cdot)$ is a norm on
$X_1+X_2$, equivalent to the norm $K(1,\cdot)
$
\cite[Section~3.1, p.~38]{BL12}. 

Let
\begin{align}\label{defC1}
q':=\frac{q}{q-1}, \quad C_1:=1+\frac{\alpha}
{1-2^{-(\log_2\alpha-1/q)}}, \quad c:=\mathop{\lim}\limits_{t\to\infty}\Phi_+'(t)<\infty.
\end{align}
For all $k\in\mathbb{Z}$, let
\begin{align}\label{defwkphi2k}
  \theta_k:=2^{k/q'},\quad   w_k:=\Phi_+'(2^k)-\Phi_+'(2^{k+1})\geq 0, \quad   I:=\{k\in\mathbb Z:w_k>0\}.
\end{align}
Note that 
\begin{equation}\label{2.2}
c+\sum_{k=m}^\infty w_k
=
\lim_{N\to\infty}
\left(
c+\Phi_+'(2^m)-\Phi_+'(2^{N+1})
\right)
=
\Phi_+'(2^m),
\quad m\in\mathbb Z.
\end{equation}

For each fixed $k$, the Banach function space
$ 
\left( L_1+L_q \right)(0,\infty)
$, equipped with the norm $K(\theta_k, \cdot ;  L_1(0,\infty), L_q(0,\infty)
)$,  is 
rearrangement-invariant.
A direct calculation shows that its fundamental
function \cite[Chapter~II, Definition~5.1]{BS}  is 
\begin{align}\label{deffundamentalf}
\varphi_k(t)
:=K(
\theta_k,\chi_{(0,t)};  L_1(0,\infty), L_q(0,\infty)
)
=\min\{t,\theta_k t^{1/q}\},
\quad t\geq 0,
\end{align}
which is concave and $
\varphi_k(0_+)
=0$. 
Let $E_k$ denote $(L_1+L_q)(\cM,\tau)$ equipped
with the norm
$ 
K(\theta_k,\cdot ;L_1(\cM,\tau),L_q(\cM,\tau)).
$
For all $x\in E_k$, we have  $$
K(\theta_k,x;  L_1(\M,\tau), L_q(\M,\tau))
\stackrel{\text{\cite[Prop.7.4.2]{DPS}}}{=}
K(\theta_k, \mu(x);  L_1(0,\infty), L_q(0,\infty)
).$$

Set
\[
Z:=\mathop{\bigoplus}\limits_{k\in I}^{\ell_1}E_k
   \mathop{\bigoplus}\limits^{\ell_1}L_1(\cM,\tau).
\]
For $x\in\Lambda_\Phi(\cM,\tau)$, let $D:\Lambda_\Phi(\cM,\tau)\to Z$ defined by
\begin{align}\label{defDtoZ}
D(x):=\bigl((w_kx)_{k\in I},cx\bigr).
\end{align}

 
The following estimate, together with
Lemma~\ref{subset}, will be used to prove that
$D$ is well-defined and bounded.

\begin{lemma}\label{le C1}
Let $\alpha, q$ be as in Lemma \ref{subset}.
For every $x\in\Lambda_\Phi(\cM,\tau)$
and $s>0$, we have 
\[
c\,d\left(s;\mu(x)\right)
+\sum_{k\in\mathbb Z}w_k
\varphi_k\left(d\left(s;\mu(x)\right)\right)
\le C_1\Phi\left(d\left(s;\mu(x)\right)\right), 
\] 
 where 
$C_1$ and $\varphi_k$ are defined  in \eqref{defC1} and  \eqref{deffundamentalf}, respectively.
\end{lemma}
\begin{proof}
Fix $x\in\Lambda_\Phi(\cM,\tau)$ and $s>0$.
Since 
$\Phi(\infty)=\infty$ (see Lemma \ref{order continuous}), it follows from \eqref{2.3} that $\mu(t;x)\to0$
as $t\to\infty$. 
Consequently,
$d(s;\mu(x))<\infty$ for every $s>0$ \cite[p.~133]{DPS}. 
If $d(s;\mu(x))=0$, the assertion is immediate.
Otherwise, let $m\in\mathbb Z$ be the least integer
such that $2^m\ge d(s;\mu(x))$.
Since $$
\frac{\theta_k d(s;\mu(x))^{1/q}}{d(s;\mu(x))}
\stackrel{\eqref{defwkphi2k}}{=}\left(\frac{2^k}{d(s;\mu(x))}\right)^{1/q'},$$
it follows that 
\begin{align} \label{kgem}
\varphi_k(d(s;\mu(x)))
&\overset{\eqref{deffundamentalf}}{=}  \min\{d\left(s;\mu(x)\right),\, \theta_k d\left(s;\mu(x)\right)^{\frac1q}\}
\notag\\ &\,\, =
\begin{cases}
d\left(s;\mu(x)\right), & k\ge m,\\
\theta_k d\left(s;\mu(x)\right)^{1/q}, & k<m.
\end{cases}
\end{align}
Observe that 
\begin{eqnarray}\label{2k ge}
&\, & c\,d\left(s;\mu(x)\right)
+
\sum_{k\ge m}w_k\varphi_k(d(s;\mu(x)))\notag
\\ &\notag \stackrel{\eqref{kgem}}{=} & 
c\,d(s;\mu(x))
+\sum_{k\ge m}w_k \cdot  d(s;\mu(x))
\\ &\notag  = &
d(s;\mu(x))
\left(c+\sum_{k\ge m}w_k\right)\stackrel{\eqref{2.2}}{=}
d(s;\mu(x))\Phi_+'(2^m)
\\
&
\stackrel{\text{$2^m\ge d(s;\mu(x))$}}{\le} & 
d(s;\mu(x))\Phi_+'(d(s;\mu(x)))
\le 
\Phi(d(s;\mu(x))),
\end{eqnarray} 
where we used the fact that  
$t\Phi_+'(t)\le\int_0^t\Phi_+'(s)ds=\Phi(t)$.

Since 
\[
\Phi_+'(2^k)-\Phi_+'(2^{k+1})
\stackrel{\eqref{defwkphi2k}}{=}
w_k\le\Phi_+'(2^k)\stackrel{\tiny t\Phi_+'(t)\le\Phi(t)}{\le}2^{-k}\Phi(2^k), 
\]
it follows  that 
\begin{eqnarray}\label{2k <}
&&\sum_{k<m}w_k\varphi_k(d(s;\mu(x)))\notag
\stackrel{\eqref{kgem}}{=}
\sum_{k<m}w_k\theta_k d(s;\mu(x))^{1/q}\notag\\
&\le &
 \sum_{k<m} 2^{-k}\Phi(2^k)\,2^{k/q'}d(s;\mu(x))^{1/q}=
\sum_{k<m}
\Phi(2^k)
\left(\frac{d(s;\mu(x))}{2^k}\right)^{1/q}\notag\\
&
\stackrel{\text{Lem.\ref{order continuous}}}{\le}  & 
\sum_{k<m} 
\left[ \alpha
\left(\frac{2^k}{d(s;\mu(x))}\right)^{\log_2\alpha} \Phi\left(2^k \cdot  \frac{d(s;\mu(x))}{2^k}
\right)
\right]
\left(\frac{d(s;\mu(x))}{2^k}\right)^{1/q}
\notag\\ 
&  = &
\alpha\Phi(d(s;\mu(x)))
\sum_{k<m}
\left(\frac{2^k}{d(s;\mu(x))}\right)^{\log_2\alpha-1/q}
\notag\\ & =& \alpha\Phi\left(d\left(s;\mu(x)\right)\right) 
\left(\frac{2^{m-1}}{d\left(s;\mu(x)\right)}\right)^{\log_2\alpha-1/q}
\sum_{j=0}^\infty2^{-j(\log_2\alpha-1/q)}
 \notag\\
&
\le & 
\frac{\alpha}
{1-2^{-(\log_2\alpha-1/q)}}
\Phi(d(s;\mu(x)))
=
(C_1-1)\Phi(d(s;\mu(x))), 
\end{eqnarray} 
where we used the fact that 
$
2^{m-1}<d(s;\mu(x)) $ and $\log_2\alpha-1/q>0$. 
Combining \eqref{2k ge} and \eqref{2k <}
completes the proof.
\end{proof}

\begin{lemma} \label{long} 
Let $\alpha, q$ be as in Lemma \ref{subset}.
The linear map  $D:\Lambda_\Phi(\cM,\tau)\to Z$ defined by \eqref{defDtoZ} is  bounded, injective, and  its range is closed. 
    More precisely,
    \begin{align} \label{pro}
        \frac12 \norm{x}_{\Lambda_\Phi\left(\cM,\tau\right)}\le \norm{D(x)}_{Z}\le C_1\norm{x}_{\Lambda_\Phi\left(\cM,\tau\right)}, ~\forall x\in \Lambda_\Phi\left(\cM,\tau\right), 
    \end{align} 
    where $C_1$ is defined  in \eqref{defC1}. 
\end{lemma}
\begin{proof}
    Fix $x\in\Lambda_\Phi(\cM,\tau)$.
By Lemma~\ref{subset},
$x\in(L_1+L_q)(\cM,\tau)$.
Let $c$ be the constant defined in \eqref{defC1}.
Moreover, $cx\in L_1(\cM,\tau)$: this is immediate
when $c=0$, while for $c>0$ it follows from 
$c\int_0^\infty\mu(t;x)dt\le \int_0^\infty \mu(t;x)\Phi_+'(t)dt\stackrel{\eqref{deffLambdaPhi}}{=} \norm{x}_{\Lambda_{\Phi}\left(\cM,\tau\right)}$.
By Tonelli's Theorem \cite[Theorem~2.37]{F99}, we have  
   \begin{eqnarray} \label{phik}
\int_0^\infty\mu(t;x)\,d\varphi_k(t)
&=& \notag \int_0^\infty \left( 
\int_0^\infty
\chi_{\{s<\mu(t;x)\}}(s)\,ds
\right)
\,d\varphi_k(t)\notag\\ 
& 
=& \notag \int_0^\infty \left(\int_0^\infty \chi_{\{\mu(t;x)>s\}}(t)d \varphi_k(t)\right)ds\notag\\
&=&
\int_0^\infty
\left(
\int_0^{m\left( \{t:\mu(t;x)>s\}
\right)
} 
d\varphi_k(t)
\right)\,ds\notag\\
&
= & \notag \int_0^\infty
\left(
\int_0^{d(s;\mu(x))}d\varphi_k(t)
\right)\,ds\notag\\
&= & \notag
\int_0^\infty 
\big( 
\varphi_k(d(s;\mu(x)))-
\varphi_k(0)
\big) \,ds
\\
& \overset{\varphi_k(0)=0}{=} &\notag
\int_0^\infty 
\varphi_k(d(s;\mu(x)))  \,ds
\notag\\
& \stackrel{\eqref{deffundamentalf}}{=} & 
\int_0^\infty
\min\left\{
d(s;\mu(x)),
\theta_k d(s;\mu(x))^{1/q}
\right\}\,ds.
   \end{eqnarray} 
Thus, 
\begin{eqnarray*}
K(\theta_k,x)
&=&K(\theta_k,\mu(x))  
\stackrel{\text{\cite[Chapter~II, Thm.~5.13]{BS}}}{\le}
\norm{\mu(x)}_{\Lambda_{\varphi_k}(0,\infty)}
\\
&=&\int_0^\infty \mu(t;x)\,d\varphi_k(t) 
\\
&\stackrel{\eqref{phik}}{=}&
\int_0^\infty
\min\left\{
d(s;\mu(x)),
\theta_k d(s;\mu(x))^{1/q}
\right\}\,ds,
\end{eqnarray*} 
where $\Lambda_{\varphi_k}(0,\infty)$ is  the
Lorentz function space associated with $\varphi_k$ \cite[Chapter~II, Def.~5.12]{BS}, with   norm
$ 
\norm{f}_{\Lambda_{\varphi_k}(0,\infty)}
:=
\int_0^\infty\mu(t;f)\,d\varphi_k(t).
$ 
Hence, 
\begin{eqnarray*}
   \norm{Dx}_Z 
&=&\sum_{k\in I}K(\theta_k,w_kx)
+\norm{cx}_{L_1(\cM,\tau)}\\
&  =
&
\sum_{k\in I}w_k K(\theta_k,x)
+\norm{cx}_{L_1(\cM,\tau)}\\
& \overset{\text{Tonelli's Thm.}}{\leq} &
\int_0^\infty \sum_{k\in I}w_k \cdot 
\min\left\{
d(s;\mu(x)),
\theta_k d(s;\mu(x))^{1/q} 
\right\} \,ds \\ 
&& 
+\int_0^\infty 
c\,d(s;\mu(x))
\,ds\\
& 
\stackrel{\eqref{deffundamentalf},~\text{Lem.~\ref{le C1}}}{\le} &
C_1\int_0^\infty\Phi(d(s;\mu(x)))\,ds
\stackrel{\substack{d\mu=\Phi_+'(t)dt\\\tiny\mbox{\cite[Prop. 1.8]{BS}}}}{=}C_1\int_0^\infty\mu(t;x)\Phi_+'(t)\,dt
\\
& \stackrel{\eqref{deffLambdaPhi}}{=} &
C_1\norm{x}_{\Lambda_\Phi(\cM,\tau)}
<\infty,
\end{eqnarray*}
which together with the fact that $x\in \Lambda_\Phi(\cM,\tau)\subset (L_1+L_q)(\cM,\tau)$
implies that $((w_kx)_{k\in I},cx)\in Z$, i.e., 
 $D$ is well-defined. 
Its linearity follows directly from the definition.

We next prove $ \norm{D(x)}_{Z}\geq \frac12 \norm{x}_{\Lambda_\Phi\left(\cM,\tau\right)}$ for all $x\in\Lambda_\Phi(\cM,\tau)$.
For each $k\in I$ and every decomposition
$\mu(x)=a+b$, where $a\in L_1(0,\infty)$ and
$b\in L_q(0,\infty)$, by H\"older's inequality, we have that
\[
\begin{aligned}
\int_0^{2^k}\mu(t;x)\,dt
&\le \int_0^{2^k}|a(t)|\,dt
+\int_0^{2^k}|b(t)|\,dt\\
&\le \norm{a}_{L_1(0,\infty)}
+\left(\int_0^{2^k}|b(t)|^q\,dt\right)^{1/q}
 \left(\int_0^{2^k}{\bf 1}^{q'}\,dt\right)^{1/q'}\\
&\le \norm{a}_{L_1(0,\infty)}
+(2^k)^{1/q'}\norm{b}_{L_q(0,\infty)}
\\ & \stackrel{\eqref{defwkphi2k}}{=}\norm{a}_{L_1(0,\infty)}
+\theta_k\norm{b}_{L_q(0,\infty)}.
\end{aligned}
\]
Taking the infimum over all such decompositions of $\mu(x)$,
we obtain
\begin{align}\label{2.1}
\int_0^{2^k}\mu(t;x)\,dt
\le K(\theta_k,\mu(x))
=K(\theta_k,x).
\end{align}
Fix $t>0$ and choose $m\in\mathbb Z$ such that
$2^m\le t<2^{m+1}$.
Since $w_k=0$ for $k\notin I$ (by the definition of $I$, see \eqref{defwkphi2k}), it follows from the monotonicity
of $\Phi_+'$ that 
\begin{eqnarray} \label{t<2 k}
c+\sum_{\substack{k\in I\\t<2^k}}w_k
= c+\sum_{k=m+1}^\infty w_k
\stackrel{\eqref{2.2}}{=}
\Phi_+'(2^{m+1})
\stackrel{\text{$
2^{m+1}\le2t
$}}{\geq}
\Phi_+'(2t).
\end{eqnarray}
Consequently, by Tonelli's theorem,
\begin{align*} 
\norm{D(x)}_Z
&\notag \,\, =\sum_{k\in I}w_kK(\theta_k,x)
  +\norm{cx}_{L_1(\cM,\tau)}
 \\ &  \notag \stackrel{\eqref{2.1}}{\ge}
\sum_{k\in I}w_k\int_0^{2^k}\mu(t;x)\,dt
+\int_0^\infty c\,\mu(t;x)\,dt\\
&\notag \,\, =
\int_0^\infty
\mu(t;x)
\left(
c+\sum_{\substack{k\in I\\t<2^k}}w_k
\right)\,dt
\stackrel{\eqref{t<2 k}}{\ge}
\int_0^\infty\mu(t;x)\Phi_+'(2t)\,dt\\
&\notag \,\, =\frac12
\int_0^\infty\mu(s/2;x)\Phi_+'(s)\,ds
\ge\frac12
\int_0^\infty\mu(s;x)\Phi_+'(s)\,ds
\\ & \,\, =\frac12\norm{x}_{\Lambda_\Phi(\cM,\tau)}.
\end{align*}
This implies that $D$ is an isomorphism
from $\Lambda_\Phi(\cM,\tau)$ onto
$D(\Lambda_\Phi(\cM,\tau))$, which completes the proof.
\end{proof}

\subsection{Proof of Theorem \ref{thm:main}} 

We are now ready to prove Theorem~\ref{thm:main}.


\begin{proof}[Proof of Theorem~\ref{thm:main}]
(1) First we assume that $\gamma= \tau(\mathbf{1})=\infty$. 
Under condition \eqref{condition 1}, we have $\Phi(0_+)=0$.
Moreover, $\Phi(\infty)=\infty$ by
Lemma~\ref{order continuous}. Hence, 
$\Lambda_\Phi(\cM,\tau)$ has order continuous norm
by \cite[Theorem~6.3.8]{DPS}. 
Therefore,   
$$\Lambda_\Phi(\cM,\tau)^*
\stackrel{
\text{
\cite[Proposition 5.3.2]{DPS}
}
}{=}
\Lambda_\Phi(\cM,\tau)^{\times}
\stackrel{
\text{\cite[Theorem 6.4.9]{DPS}}
}{=}
M_\Phi\left(\cM,\tau\right)$$ with equality of  norms. 

Let $\alpha$ be given by
Lemma~\ref{order continuous}, and choose
\[
q>\max\left\{2,\frac{1}{\log_2\alpha}\right\}.
\]
Since $q>2$, its conjugate exponent $q':=q/(q-1)$
satisfies $1<q'<2$. 
For this choice of $q$, let $C_1$,  $\theta_k$, 
$I$, $Z$, $\left(E_k, K(\theta_k,\cdot)\right)$, and $D$
be as  defined before Lemma \ref{le C1}. 
By \cite[Proposition 6.1.9]{DPS}
and  \cite[Theorem 2.7.1]{BL12}, we have that  for every $k\in I$, 
\begin{align*}
\left(E_k,K(\theta_k,\cdot)\right)^* &=\left(\left(L_1+L_q\right)\left(\cM,\tau\right),K(\theta_k,\cdot)\right)^*\\ 
& =\left(L_{q'}\left(\cM,\tau\right)\cap\cM,\norm{\cdot}_{L_{q'}\left(\cM,\tau\right)\cap\cM}\right),
\end{align*}
where $\norm{\cdot}_{L_{q'}\left(\cM,\tau\right)\cap\cM}=\max\{\norm{\cdot}_\cM,\theta_k^{-1}\norm{\cdot}_{L_{q'}\left(\cM,\tau\right)}\}$.
By \cite[p.61]{L04},
we have that
\begin{align*}
Z^*=
\left(\mathop{\bigoplus}\limits_{k\in I}^{\ell_1}E_k\mathop{\bigoplus}\limits^{\ell_1} L_1\left(\cM,\tau\right)\right)^*=\mathop{\bigoplus}\limits_{k\in I}^{\ell_\infty}E_k^*\mathop{\bigoplus}\limits^{\ell_\infty} \cM
\end{align*}
Applying Proposition  \ref{longlong} to $r=q'$ and $b_k=\theta_k^{-1}$, we obtain that $\mathop{\bigoplus}\limits_{k\in I}^{\ell_\infty}E_k^*$ is Grothendieck. By \cite[Corollary~7]{P94}, the von Neumann algebra
$\cM$ has the Grothendieck property. 
By  Proposition ~\ref{surjective image}(1),   $Z^*$ is Grothendieck.

To prove that $M_\Phi(\cM,\tau)$ is Grothendieck,
it suffices, by Proposition~\ref{surjective image}(2),
to show that the adjoint operator 
\[
D^*:Z^*\rightarrow
\Lambda_\Phi(\cM,\tau)^*=M_\Phi(\cM,\tau)
\]
is surjective.
Indeed, fix $x\in M_\Phi\left(\cM,\tau\right)$.
Note that $D\left(\Lambda_\Phi\left(\cM,\tau\right)\right)$ is a (closed) subspace of $Z$ (see Lemma \ref{long}).
Let $D(y)\mapsto \tau(xy)$ be a linear functional on $D\left(\Lambda_\Phi\left(\cM,\tau\right)\right)\subseteq Z$. 
Then, for all $y\in \Lambda_\Phi\left(\cM,\tau\right)$,
$$|\tau(xy)|\le \norm{x}_{M_\Phi\left(\cM,\tau\right)}\norm{y}_{\Lambda_\Phi\left(\cM,\tau\right)}\stackrel{\eqref{pro}}{\le} 2\norm{x}_{M_\Phi\left(\cM,\tau\right)}\norm{D(y)}_{Z}.$$
By the Hahn--Banach theorem, there exists a bounded linear functional $z^*\in Z^*$ such that  
\[
z^*(D(y))=\tau(xy),
\qquad y\in\Lambda_\Phi(\cM,\tau).
\]
In particular,  
$$\left( D^* (z^*)\right)(y)=z^*(D(y))=\tau(xy),  \quad   y\in \Lambda_\Phi\left(\cM,\tau\right),$$
which implies that  $D^*(z^*)=x$.
This implies that $D^*$ is surjective.

(2) Suppose that $\gamma=\tau(\mathbf1)<\infty$.
By Lemma~\ref{lem:growth}, there exists an increasing
concave function $\Psi:[0,\infty)\to[0,\infty)$
with $\Psi(0)=0$ satisfying \[\liminf_{t\downarrow0}\frac{\Psi(2t)}{\Psi(t)}>1,
\qquad
\liminf_{t\to\infty}\frac{\Psi(2t)}{\Psi(t)}>1
\] and
\[
\Phi(t)\le\Psi(t)\le2\Phi(t),
\qquad 0<t<\gamma.
\]

Equip the von Neumann algebra 
\[
\cN:=\cM\mathop{\bigoplus}^{\ell_\infty}
L_\infty(0,\infty)
\]
 with the semifinite faithful normal trace
\begin{align*} 
\widetilde\tau(x,f)
:=\tau(x)+\int_0^\infty f(t)\,dt,
\qquad (x,f)\in\cN_+.
 \end{align*}
By the above result, we obtain that $M_\Psi(\cN,\tilde{\tau})$ is Grothendieck. 
Observing  that $M_\Phi(\M,\tau)$ is a complemented subspace of $M_\Psi(\cN,\tilde{\tau})$, we obtain that  $M_\Phi(\cM,\tau)$ is Grothendieck. 
\end{proof}

\begin{proof}[Proof of Corollary~\ref{cor:characterization}]
Sufficiency follows from Theorem~\ref{thm:main}.
The necessity is a consequence of  Lemma~\ref{lem:expectation-quotient}  and \cite[Theorem~5.6]{PS20}, by 
 observing  that $\cM$ contains an atomless abelian reduced von Neumann algebra $\cA$ which is $*$-isomorphic to $L_\infty (0,\tau({\bf 1}))$ \cite[Lemma~4.17]{MNS} (see also \cite[Lemma~1.3]{CKS}; the existence of conditional expectation follows from \cite[Proposition~7.1.1]{DPS}) and $M_\Phi(\cM,\tau)$ is a fully symmetric space. 
\end{proof}

\section{The atomic case: counterexamples and two necessary
growth conditions}\label{sec:atomic-traces}
We begin with two examples of 
atomic
von Neumann algebras for which $M_\Phi(\M,\tau)$ is Grothendieck,
but conditions \eqref{condition 1} and \eqref{condition 2} of
Theorem~\ref{thm:main} fail in the infinite and finite trace
cases, respectively.

\begin{example}\label{atomic}
  (1)  Let  
  $\Phi(0)=0$ and $\Phi(t)=1+t$ for $t>0$.  Then,  $\Phi$ is increasing,  concave, and 
    $$\liminf_{t\downarrow 0}\frac{\Phi(2t)}{\Phi(t)}=\liminf_{t\downarrow 0} \frac{1+2t}{1+t}=1, \quad 
\lim_{t\to\infty}\frac{\Phi(2t)}{\Phi(t)}=2.$$ 
    In particular, $\Phi$ does not satisfy the condition \eqref{condition 1} of  Theorem \ref{thm:main}, and so,  $M_\Phi(0,\infty)$ fails the Grothendieck
property by \cite[Theorem~5.6]{PS20}.
However, if 
 $\cM=\ell_\infty$ is equipped with the counting trace, then $M_\Phi(\cM,\tau) =\ell_\infty$ with equivalent norms, which has the Grothendieck property. 

(2)  If $\cM$ is a finite-dimensional von Neumann algebra, then  
 for any nonzero increasing concave function $\Phi:[0,\infty)\to [0,\infty)$, $M_\Phi\left(\cM,\tau\right)$ is finite-dimensional, and so $M_\Phi\left(\cM,\tau\right)$ is Grothendieck.
\end{example}

Let $(\mathcal M,\tau)$ be a nonzero semifinite von Neumann algebra.
The next two propositions give two different assumptions on projection
traces under which the growth conditions at zero and at infinity,
respectively, are necessary for $M_\Phi(\mathcal M,\tau)$ to have
the Grothendieck property.


	Let $E$ be a fully symmetric  Banach function space on $[0,\infty)$.
	A positive linear functional $\varphi\in E^*$ is called
	\emph{fully symmetric} if
	\[
	f\prec\prec g \quad\Longrightarrow\quad \varphi(f)\le\varphi(g)
	\qquad \text{for all } f,g\in E_+;
	\]
	see, e.g., \cite[Definition~2.1]{DPSS98}, \cite{LSZ,LSZ21,LMSZ23}.  
	For $s>0$, define the dilation operator $D_s$ by
	\[
	(D_s f)(t):=f(t/s),\qquad t>0.
	\]
	Note that the operator $D_s$ is bounded on $E$
	(see, e.g., \cite[p.~50]{DPSS98}, \cite[p. 392]{DPS}).
	\begin{prop} \label{prop:atomic-zero}
    Let $\Phi:[0,\tau(\mathbf 1))\to[0,\infty)$ be a nonzero nondecreasing concave function with $\Phi(0)=0$.
	Suppose that $M_\Phi(\M,\tau)$ has the Grothendieck property.
    If for some $c>0$, there is a family $\{e_{n,k}, ~k\ge 1, 1\le n\le k\}$ of mutually orthogonal projections in $\cM$  
    with
	\begin{equation}\label{eq:triangular-traces}
		\tau(e_{n,k})=c2^{-k},
		\qquad k\ge1,\quad 1\le n\le k,
	\end{equation}
	then
	\begin{equation*}
		\liminf_{t\downarrow0}\frac{\Phi(2t)}{\Phi(t)}>1.
	\end{equation*}
	\end{prop}

	\begin{proof}
	Set
	$ p:=\mathop{\sum}\limits_{k\ge1}
    \mathop{\sum}\limits_{n=1}^k  e_{n,k},
$
	where the sum converges in the strong operator topology.
	Let 
\[
	\cA:=
	\mathop{\bigoplus}\limits_{\substack{k\ge1\\1\le n\le k}}^{\ell_\infty}
	\mathbb C e_{n,k},
	\]
	with unit \(p\). 
	By normality of $\tau$, we have that
	\[
	\tau(p)
	\overset{\eqref{eq:triangular-traces}}{=}
	c\sum_{k\ge1}k2^{-k}=2c<\infty.
	\]
	Thus, $\tau|_\cA$ is finite, and so	
	$M_\Phi(\cA,\tau|_{\cA})$ is Grothendieck by
	Lemma~\ref{lem:expectation-quotient}. 

    We claim that $\Phi(0_+)=0$.
	Indeed, assume by contradiction that $\Phi(0_+)>0$.
    Then, $0<\Phi(0_+)\le\Phi(t)\le\Phi(2c)$
	for $0<t<2c$. 
    Since $\tau(p)=2c$, we have that for all $x\in M_\Phi\left(\cA,\tau|_\cA\right)$,
	\begin{align*}
		 \frac{\norm{x}_{L_1(\cA,\tau|_{\cA})}}{\Phi(2c)}
		 &=	\frac{\int_0^{2c}\mu(s;x)\,ds}{\Phi(2c)}
	   =\sup_{0<t<2c}
		 \frac{\int_0^t\mu(s;x)\,ds}{\Phi(2c)}\\
		 &\le \sup_{0<t<2c}
	   \frac{\int_0^t\mu(s;x)\,ds}{\Phi(t)}
	   \overset{\eqref{corneranormMphi}}{=}
		 	\norm{x}_{M_\Phi(\cA,\tau|_{\cA})}\\
		 	&\le \sup_{0<t<2c}
		 	\frac{\int_0^t\mu(s;x)\,ds}{\Phi(0_+)}
		 	=\frac{\norm{x}_{L_1(\cA,\tau|_{\cA})}}{\Phi(0_+)}.
		 \end{align*}
		Thus, \( M_\Phi(\cA,\tau|_{\cA})=L_1(\cA,\tau|_{\cA})\) as sets, with equivalent norms. 
		Since \(\cA\) has countably infinitely many atoms, it follows that  $L_1(\cA,\tau|_{\cA})$ fails to have the Grothendieck
		property \cite[discussion following Proposition~2.3]{Wnuk}. 
        This proves the claim.

         Below, we will prove that $	\mathop{\liminf}\limits_{t\downarrow0}\frac{\Phi(2t)}{\Phi(t)}>1$.
         Assume by contradiction that
     $	\mathop{\liminf}\limits_{t\downarrow0}\frac{\Phi(2t)}{\Phi(t)}=1$.
	A combination of $\Phi(0_+)=0$, $\mathop{\liminf}\limits_{t\downarrow0}\frac{\Phi(2t)}{\Phi(t)}=1$ and \cite[Theorem~3.4(ii) and Definition 2.7]{DPSS98} yields a nonzero positive fully symmetric functional $\varphi\in M_\Phi(0,\infty)^*$ that vanishes on
	$M_\Phi(0,\infty)\cap L_\infty(0,\infty)$.
	Since $\varphi$ is nonzero and positive, it follows that  there exists
	$h\in M_\Phi(0,\infty)_+$ with $\varphi(h)>0$.
 Since $\mu(h)\prec\prec h$, $h\prec\prec \mu(h)$ and $\varphi$ is fully symmetric, it follows that 
 $\varphi\left(\mu(h)\right)=\varphi(h)>0$.
	
 Note that for $t\ge c$, 
 \[
 \mu(t;h)\le \frac1c\int_0^c \mu(s;h)\,ds
 \le \frac{\Phi(c)}{c}\norm{\mu(h)}_{M_\Phi(0,\infty)}.
 \]
 This implies that $\mu(h)\chi_{[c,\infty)}\in M_\Phi(0,\infty)\cap L_\infty(0,\infty)$,
 and hence
 $\varphi(\mu(h)\chi_{[c,\infty)})=0$. 
 Then, 
 $\varphi(\mu(h)\chi_{(0,c)})=\varphi(\mu(h))>0$.

  Set
		 \[
		 g:=\sum_{k\ge1} \mu\left(c2^{1-k};h\right)\chi_{[c2^{-k},\,c2^{1-k})}.
		 \]
	For $t\in[c2^{-k},c2^{1-k})$, we have $t<c2^{1-k}\le 2t$.
		  Since $\mu(h)$ is decreasing, it follows that 
		  \[
		  \mu(2t;h)\le \mu\left(c2^{1-k};h\right)=g(t)\le \mu(t;h).
		  \]
        This together with $g(t)=0$ for all $t\ge c$ yields that for all $t>0$,
	   $$\mu(2t;h)\chi_{(0,c)}(2t)\le \mu(2t;h)\chi_{(0,c)}(t)\le g(t)\chi_{(0,c)}(t)=g(t)\le \mu(t;h)\chi_{(0,c)}(t),$$
	   which implies that $g\in M_\Phi(0,\infty)$.
      Then,
		 \begin{eqnarray} \label{>0}
		 \varphi(g)& \ge  & \varphi\left(\mu\left(2\,\cdot;h\right)\chi_{(0,c)}(2\cdot)\right) =
		 \varphi \left(D_{\frac{1}{2}}\left(\mu(h)\chi_{(0,c)}\right)\right)\notag\\
		 &\overset{\text{\cite[Prop.2.3]{DPSS98}}
		 }{=}&
		 \tfrac12\varphi\left(\mu(h)\chi_{(0,c)}\right)\stackrel{\varphi(\mu(h)\chi_{(0,c)})=\varphi(\mu(h))}{>}0.
		 \end{eqnarray}

For each $n\ge1$ and $N\ge n$, set
\[
v_n^{(N)}:=\sum_{k=n}^{N}\mu\left(c2^{1-k};h\right)e_{n,k}.
\]
For any $M>N$ and $\varepsilon>0$, since $e_{n,k}$ are mutually orthogonal, we have that
	\begin{align*}
    \tau\!\left(e^{|v_n^{(M)}-v_n^{(N)}|}
	(\varepsilon,\infty)\right)
		&=
	\sum_{\substack{N<k\le M\\ \mu\left(c2^{1-k};h\right)>\varepsilon}}\tau(e_{n,k}) 
	\overset{\eqref{eq:triangular-traces}}{=}
	\sum_{\substack{N<k\le M\\ \mu\left(c2^{1-k};h\right)>\varepsilon}}c2^{-k}
	\\ & \le c\sum_{k>N}2^{-k}
	=c2^{-N}\rightarrow0.
		\end{align*}
Thus,  $(v_n^{(N)})_{N\ge n}$ is Cauchy in the measure topology \cite[p.4]{BS19}.
By the completeness of $S\left(\cA,\tau|_{\cA}\right)$, it converges in measure to 
\[
v_n:=\sum_{k\ge n}\mu\left(c2^{1-k};h\right)e_{n,k}\in S\left(\cA,\tau|_{\cA}\right).
\]

     Observe that 
		 \begin{align}\label{muvngchi}
		 \mu(v_n)=g\chi_{(0,c2^{1-n})}.
		 \end{align}
		 In particular, 
    $\norm{v_n}_{M_\Phi(\cA,\tau|_{\cA})}\le \norm{g}_{M_\Phi(0,\infty)}$ 
    and $v_nv_m=v_mv_n=0$ when $m\ne n$.
    Since $s(v_n)\le\sum_{k\ge n}e_{n,k}$, it follows that the support projections $s(v_n)$ of $v_n$ are pairwise orthogonal. 

	By \cite[Theorem~4.2]{DPSS98},
		the map
		\begin{align} \label{extension}
		x\mapsto\varphi(\mu(x)),
		\qquad x\in M_\Phi(\cA,\tau|_{\cA})_+,
		\end{align}
		extends uniquely to a bounded positive linear functional
		$\psi\in M_\Phi(\cA,\tau|_{\cA})^*$.
        Noting $g\in M_\Phi\left(0,\infty\right)$, we have that 
        $$g\chi_{[c2^{1-n},\infty)}
=\sum_{k=1}^{n-1} \mu\left(c2^{1-k};h\right)\,
\chi_{[c2^{-k},\,c2^{1-k})}\in M_\Phi(0,\infty)\cap L_\infty(0,\infty),$$
which implies that $\varphi(g\chi_{[c2^{1-n},\infty)})=0$. 
Hence,
\begin{eqnarray} \label{psi(vn)>0}
	\psi(v_n)
	&\stackrel{\eqref{extension}}{=}&\varphi(\mu(v_n))
	\overset{\eqref{muvngchi}}{=}
	\varphi\bigl(g\chi_{(0,c2^{1-n})}\bigr)\notag\\ 
    & =&\varphi(g)-\varphi\bigl(g\chi_{[c2^{1-n},\infty)}\bigr)
	=\varphi(g)\stackrel{\eqref{>0}}{>}0,\quad n\ge1.
\end{eqnarray}

 Define
	\begin{align}\label{defQatomic}
	Q\colon M_\Phi(\cA,\tau|_{\cA})\rightarrow \ell_1,
	\qquad Q(x)=\bigl(\psi\left(xs(v_n)\right)\bigr)_{n\ge1}.
	\end{align}
    Since $\cA$ is abelian, it follows that $|xs(v_n)|=|x|s(v_n)$ for all $x\in M_\Phi\left(\cA,\tau|_\cA\right)$. 
    By positivity of $\psi$
	and $\sum_{n=1}^N s(v_n)\le p$, we have 
	\[
	\sum_{n=1}^N |\psi(xs(v_n))|
	\le \psi\!\left(|x|\sum_{n=1}^N s(v_n)\right)
	\le \psi(|x|)
	\le \norm{\psi}\norm{x}_{M_\Phi(\cA,\tau|_{\cA})}.
	\]
	Letting $N\to\infty$, we obtain that $\bigl(\psi\left(xs(v_n)\right)\bigr)_{n\ge1}\in \ell_1$, that is, $Q$ is well-defined and  is a bounded linear map
	into $\ell_1$.

By $s(v_n)v_n=v_ns(v_n)=v_n$ and $v_ns(v_m)=0$ for $m\ne n$, we have
$Q(v_n/\psi(v_n))$ is the $n$-th unit vector of $\ell_1$.
Since $\norm{v_n}_{M_\Phi(\cA,\tau|_{\cA})}\le \norm{g}_{M_\Phi(0,\infty)}$
and $\psi(v_n)\stackrel{\eqref{psi(vn)>0}}{=}\varphi(g)>0$ for every $n$, it follows that
$(v_n/\psi(v_n))_{n\ge1}$ is norm bounded.
Thus, for any $b=(b_n)\in\ell_1$, the series
\[
x=\sum_{n\ge1}b_n\frac{v_n}{\psi(v_n)}
\]
converges in $M_\Phi(\cA,\tau|_{\cA})$.
Moreover, by the continuity of $\psi$, we have that $Q(x)=b$.
Therefore,
$Q$ is surjective. 
Since $M_\Phi\left(\cA,\tau|_{\cA}\right)$ is Grothendieck, it follows from Proposition~\ref{surjective image}(2) that $\ell_1$ is Grothendieck, which contradicts \cite[p. 764]{Wnuk}.
Therefore, $\mathop{\liminf}\limits_{t\downarrow0}\frac{\Phi(2t)}{\Phi(t)}>1$.
	\end{proof}

\begin{prop}\label{prop:atomic-infty} 	
Let $\Phi:[0,\tau(\mathbf 1))\to[0,\infty)$ be a nonzero nondecreasing concave function with $\Phi(0)=0$.
Suppose that $M_\Phi(\M,\tau)$ has the Grothendieck property. 
If for some $a>0$, there are mutually orthogonal projections
      $r_j\in\cM$ with
      $\tau(r_j)=a,$ $ j\ge0,
    $ 
      then \begin{equation*}
      \liminf_{t\to\infty}\frac{\Phi(2t)}{\Phi(t)}>1.
      \end{equation*}
\end{prop}
\begin{proof}
Set
\[
p:=\sum_{k\in\mathbb N}r_k,
\qquad
\cA:=\mathop{\bigoplus}\limits_{k\in\mathbb N}^{\ell_\infty}
\mathbb C r_k\subseteq p\M p,
\]
where the sum defining $p$ converges in the strong operator topology
and $\cA$ has unit $p$. 
Since $\tau(r_k)=a$ for every $k$ and $r_j$ are mutually orthogonal, we have that
$\tau(p)=\infty$. 
By \cite[Remark 1.15.5]{DPS}, we have that
$\tau|_{\cA}$ is semifinite. 
This together with Lemma~\ref{lem:expectation-quotient} yields that $M_\Phi(\cA,\tau|_{\cA})$ is Grothendieck.

We claim that $\Phi(\infty)=\infty$.
Indeed, assume by contradiction that
$\Phi(\infty)<\infty$. 
For every
$x\in M_\Phi(\cA,\tau|_{\cA})$ and $t>0$, we have
\[
\int_0^t\mu(s;x)\,ds
\le \Phi(t)\norm{x}_{M_\Phi(\cA,\tau|_{\cA})}
\le \Phi(\infty)\norm{x}_{M_\Phi(\cA,\tau|_{\cA})}<\infty.
\]
Taking $t\to\infty$, we obtain that
\[
\tau(|x|)
=\int_0^\infty\mu(s;x)\,ds
\le \Phi(\infty)\norm{x}_{M_\Phi(\cA,\tau|_{\cA})}<\infty.
\]
Since $|\tau(x)|\le\tau(|x|)$, it follows that the map $x\in M_\Phi(\cA,\tau|_{\cA})\mapsto \tau(x)$ is a bounded positive linear functional on $M_\Phi(\cA,\tau|_{\cA})$.
Note that
\begin{align} \label{rk/tau(rk)}
\norm{\frac{r_k}{\tau(r_k)}}_{M_\Phi(\cA,\tau|_{\cA})}
\overset{\eqref{eq:projection-norm}}{=}
\frac{1}{\Phi(\tau(r_k))}
= \frac{1}{\Phi(a)}.
\end{align}
Applying the construction in \eqref{defQatomic} to $\psi=\tau$ and $s(v_n)=r_k$,
we obtain that there is a bounded linear mapping
\[
Q:M_\Phi\left(\cA,\tau|_{\cA}\right)\rightarrow\ell_1,
\qquad x\mapsto\bigl(\tau(xr_k)\bigr)_{k\ge0}.
\]
Observe that $Q\left(\frac{r_k}{\tau(r_k)}\right)=e_k$, where  $\{e_n\}_{n=0}^\infty$ be the standard basis of $c_0$.
For any $b=(b_k)_{k=0}^\infty\in \ell_1$, we have
$$\norm{\sum_{k=0}^\infty b_k\frac{r_k}{\tau(r_k)}}_{M_\Phi\left(\cA,\tau|_{\cA}\right)}\le \sum_{k=0}^\infty\norm{b_k\frac{r_k}{\tau(r_k)}}_{M_\Phi\left(\cA,\tau|_{\cA}\right)}\stackrel{\eqref{rk/tau(rk)}}{=}\frac{1}{\Phi(a)}\sum_{k=0}^\infty |b_k|<\infty.$$
Then, by the continuity of $Q$, $Q\left(\sum_{k=0}^\infty b_k\frac{r_k}{\tau(r_k)}\right)\stackrel{Q\left(\frac{r_k}{\tau(r_k)}\right)=e_k}{=}b,$
which shows that $Q$ is surjective.
This contradicts the Grothendieck property of
$M_\Phi(\cA,\tau|_{\cA})$. 
This proves the claim.

Suppose by contradiction now that 
\[
\liminf_{t\to\infty}\frac{\Phi(2t)}{\Phi(t)}=1.
\]
By \cite[Theorem~3.4(i) and Definition~2.7]{DPSS98}, there is a
nonzero positive fully symmetric functional
$\varphi\in M_\Phi(0,\infty)^*$ that vanishes on
$M_\Phi(0,\infty)\cap L_1(0,\infty)$. 
By positivity and symmetry of $\varphi$, we may choose a nonnegative decreasing
$f=\mu(f)\in M_\Phi(0,\infty)\subseteq \left(L_1+L_\infty\right)\left(0,\infty\right)$ with $\varphi(f)>0$. 
Set
\begin{align}\label{defggkatomic}
g_k:=\frac1a\int_{ka}^{(k+1)a}f(s)\,ds<\infty,
\qquad
g:=\sum_{k\ge0}g_k\chi_{[ka,(k+1)a)}.
\end{align}
Since $f$ is decreasing, it follows that  
\begin{align*}
g_{k+1}
&=\frac1a\int_{(k+1)a}^{(k+2)a}f(t)\,dt
=\frac1a\int_{ka}^{(k+1)a}f(s+a)\,ds
\\ & \le\frac1a\int_{ka}^{(k+1)a}f(s)\,ds
=g_k,
\end{align*} 
which implies that 
$(g_k)$ is decreasing. 
By the definition of $g_j$, we have that  for every integer $k\ge1$,
\begin{align} \label{g=f ka}
\int_0^{ka}g(s)\,ds
=\sum_{j=0}^{k-1}a g_j
=\int_0^{ka}f(s)\,ds.
\end{align}
Fix $t>0$ and choose $k\ge0$ such that $ka\le t<(k+1)a$. Then, 
\begin{eqnarray*}
&&\int_0^t g(s)\,ds
\\ &\overset{\eqref{defggkatomic},\eqref{g=f ka}}{=}&
\int_0^{ka}f(s)\,ds+(t-ka)g_k
\\ &  
\overset{\eqref{defggkatomic}}{=} & 
\left(1-\frac{t-ka}{a}\right)\int_0^{ka}f(s)\,ds
+\frac{t-ka}{a}\int_0^{(k+1)a}f(s)\,ds
\\
& \le  &
\int_0^{\left(1-\frac{t-ka}{a}\right) ka+ \frac{t-ka}{a} 
(k+1)a
} f(s)\,ds 
=\int_0^t f(s)\,ds,
\end{eqnarray*}
where the inequality follows from the concavity of
$t\mapsto\int_0^t f(s)\,ds$ \cite[Theorem 24.2]{R70}. 
This implies that $g\prec\prec f$, and so $g\in M_\Phi(0,\infty)$.

Define
\[
\widetilde g(t):=
\begin{cases}
0, & 0<t\le a,\\
g(t-a), & t>a.
\end{cases}
\] 
Given $t>a$, choose an integer $k\ge0$ such that
$(k+1)a\le t<(k+2)a$, and so $ka\le
t-a<(k+1)a$.  
Since $f$ is decreasing, it follows that 
\[
\widetilde g(t)=g(t-a)\overset{\eqref{defggkatomic}}{=}
g_k
\overset{\eqref{defggkatomic}}{=}
\frac1a\int_{ka}^{(k+1)a}f(s)\,ds
\ge f((k+1)a)\ge f(t).
\]
Consequently, $0 \leq f\le f\chi_{(0,a)}+\widetilde g$ almost everywhere.
Since  $f\chi_{(0,a)}\in M_\Phi(0,\infty)\cap L_1(0,\infty)$ and $\widetilde g$ and $g$ are equimeasurable,  it follows from positivity and
symmetry of $\varphi$ in combination with $\varphi|_{M_\Phi\left(0,\infty\right)\cap L_1\left(0,\infty\right)}=0$ that 
\[
0<\varphi(f)
\le \varphi(f\chi_{(0,a)})+\varphi(\widetilde g)
=\varphi(g).
\]

Since $\mathbb N\times\mathbb N$ is countably infinite, we may
reindex the family $(r_j)_{j\in\mathbb N}$ as $(r_{n,k})_{n,k\ge0}$.  
For each $n\ge0$, set
\begin{align} \label{defvn}
v_n:=\sum_{k\ge0}g_kr_{n,k}.
\end{align}
Since $0\le g_k\le g_0<\infty$ and $r_{n,k}$ are mutually orthogonal, it follows that the sum in \eqref{defvn} converges in the strong operator topology and $0\le v_n\in \cA$.
Hence,
\begin{align} \label{mu(vn)=g}
\mu(v_n)\stackrel{\eqref{defggkatomic}}{=}g,
\qquad
\norm{v_n}_{M_\Phi(\cA,\tau|_{\cA})}
=\norm{g}_{M_\Phi(0,\infty)}.
\end{align}

By \cite[Theorem~4.2]{DPSS98}, the map
$x\mapsto\varphi(\mu(x))$ on $M_\Phi(\cA,\tau|_{\cA})_+$
extends to a bounded positive fully symmetric linear functional
$\psi\in M_\Phi(\cA,\tau|_{\cA})^*$. 
In particular,
\[
\psi(v_n)=\varphi(\mu(v_n))\stackrel{\eqref{mu(vn)=g}}{=}\varphi(g)>0,\qquad n\ge0.
\]
Since $(v_n)$ is a disjoint norm-bounded sequence of positive
elements and $\inf_{n\ge0}\psi(v_n)=\varphi(g)>0$, it follows from the construction
in \eqref{defQatomic} that there exists a bounded linear surjection from
$M_\Phi(\cA,\tau|_{\cA})$ onto $\ell_1$.
By Proposition~\ref{surjective image}(2), this contradicts the
Grothendieck property of $M_\Phi(\cA,\tau|_{\cA})$.
Therefore, $\mathop{\liminf}\limits_{t\to\infty}\frac{\Phi(2t)}{\Phi(t)}>1$.
\end{proof}

As Example~\ref{atomic}(1) shows, 
the assumptions of Proposition~\ref{prop:atomic-infty}
are insufficient to guarantee
\[
\liminf_{t\downarrow0}\frac{\Phi(2t)}{\Phi(t)}>1.
\]

Combining the preceding propositions with Theorem~\ref{thm:main},
we obtain the following two characterizations in the atomic case.

\begin{corollary}\label{cor:atomic-finite}
Let $(\mathcal M,\tau)$ be a nonzero semifinite atomic
von Neumann algebra with $\tau(\mathbf1)<\infty$.
Suppose that the projection hypothesis of
Proposition~\ref{prop:atomic-zero} holds.
Then, $M_\Phi(\mathcal M,\tau)$ has the Grothendieck property
if and only if
\[
\liminf_{t\downarrow0}\frac{\Phi(2t)}{\Phi(t)}>1.
\]
\end{corollary}

\begin{corollary}\label{cor:atomic-infinite}
Let $(\mathcal M,\tau)$ be a semifinite atomic
von Neumann algebra with $\tau(\mathbf1)=\infty$.
Suppose that the projection hypotheses of
Propositions~\ref{prop:atomic-zero} and~\ref{prop:atomic-infty}
both hold.
Then, $M_\Phi(\mathcal M,\tau)$ has the Grothendieck property
if and only if
\[
\liminf_{t\downarrow0}\frac{\Phi(2t)}{\Phi(t)}>1,
\qquad
\liminf_{t\to\infty}\frac{\Phi(2t)}{\Phi(t)}>1.
\]
\end{corollary}


\section{The Grothendieck property of intersections of Banach bimodules}
\label{s6}
The following theorem is the main result in this section.
\begin{theorem}\label{main}
Let $E,F\subseteq S(\M,\tau)$ be two Grothendieck Banach $\M$-bimodules.
Then,
\[
 X:=E\cap F,\qquad \|x\|_X:=\max\{\|x\|_E,\|x\|_F\},
\]
is Grothendieck. 
In particular, any finite intersection of
Grothendieck Banach $\M$-bimodules, equipped with the maximum norm,
is Grothendieck.
\end{theorem}

The following corollary is a direct consequence of Theorem \ref{main} when $\cM=L_\infty(0,\infty ) $, which answers questions in \cite{HL}.

\begin{corollary} \label{affirmative}
Let $E,F\subseteq S\left(0,\infty\right)$ be two Grothendieck Banach function spaces.
Then,
\[
 X:=E\cap F,\qquad \|x\|_X:=\max\{\|x\|_E,\|x\|_F\},
\]
is Grothendieck. 
\end{corollary}

\subsection{Some facts and symbols} 
Let $(\M,\tau)$ be a semifinite von Neumann algebra with a faithful
normal semifinite trace $\tau$.
Throughout this section, $E$ and $F$ are
two Grothendieck Banach $\cM$-bimodules.
We equip $X=E\cap F$ with
$\|x\|_X=\max\{\|x\|_E,\|x\|_F\}$. 
For more related commutative intersection results, the reader is referred to
\cite[Section~2]{HL}.

For a Banach $\cM$-bimodule $Y$, set
\[
 Y_{\mathrm{sa}}:=\{x\in Y:x=x^*\},\qquad
 Y_+:=\{x\in Y:x\ge0\},\qquad B_Y^+:=B_Y\cap Y_+.
\]
Here  $Y^*$ and $Y^{**}$ denote
the Banach dual and bidual. 
For any $\phi\in Y^*$, the conjugate functional $\overline \phi\in Y^*$ is defined by setting
$$\overline \phi(x)=\overline{\phi(x^*)}, \quad x\in Y.$$
A linear functional $\phi\in Y^*$ is called {\it Hermitian} if $\overline \phi=\phi.$
We denote the set of all Hermitian functionals $\phi\in Y^*$ by $Y_h^*$. 
Its positive cone is
\[
 Y^*_+:=\{h\in Y^*_h:h(x)\ge0\text{ for any }x\in Y_+\}.
\]
For $f,g\in Y^*_{h}$, the notation $f\le g$ means that
$g-f\in Y^*_+$. 
Define the square-root space of $Y$ by
\begin{align} \label{the def of D(2)}
 Y^{(2)}:=\{z\in S(\M,\tau):z^*z\in Y\},\qquad
 \|z\|_{Y^{(2)}}:=\|z^*z\|_Y^{1/2}.
\end{align}
Note that it is a linear $\M$-bimodule closed under adjoints
\cite[Lemma~4.2.4]{DPS}. 
Moreover, if $x,y\in Y^{(2)}$, then $xy\in Y$ \cite[Lemma 4.2.4]{DPS} and
\begin{equation}\label{mixed}
 \|xy\|_Y\le\|x\|_{Y^{(2)}}\|y\|_{Y^{(2)}},
 \quad x,y\in Y^{(2)}.
\end{equation}
Indeed, let $xy=v|xy|$ be the polar decomposition of $xy$.
For any $t>0$, by 
$$0\le \left(t^{1/2}x^*v-t^{-1/2}y\right)^*\left(t^{1/2}x^*v-t^{-1/2}y\right)
=tv^*xx^*v+t^{-1}y^*y-2|xy|,$$
we have that
\[
  2|xy|\le t\,v^*xx^*v+t^{-1}y^*y.
\]
Then,
\[
  2\|xy\|_Y
  \le t\|xx^*\|_Y+t^{-1}\|y^*y\|_Y
  = t\|x\|_{Y^{(2)}}^2+t^{-1}\|y\|_{Y^{(2)}}^2.
\]
Taking $t=\frac{\norm{y}_{Y^{(2)}}}{\norm{x}_{Y^{(2)}}}$ when $\norm{x}_{Y^{(2)}}\ne 0$,
we obtain that 
$\|xy\|_Y\le\|x\|_{Y^{(2)}}\|y\|_{Y^{(2)}}$.

\subsection{GNS representation} \label{small pre}
For a Hilbert space $H$ and a self-adjoint set
$\mathcal A\subseteq B(H)$, its \emph{commutant} is
\[
 \mathcal A':=\{T\in B(H):Ta=aT\text{ for every }a\in\mathcal A\}.
\]

Let $H$ be a Hilbert space and $\cA\subset B(H)$ be a unital von Neumann algebra. 
A {\it $*$-representation} of $\mathcal A$ on $H$
is a complex-linear map $\pi:\mathcal A\to B(H)$ such that
\[
 \pi(xy)=\pi(x)\pi(y),\qquad \pi(x^*)=\pi(x)^*,
 \qquad x,y\in\mathcal A.
\]
We say that $\pi$ is {\it unital} if $\pi(\mathbf 1)=\mathbf 1_H$
\cite[Definition~II.6.1.1]{Blackadar}.

A vector $h\in H$ is called {\it cyclic} for $\pi$ if
\[
 \overline{\{\pi(x)h:x\in\mathcal A\}}^{\,\|\cdot\|_H}=H.
\]
A representation admitting a cyclic vector is called a {\it cyclic
representation} \cite[Sections~II.6.4.2--II.6.4.3]{Blackadar}.

Let $\omega$
be a positive linear functional on $\mathcal \cA$. 
Set
\[
 N_\omega:=\{a\in\cA:\omega(a^*a)=0\},
 \qquad [a]:=a+N_\omega.
\]
Note that $a\in N_\omega$ implies that $\omega(a)=0$ \cite[Proposition II.6.2.4]{Blackadar}.
Equip $\cA/N_\omega$ with the inner product
$\langle[a],[b]\rangle:=\omega(b^*a)$, and let $H_\omega$ be its
Hilbert space completion of $\cA/N_\omega$ with respect to $\langle[\cdot],[\cdot]\rangle$. 

For each $b\in\mathcal A$, define
\[
 L_b:\cA/N_\omega\to \cA/N_\omega,\qquad L_b[a]:=[ba].
\]
By the GNS construction \cite[Section~II.6.4.1]{Blackadar},
every $L_b$ has a unique bounded extension, denoted by $\pi_\omega(b)$, on $H_\omega$.
The representation $\pi_\omega:\cA\to B(H_\omega)$ is called the {\it GNS representation} of $\cA$ associated with $\omega$.
Moreover, $\pi_\omega$ is a unital $*$-representation.

Set $\gamma_\omega:=[\mathbf 1]$. 
By \cite[Section~II.6.4.2]{Blackadar}, we have that
$\gamma_\omega$ is a cyclic vector for $\pi_\omega$ and
\[
 \omega(b)=\langle\pi_\omega(b)\gamma_\omega,\gamma_\omega\rangle,
 \qquad b\in\mathcal A.
\]
The triple $(\pi_\omega,H_\omega,\gamma_\omega)$ is called the
{\it GNS triple} of $\omega$.


\subsection{Hilbert space representation} \label{6.3}
Let $Y\subset S(\cM,\tau)$ be a Banach $\cM$-bimodule and 
suppose that $h:Y\to \mathbb C$ is a positive linear functional.
Define a semi-inner product $\langle \cdot,\cdot \rangle_h$ on $Y^{(2)}$ \cite[p. 255]{DPS} by
\begin{align} \label{the def of inner product}
\langle x,y \rangle_h=h(y^*x), \quad x,y \in Y^{(2)}.
\end{align}
Moreover,
\begin{align} \label{Cauchy-Schwarz}
    |h(y^*x)|^2\le h(x^*x)h(y^*y), \quad x,y \in Y^{(2)}.
\end{align}
Set
\begin{align} \label{the def of Nhj}
 N_h:=\{x\in Y^{(2)}:h(x^*x)=0\},\qquad [x]_h:=x+N_h, \quad x\in Y^{(2)}.
\end{align}
Let $H_h$ be the Hilbert space completion of the
quotient $Y^{(2)}/N_h$ and let
\[
 j_h:Y^{(2)}\longrightarrow H_h,\qquad j_h(x):=[x]_h .
\]
In this case,
$\langle[x]_h,[y]_h\rangle_{H_h}:=h(y^*x)$ is an inner product
on $Y^{(2)}/N_h$. 
Note that $j_h(Y^{(2)})$ is dense in $H_h$ and $\|j_h(x)\|_{H_h}^2=h(x^*x)$.

For any $a\in\M$, let $L_a[x]_h:=[ax]_h$.
By $$\|[ax]_h\|_{H_h}^2=h(x^*a^*ax)\le\|a\|_{\M}^2h(x^*x)=\|a\|_{\M}^2\|[x]_h\|_{H_h}^2,$$  we have that
$L_a$ is
bounded on $Y^{(2)}/N_h$. 
Denote its extension
to $H_h$ by $\pi_h(a)$. 
Note that the map $\pi_h:\M\to B(H_h)$ is a unital $*$-representation.


We have that $(Y^{(2)})^2=Y$ \cite[Lemma 4.2.4]{DPS}.
For any $A\in(\pi_h(\M)')_+$, set
\begin{align} \label{the def of hA on Y+}
h_A(t):=\|A^{1/2}[t^{1/2}]_h\|_{H_h}^2, \quad  t\in Y_+.
\end{align}
We claim that $h_A$ is additive on $Y_+$.
Indeed, fix two operators $x,y\in Y_+$ and set $d:=(x+y)^{1/2}\in Y^{(2)}$.
Let $\left(\begin{smallmatrix}x^{1/2}&0\\y^{1/2}&0\end{smallmatrix}\right)=\left(\begin{smallmatrix}u&0\\v&0\end{smallmatrix}\right)\operatorname{diag}(d,0)$ be the polar decomposition of $\left(\begin{smallmatrix}x^{1/2}&0\\y^{1/2}&0\end{smallmatrix}\right)$.
Then,
\begin{align} \label{x=ud}
    x^{1/2}=ud,\qquad y^{1/2}=vd,\qquad u^*u+v^*v=s(d).
\end{align}
Set $\xi:=A^{1/2}\left[d\right]_h$. 
Since $A^{1/2}\in\pi_h(\mathcal M)'$
and $s(d)d=d$, we have 
\begin{align} \label{invariant}
\pi_h\left(s\left(d\right)\right)\xi=A^{1/2}\pi_h(s(d))[d]_h=A^{1/2}[s(d)d]_h=\xi. 
\end{align}
Since $\pi_h$ is a $*$-representation, we have that
\begin{eqnarray} \label{similarly}
h_A(x)+h_A(y) &=&\|A^{1/2}[x^{1/2}]_h\|_{H_h}^2+\|A^{1/2}[y^{1/2}]_h\|_{H_h}^2\notag\\
 &\stackrel{\eqref{x=ud}}{=}&\|A^{1/2}[ud]_h\|_{H_h}^2+\|A^{1/2}[vd]_h\|_{H_h}^2\notag\\
 &=&\|A^{1/2}\pi_h(u)[d]_h\|_{H_h}^2+\|A^{1/2}\pi_h(v)[d]_h\|_{H_h}^2\notag\\
 &\stackrel{A^{1/2}\in\pi_h(\mathcal M)'}{=}&\|\pi_h(u)\xi\|_{H_h}^2+\|\pi_h(v)\xi\|_{H_h}^2\notag\\
 &=&\langle\pi_h(u^*u+v^*v)\xi,\xi\rangle_{H_h}\notag\\
 &\stackrel{\eqref{x=ud}, \eqref{invariant}}{=}&\|\xi\|_{H_h}^2\stackrel{\eqref{the def of hA on Y+},~d=(x+y)^{1/2}}{=}h_A(x+y).
\end{eqnarray}
This proves the claim.
Then, $h_A$ extends to a positive linear functional on $Y$. 
By \cite[Lemma~4.2.18]{DPS}, we have that $h_A$ is bounded.

For any $z\in Y^{(2)}$, let $z=u|z|$ be the polar decomposition of $z$, and
set $\eta:=A^{1/2}[|z|]_h$. 
Arguing similarly as \eqref{similarly}, we have that 
\[
 \langle A[z]_h,[z]_h\rangle_{H_h}
 =\|\pi_h(u)\eta\|_{H_h}^2
 =\langle\pi_h(u^*u)\eta,\eta\rangle_{H_h}
 \stackrel{\pi_h(u^*u)\eta=\eta}{=}\|\eta\|_{H_h}^2
 =h_A(z^*z).
\]
By polarization identity, we have that
\begin{align*} 
 h_A(y^*x)=\langle A[x]_h,[y]_h\rangle_{H_h},
 \qquad x,y\in Y^{(2)}, \quad A\in \left(\pi_h(\cM)'\right)_+.
\end{align*}

Since $(Y^{(2)})^2=Y$ by \cite[Lemma~4.2.4]{DPS}, the above formula
determines $h_A$ uniquely and extends $A\mapsto h_A$ linearly
from $(\pi_h(\M)')_+$ to $\pi_h(\M)'$.
Now, we have that for any $A\in \pi_h(\M)'$,
\begin{equation}\label{GNS}
 h_A(y^*x)=\langle A[x]_h,[y]_h\rangle_{H_h},
 \qquad x,y\in Y^{(2)}.
\end{equation}
Note that the map $A\in \pi_h(\M)'\mapsto h_A\in Y^*$ is linear and positive, and $h_\mathbf 1=h$ by $\langle[x]_h,[y]_h\rangle_{H_h}=h(y^*x)$.

Conversely, let $k\in Y^*_+$ with $0\le k\le h$. 
Consider the form on $Y^{(2)}/N_h\subseteq H_h$
\[
 ([x]_h,[y]_h)\longmapsto k(y^*x).
\]
By \eqref{Cauchy-Schwarz}, $|k(y^*x)|^2\le k(x^*x)k(y^*y)
\stackrel{0\le k\le h=h_\mathbf 1}{\le}\|[x]_h\|_{H_h}^2\|[y]_h\|_{H_h}^2$. 
Then, the form is well-defined and extends to a bounded sesquilinear form on $H_h$ (see the definition of \cite[p. 19]{DPS}). 
By the Riesz representation
theorem \cite[I.2.2.2]{Blackadar}, there exists a unique operator $C\in B(H_h)_+$ satisfying
\begin{align} \label{construction of C}
 k(y^*x)=\langle C[x]_h,[y]_h\rangle_{H_h}, \qquad x,y\in Y^{(2)}.
\end{align}
Taking $y=x$, we obtain that 
$$k(x^*x)=\langle C[x]_h,[x]_h\rangle_{H_h}\le h(x^*x)=h_\mathbf 1(x^*x)=\langle [x]_h,[x]_h\rangle_{H_h},$$
which together with the density of $Y^{(2)}/N_h$ in $H_h$ yields that
$0\le C\le \mathbf 1$. 
Given any $a\in\M$, by \eqref{construction of C}, we have that
\begin{align*}
    \langle C\pi_h(a)[x]_h,[y]_h\rangle_{H_h}
    & =\langle C[ax]_h,[y]_h\rangle_{H_h}=k(y^*ax)=k((a^*y)^*x)\\&=\langle C[x]_h,[a^*y]_h\rangle_{H_h}=\langle \pi_h(a) C[x]_h,[y]_h\rangle_{H_h},
\end{align*}
which together with the density of $Y^{(2)}/N_h$ in $H_h$ yields that $C\pi_h(a)=\pi_h(a)C$.
Thus, $C\in\pi_h(\M)'$ and $k=h_C$ by \eqref{GNS}. 

\subsection{Proof of Theorem \ref{main}}
The following lemma plays an important role in proving Theorem~\ref{main}.

\begin{lemma}\label{hilbert}
Let $\mathcal A\subseteq B(H_0)$ be a von Neumann algebra,
and let $\xi,\eta\in H_0$. Set
\begin{align} \label{teh def of m:}
 m:=\sup_{x\in B_{\mathcal A'}}|\langle x\xi,\eta\rangle_{H_0}|.
\end{align}
Then, there is an orthogonal projection $P\in\mathcal A$ such that
\[
 \|P\xi\|_{H_0}^2+\|(\mathbf 1-P)\eta\|_{H_0}^2\le m.
\]
If $\mathcal A=\pi(\M)'$ for a unital $*$-representation
$\pi:\M\to B(H_0)$, then
$m=\sup_{x\in B_\M}|\langle\pi(x)\xi,\eta\rangle_{H_0}|$.
\end{lemma}
\begin{proof}
Define the normal functional 
$\phi(x):=\langle x\xi,\eta\rangle_{H_0}$ for all $x\in \mathcal A'$.
By \cite[Proposition~III.2.3.1 and Proposition II.6.2.5(i)]{Blackadar}, there exist a unique positive linear functional, denoted by $|\phi|$, and a unique partial isometry $v\in \cA'$ such that for all $x\in \cA'$,
\begin{align} \label{|phi|}
\phi(x)=|\phi|(xv), \quad |\phi|(x)=\phi(xv^*), \quad |\phi|(\mathbf 1)=\norm{|\phi|}_{\cA'^*}=\norm{\phi}_{\cA'^*}=m, \quad v^*v= s(|\phi|).
\end{align}
Let $(\pi_\phi,H_\phi,\gamma)$ be a GNS triple obtained from Subsection \ref{small pre}, applied to $ \cA'$ and $|\phi|$, satisfying
\begin{align} \label{the def of |phi|}
 |\phi|(x)=\langle \pi_\phi(x)\gamma,\gamma\rangle_{H_\phi},
 \qquad x\in\mathcal A'.
\end{align}
Taking $x=\mathbf 1$, from $\pi_\phi(\mathbf 1)=\mathbf 1$ it follows that 
$\|\gamma\|_{H_\phi}^2=|\phi|(\mathbf 1)=m$.

Considering the Hilbert space $H_0\oplus H_\phi$, let 
\begin{align} \label{the def of sigma}
\sigma(x):=x\oplus\pi_\phi(x), \quad x\in \cA', \qquad
 w_1:=\left(\xi,\pi_\phi(v)\gamma\right),\qquad
 w_2:=\left(\eta,-\gamma\right).
\end{align}
Note that $\sigma$ is a $*$-representation.
For any $x\in\mathcal A'$, we have that
$$\sigma(x)w_1=\left(x\xi,\pi_\phi(x)\left(\pi_\phi(v)\gamma\right)\right)=\left(x\xi,\pi_\phi(xv)\gamma\right). 
$$
Then, for any $x\in \cA'$,
$$\langle\sigma(x)w_1,w_2\rangle_{H_0\oplus H_\phi}=\langle x\xi,\eta\rangle_{H_0}-\langle \pi_\phi(xv)\gamma,\gamma\rangle_{H_\phi}\stackrel{\eqref{the def of |phi|}}{=}\phi(x)-|\phi|(xv)\stackrel{\eqref{|phi|}}{=}0.$$
For any $y\in \cA'$, replacing $x$ with $y^*x$, we obtain that
\[
 \langle \sigma(x)w_1,\sigma(y)w_2\rangle_{H_0\oplus H_\phi}
 = \langle \sigma(y)^*\sigma(x)w_1,w_2\rangle_{H_0\oplus H_\phi}
 = \langle \sigma(y^*x)w_1,w_2\rangle_{H_0\oplus H_\phi}
 = 0,
\]
which implies that $\overline{\{\sigma(x)w_1:x\in\mathcal A'\}}\perp \overline{\{\sigma(x)w_2:x\in\mathcal A'\}}$.

Let $P_0$ be the orthogonal projection onto
$\overline{\sigma(\mathcal A')w_2}$.
Then, by $\sigma(\mathbf1)=\mathbf 1$, we have that 
\begin{align} \label{P0w1=0}
P_0w_1=P_0\mathbf 1w_1=0, \quad (\mathbf1-P_0)w_2=(\mathbf1-P_0)\mathbf 1w_2=0.
\end{align}

We claim that $P_0\in \sigma(\cA')'$.
Indeed, fix $y\in \cA'$.
Note that for any $x\in \cA'$,
$\sigma(y)\sigma(x)w_2=\sigma(yx)w_2\in \overline{\sigma(\mathcal A')w_2}$,
which implies that $\sigma(y)\overline{\sigma(\mathcal A')w_2}\subseteq \overline{\sigma(\mathcal A')w_2}$.
Then, 
\begin{align} \label{sigma(y)P}
P_0\sigma(y)P_0=\sigma(y)P_0.
\end{align}
Replacing $y$ with $y^*$, we obtain that
$P_0\sigma(y^*)P_0=\sigma(y^*)P_0$.
Taking adjoints in both sides, from $\sigma(y^*)^*=\sigma(y)$ it follows that
$P_0\sigma(y)P_0=P_0\sigma(y)$.
This together with \eqref{sigma(y)P} proves the claim.

By $v^*v=s(|\phi|)$, we have that $\sigma(v^*v)=\sigma(v)^*\sigma(v)$ is also a projection, 
and so $\norm{\sigma(v)}^2_{B\left(H_0\oplus H_\phi\right)}=\norm{\sigma(v)^*\sigma(v)}_{B\left(H_0\oplus H_\phi\right)}=\norm{\sigma(v^*v)}_{B\left(H_0\oplus H_\phi\right)}\le 1$.
Note that 
\begin{align*} 
P_0w_1 & =P_0(\xi,0)+P_0\left(0,\pi_\phi(v)\gamma\right)\stackrel{\eqref{P0w1=0}}{=}0, \\  \left(\mathbf1-P_0\right)w_2 &
=\left(\mathbf1-P_0\right)(\eta,0)+\left(\mathbf1-P_0\right)(0,-\gamma)\stackrel{\eqref{P0w1=0}}{=}0.
\end{align*}
Then, $\|(\mathbf 1-P_0)(\eta,0)\|_{H_0\oplus H_\phi}=\|(\mathbf 1-P_0)(0,\gamma)\|_{H_0\oplus H_\phi}$ by $(0,-\gamma)=-(0,\gamma)$ and
\begin{eqnarray} \label{le 0, gamma}
 \|P_0(\xi,0)\|_{H_0\oplus H_\phi}
 &=&\|P_0\left(0,\pi_\phi(v)\gamma\right)\|_{H_0\oplus H_\phi}\stackrel{\eqref{the def of sigma}}{=}\norm{P_0\sigma(v)(0,\gamma)}_{H_0\oplus H_\phi}
 \\ & \stackrel{P_0\sigma(v)=\sigma(v)P_0}{=} & \norm{\sigma(v)P_0(0,\gamma)}_{H_0\oplus H_\phi}
 \stackrel{\|\sigma(v)\|_{B(H_0\oplus H_\phi)}\le1}{\le}\|P_0(0,\gamma)\|_{H_0\oplus H_\phi}\notag.
\end{eqnarray}
Moreover, $\norm{P_0(0,\gamma)}^2_{H_0\oplus H_\phi}+\norm{(\mathbf 1-P_0)(0,\gamma)}^2_{H_0\oplus H_\phi}=\norm{\gamma}_{H_\phi}^2=m$.
Let $J:H_0\to H_0\oplus H_\phi$ given by $J\zeta=(\zeta,0), ~\zeta\in H_0$.
Note that $J^*J=\mathbf 1$ and $J^*(u,w)=u$ for all $(u,w)\in H_0\oplus H_\phi$.
Then, $J^*P_0J$ is a positive contraction.
For any $\zeta\in H_0$, we have that $\sigma(x)J\zeta\stackrel{\eqref{the def of sigma}}{=}(x\zeta,0)=Jx \zeta$, which implies that $\sigma(x)J=Jx$.
Furthermore, for all $(u,w)\in H_0\oplus H_\phi$,
\begin{align} \label{commute}
J^*\sigma(x)(u,w)\stackrel{\eqref{the def of sigma}}{=}J^*\left(xu,\pi_\phi(x)w\right)=xu=xJ^*(u,w).
\end{align}
This implies that for all $x\in \cA'$,
$$\left(J^*P_0J\right)x\stackrel{\sigma(x)J=Jx}{=}J^*P_0\sigma(x)J\stackrel{P_0\sigma(x)=\sigma(x)P_0}{=}J^*\sigma(x)P_0J\stackrel{\eqref{commute}}{=}xJ^*P_0J,$$
which shows that $J^*P_0J\in \cA''=\cA$.
Then, $J^*P_0J\le \norm{J^*P_0J}_{\infty}\mathbf 1\le \mathbf 1$.
Moreover,
\begin{eqnarray} \label{le m}
 &&  \langle J^*P_0J\xi,\xi\rangle_{H_0}+
 \langle(\mathbf 1-J^*P_0J)\eta,\eta\rangle_{H_0}\notag\\ &
  \stackrel{J^*J=\mathbf 1}
 {=}&\langle P_0(\xi,0),(\xi,0)\rangle_{H_0\oplus H_\phi}+
 \langle J^*(\mathbf 1-P_0)J\eta,\eta\rangle_{H_0}\notag\\
 &=& \norm{P_0(\xi,0)}^2_{H_0\oplus H_\phi}+\norm{\left(\mathbf 1-P_0\right)(\eta,0)}^2_{H_0\oplus H_\phi}\notag\\
 &\stackrel{\eqref{le 0, gamma}}{\le}&\|P_0(0,\gamma)\|^2_{H_0\oplus H_\phi}+\|(\mathbf 1-P_0)(0,\gamma)\|^2_{H_0\oplus H_\phi}= m.
\end{eqnarray}

For any $0<t<1$,
let $P_t:=e^{J^*P_0J}(t,1]\in\mathcal A$.
Applying \cite[Theorem 8.16]{R87} to $\varphi(t)=t$, $X=[0,1]$, $\mu_\nu(B)=\langle e^{J^*P_0J}(B)\nu,\nu \rangle$ and $f(\lambda)=\lambda$, we have that
\[
 \int_0^1\bigl(\|P_t\xi\|^2_{H_0}+\|(\mathbf 1-P_t)\eta\|_{H_0}^2\bigr)\,dt
 \stackrel{\substack{J^*P_0J\le \mathbf 1\\ \tiny\mbox{\cite[Thm.1.5.7 (ii)]{DPS}}}}{=}\langle J^*P_0J\xi,\xi\rangle_{H_0}+
   \langle\left(\mathbf 1-J^*P_0J\right)\eta,\eta\rangle_{H_0}\stackrel{\eqref{le m}}{\le} m.
\]
Then, there exists $t\in (0,1)$ such that $\|P_t\xi\|^2_{H_0}+\|(\mathbf 1-P_t)\eta\|^2_{H_0}\le m$ and in this case, 
take $P=P_t$.

For the second statement, fix a unital $*$-representation $\pi:\cM\to B(H_0)$.
Take $\cA=\pi(\cM)'$ and in this case, $\cA'=\pi(\M)''$.
Recall that the  map
\[
\mathcal M/\ker\pi \longrightarrow \pi(\mathcal M),
\qquad x+\ker\pi\longmapsto\pi(x),
\]
is an isometric $*$-isomorphism
\cite[Corollary~II.2.2.9]{Blackadar}.
By 
 \cite[Proposition~II.5.1.5]{Blackadar},
we have that for any $b\in \pi(\cM)$,
there exists $x\in \cM$ such that
\[
\pi(x)=b,\qquad \|x\|_\cM=\|b\|_{B(H_0)}.
\]
Since $\pi$ is contractive \cite[Corollary II.1.6.6]{Blackadar}, it follows that
$\pi(B_{\mathcal M})=B_{\pi(\mathcal M)}$.

Now, we may fix $b\in B_{\mathcal A'}=B_{\pi(\cM)''}$.
Since $\pi$ is unital, it follows from \cite[Theorem~I.9.1.3]{Blackadar} and $\pi(B_{\mathcal M})=B_{\pi(\mathcal M)}$ that
there exists a net
$(x_\lambda)\subset B_{\mathcal M}$ such that
$\pi(x_\lambda)\to b$ in the uso-topology.
Consequently,
\[
\langle\pi(x_\lambda)\xi,\eta\rangle_{H_0}
\longrightarrow \langle b\xi,\eta\rangle_{H_0},
\]
and so
\[
\sup_{b\in B_{\cA'}}|\langle b\xi,\eta\rangle_{H_0}|
\le
\sup_{x\in B_{\mathcal M}}
|\langle\pi(x)\xi,\eta\rangle_{H_0}|.
\]
The reverse inequality follows from
$\pi(B_{\mathcal M})\stackrel{\pi(B_{\mathcal M})=B_{\pi(\mathcal M)}}{\subseteq} B_{\mathcal A'}$.
Therefore,
\[
m\stackrel{\eqref{teh def of m:}}{=}\sup_{b\in B_{\mathcal A'}}|\langle b\xi,\eta\rangle_{H_0}|
 =\sup_{x\in B_{\mathcal M}}
   |\langle\pi(x)\xi,\eta\rangle_{H_0}|.
\]
The proof is complete.
\end{proof}


Let $Y$ be a Banach $\cM$-bimodule.
We say that $Y$ is \emph{positively Grothendieck} if every positive
weak$^*$-null sequence in $Y^*$ is weakly null. 
Write $(G_j)\prec(g_n)$ if
\[
 G_j=\sum_{n=p_j}^{q_j}\alpha_{j,n}g_n,\qquad
 \alpha_{j,n}\ge0,\quad \sum_{n=p_j}^{q_j}\alpha_{j,n}=1,
 \qquad 1\le p_j\le q_j<p_{j+1}.
\]
We call $(G_j)$ a \emph{convex block sequence} of $(g_n)$ \cite[p. 1]{AR16}.
Note that $\prec$ is transitive.
The approaches used  in  \cite{M91} for proving the equivalence between  the Grothendieck property and positive Grothendieck property are hard to adapt
 in the noncommutative setting. 
The above notion of convex block sequence is the key tool in proving such an equivalence. 




The next lemma aims to select a convex block sequence so that, for any positive functionals dominated by further convex blocks, weak$^*$ convergence to zero can be verified by  a single element of $Y_+$.

\begin{lemma}\label{control}
Let $Y\subseteq S(\M,\tau)$ be a positively Grothendieck Banach
$\M$-bimodule, and let $(g_n)_{n\ge1}$ be a bounded sequence in $(Y^*)_+$.
Then, there exist a convex block sequence $(G_j)\prec(g_n)$ and
an element $u\in Y_+$ such that for every $(L_j)\prec(G_j)$ and
every sequence $(w_j)\subseteq Y^*$,
\begin{equation}\label{control-eq}
 0\le w_j\le L_j,\qquad w_j(u)\to0
 \quad\Longrightarrow\quad
 w_j(x)\to0\quad\text{for every }x\in Y .
\end{equation}
\end{lemma}
\begin{proof}
{\bf Step $1$}. We aim to construct a convex block sequence 
of $(g_n)$.
Set $C:=\sup_n\|g_n\|_{Y^*}$
and $R_{0,n}:=g_n$ for $n\ge1$.
Let $j=1$.
Set
\[
 \delta_1:=\sup\bigl\{\limsup_{n\to\infty}w_n^{(1)}(x):
 x\in B_Y^+,\ (L_n^{(1)})\prec(g_n),\ 0\le w_n^{(1)}\le L_n^{(1)}\bigr\}.
\]
By $(L_n^{(1)})\prec(g_n)$ and $0\le w_n^{(1)}\le L_n^{(1)}$, we have that $\|w_n^{(1)}\|_{Y^*}\le\|L_n^{(1)}\|_{Y^*}\le C$, and so $0\le \delta_1\le \sup\limsup\norm{w_n^{(1)}}_{Y^*}\le C$.
If $\delta_1>0$, then we may choose $x_1\in B_Y^+$, $(L_n^{(1)})$ and $(w_n^{(1)})$ as above with
$\limsup_n w_n^{(1)}(x_1)>\delta_1/2$. 
By the definition of $\limsup$, there exists a strictly increasing sequence $(n_{1,k})_{k\ge1}$ 
such that
$w_{n_{1,k}}^{(1)}(x_1)>\delta_1/2$ for any $k\ge1$.
Set
\[
 R_{1,k}:=L_{n_{1,k}}^{(1)},\quad b_{1,k}:=w_{n_{1,k}}^{(1)},\quad k\ge1.
\]
If $\delta_1=0$, then we may take $x_1=0$, $R_{1,k}=g_k$
and $b_{1,k}=0$ for all $k\ge1$. 
In either case, take $G_1:=R_{1,1}$.



Now, consider $j\ge2$. 
Suppose that the sequence
$(R_{j-1,k})_{k\ge1}$ and the elements $x_1,\ldots,x_{j-1}\in B_Y^+$
have been chosen. 
Let
\[
 \begin{aligned}
 \delta_j:=\sup\bigl\{\limsup_{n\to\infty}w_n^{(j)}(x):\
&x\in B_Y^+,\ (L_n^{(j)})\prec(R_{j-1,k})_{k\ge1},\\
 &0\le w_n^{(j)}\le L_n^{(j)}, ~~ w_n^{(j)}(x_r)\to0, ~\forall1\le r<j\bigr\}.
 \end{aligned}
\]
Since $0\le w_n^{(j)}\le L_n^{(j)}$ and $(L_n^{(j)})\prec (R_{j-1,k})_{k\ge1}$, we have that $\|w_n^{(j)}\|_{Y^*}\le\|L_n^{(j)}\|_{Y^*}\le C$, and so $0\le \delta_j\le C$.
If $\delta_j>0$, then we may choose $x_j$, $(L^{(j)}_n)$ and $(w^{(j)}_n)$ satisfying the above conditions and $\limsup_n w^{(j)}_n(x_j)>\delta_j/2$.
There exists a strictly increasing sequence $(n_{j,k})_{k\ge1}$ such that $w^{(j)}_{n_{j,k}}(x_j)>\delta_j/2$ for any $k\ge1$.
Set $R_{j,k}:=L^{(j)}_{n_{j,k}}$ and $b_{j,k}:=w^{(j)}_{n_{j,k}}$.
If $\delta_j=0$, then we may take $x_j=0$, $R_{j,k}=R_{j-1,k}$ and $b_{j,k}=0$.
In this case, we have that for any $j\ge1$,
\begin{align} \label{step 1}
 \begin{aligned}
 &(R_{j,k})_{k\ge1}\prec(R_{j-1,k})_{k\ge1},\qquad
       0\le b_{j,k}\le R_{j,k},\\
 &b_{j,k}(x_j)\ge\frac{\delta_j} 2, \quad k\ge1,\qquad
       b_{j,k}(x_r)\to0,\quad k\to\infty,\ 1\le r<j.
 \end{aligned}
\end{align}

For any $j\ge1$, let $G_j:=R_{j,j}$.
By the definition of $G_j$ and $(R_{j,k})_{k\ge1}\prec(R_{j-1,k})_{k\ge1}$ for any $j\ge1$, we have that
\begin{align} \label{Gmprec}
 (G_j)\prec(g_n),\qquad
 (G_m)_{m\ge j}\prec(R_{j,k})_{k\ge1},\quad j\ge1.
\end{align}

{\bf Step $2$}. Below, we prove that $\delta_j\to0$ as $j\to \infty$.
By $(G_j)\prec (g_n)$, we have that $\|G_m\|_{Y^*}\le C$ for all $m\ge j$.
This together with \cite[F.5]{AK06} and \cite[Theorem 4.29]{F99} yields that there exists a subnet $(G_{m_\alpha})_\alpha$ of $(G_m)$ converging weak$^*$ to some $g\in Y^*_+$.

Fix $j\ge1$. 
For any $m\ge j$, by $(G_m)_{m\ge j}\stackrel{\eqref{Gmprec}}{\prec}(R_{j,k})_{k\ge1}$, we may set
\[
 G_m=\sum_{k\in I_m}\alpha_{m,k}R_{j,k}\in Y^*_+,\qquad
 \alpha_{m,k}\ge0,\qquad \sum_{k\in I_m}\alpha_{m,k}=1, \quad \min I_m\ge m.
\]
Moreover, $I_m$ are mutually disjoint.
Set
\[
 B_{j,m}:=\sum_{k\in I_m}\alpha_{m,k} b_{j,k}\in Y^*_+.
\]
By \eqref{step 1}, we have that
\begin{align} \label{limit before}
 \begin{aligned}
 &0\le B_{j,m}\le G_m,\qquad B_{j,m}(x_j)\ge\frac{\delta_j}2,\\
0\le B_{j,m}(x_r)  & \le\sup_{k\ge m} b_{j,k}(x_r)\to0,\quad m\to\infty, ~1\le r<j.
 \end{aligned}
\end{align}
Since $\norm{B_{j,m}}\le \sup_{k\in I_m}\norm{b_{j,k}}=\sup_{k\in I_m}\norm{w_{n_{j,k}}^{(j)}}\le C$, it follows from \cite[F.5]{AK06} and \cite[Theorem 4.29]{F99} that there exists a subnet $(B_{j,m_{\alpha_\lambda}})_\lambda$ of $(B_{j,m_\alpha})_{\alpha}$ converging weak$^*$ to some $\beta_j\in Y^*_+$.
Then,
 $(G_{m_{\alpha_\lambda}})_\lambda$
converges weak$^*$ to $g$.
Moreover, by \eqref{limit before}, we have that
\begin{align} \label{betaj(xr)=0}
 0\le\beta_j\le g,\qquad \beta_j(x_j)\ge\frac{\delta_j}2,\qquad
 \beta_j(x_r)=0,\quad 1\le r<j.
\end{align}

Let $N_g$ be defined as in \eqref{the def of Nhj}.
Let $j_{g}:Y^{(2)}\to H_{g}$ given by $j_{g}(x)=[x]_{g}=x+N_{g}$, where $H_{g}$ denotes the Hilbert space completion of $Y^{(2)}/N_{g}$.
For any $a\in \cM$, set $L_a[x]_{g}:=[ax]_{g}$ and we denote its extension to $H_{g}$ by $\pi_{g}(a)$.
Moreover, the map $\pi_{g}:\cM \to B(H_{g})$ is a unital $*$-representation.

Applying \eqref{construction of C} to $\beta_j$, we obtain that there exists a unique operator $A_j\in \pi_g(\cM)'$ with $0\le A_j\le \mathbf 1$ such that
\begin{align} \label{the def of betaj}
\beta_j(y^*x)=\langle A_j[x]_{g}, [y]_g \rangle_{H_g}, \qquad x,y \in Y^{(2)}.
\end{align}
Then, for any $1\le r< j$,
\begin{align} \label{==0}
\Biggl\|A_j^{1/2}j_g(x_r^{1/2})\Biggr\|_{H_g}^2=\beta_j(x_r)\stackrel{\eqref{betaj(xr)=0}}{=}0.
\end{align}
Let $P_0:=0$ and let $P_j$ be the orthogonal projection onto
\[
 K_j:=\overline{\operatorname{span}}\left\{\pi_g(a)j_g(x_r^{1/2}):
                       a\in\M,\ 1\le r\le j\right\}.
\]
Since $K_j$ is invariant under $\pi_g(\M)$ ($\pi_g$ is a $*$-representation), it follows from \cite[II.6.1.4]{Blackadar} that 
$P_j\in\pi_g(\M)'$.
Note that 
$P_{j-1}\le P_j$ and $P_jj_g(x_j^{1/2})=j_g(x_j^{1/2})$ by $\pi_g(\mathbf 1)=\mathbf 1$.
By $A_j^{1/2}j_g(x_r^{1/2})\stackrel{\eqref{==0}}{=}0$ for any  $1\le r<j$,
we have 
\[
 A_j^{1/2}\pi_g(a)j_g(x_r^{1/2})\stackrel{A_j\in \pi_g(\cM)'}{=}\pi_g(a)A_j^{1/2}j_g(x_r^{1/2})=0,
 \quad a\in\M,\ 1\le r\le j-1,
\]
which implies that $A_j^{1/2}P_{j-1}=0$. 
Since $P_j-P_{j-1}\in \pi_g(\cM)'$ is a projection, we have that
$h_j:=g_{P_j-P_{j-1}}\in(Y^*)_+$, defined in \eqref{GNS}, satisfies
 \begin{eqnarray} \label{deltajto0}
 \frac{\delta_j} 2
 &\stackrel{\eqref{betaj(xr)=0}, ~\eqref{the def of betaj}}{\le}& \Biggl\|A_j^{1/2}j_g(x_j^{1/2})\Biggr\|_{H_g}^2
   \stackrel{\substack{P_jj_g(x_j^{1/2})=j_g(x_j^{1/2})\\A_j^{1/2}P_{j-1}=0}}{=}\Biggl\|A_j^{1/2}\left(P_j-P_{j-1}\right)j_g(x_j^{1/2})\Biggr\|_{H_g}^2\notag\\
 &\le& \Biggl\|\left(P_j-P_{j-1}\right)j_g(x_j^{1/2})\Biggr\|_{H_g}^2
   \overset{\eqref{GNS}}{=} h_j(x_j)\stackrel{x_j\in B_Y^+}{\le}\|h_j\|_{Y^*}.
 \end{eqnarray}
Noting that the map $A\in \pi_g(\cM)'\mapsto g_A$ is linear and positive (see \eqref{GNS}), with $g_\mathbf 1=g$, 
we have that
\[
 \sum_{j=1}^n h_j\stackrel{h_j=g_{P_j-P_{j-1}}}{=}g_{P_n}\le g,\qquad n\ge1.
\]

Set $s_n:=\sum_{j=1}^n h_j$. 
By $0\le s_n\le g$, the sequence
$s_n(x)$ increases to a finite limit for every $x\in Y_+$. 
By $Y=Y_+-Y_++i(Y_+-Y_+)$, we have that $s_n\to s$ weak$^*$ for some
$s\in(Y^*)_+$ with $s_n\le s\le g$ ($s(x):=\lim_n s_n(x)$ for all $x\in Y$).
Since $s-s_n$ are weak$^*$ null, it follows from the positive Grothendieck property of $Y$ that $s-s_n$ are weakly null.
Given $\varepsilon>0$,
by \cite[Theorem~3.13]{R91}, there exists a convex combination $r_\varepsilon$
of $s-s_1,\ldots,s-s_N$ with $\|r_\varepsilon\|_{Y^*}<\varepsilon$.
Since $(s-s_n)$ decreases,  we have that for any $j>N$,
\[
  0\le h_j\le s-s_N\le r_\varepsilon.
\]
Then, $\|h_j\|_{Y^*}\le\|r_\varepsilon\|_{Y^*}<\varepsilon$
by \cite[Corollary~4.2.8]{DPS}. 
Hence, $\|h_j\|_{Y^*}\to0$.
This together with \eqref{deltajto0} yields that
$\delta_j\to0$ as $j\to \infty$.

{\bf Step $3$}. We construct a control element in $ Y_+$ and verify 
 \eqref{control-eq}.
Since $x_r\in B_Y^+$, we have that the series $u:=\sum_{r=1}^{\infty}2^{-r}x_r$ converges in
$Y$ and $u\in Y_+$.
Fix $(L_k)\prec(G_j)$ ($G_j=R_{j,j}$)
and $0\le w_k\le L_k$ satisfying $w_k(u)\to0$ as $k\to \infty$. 
For any $r\ge1$,
\begin{align} \label{wkxrto0}
 0\le w_k(x_r)\stackrel{u\ge 2^{-r}x_r}{\le} 2^r w_k(u)\longrightarrow0 .
\end{align}

Fix $j\ge1$. 
For any $k\ge j$, the block $L_k$ uses only $G_m$ with
$m\ge k$. 
By \eqref{Gmprec} and \eqref{step 1}, we have that
$$(L_k)_{k\ge j}\stackrel{(L_k)\prec(G_j)}{\prec} (G_m)_{m\ge j} \stackrel{\eqref{Gmprec}}{\prec}(R_{j,k})_{k\ge1}\stackrel{\eqref{step 1}}{\prec}(R_{j-1,k})_{k\ge1},$$
which implies that $(L_k)_{k\ge j}\prec (R_{j-1,k})_{k\ge1}$ by transitivity of $\prec$.
This together with $w_k(x_r)\stackrel{\eqref{wkxrto0}}{\to}0$ for any $1\le r<j$ and the definition of $\delta_j$   shows that 
for any $x\in B_Y^+$,
\[
 0\le\limsup_{k\to\infty} w_k(x)\le\delta_j, \qquad j\ge1.
\]
Letting $j\to\infty$, by $\delta_j\to 0$ as $j\to \infty$, we obtain that $w_k(x)\to0$ for all $x\in B_Y^+$. 
Therefore, $w_k(x)\to0$ for any $x\in Y$.
The proof is complete.
\end{proof}

We now use Lemma \ref{hilbert} and Lemma \ref{control} to show that positive functionals suffice to test the Grothendieck property of a Banach $\cM$-bimodule.

\begin{theorem}\label{criterion}
Let $Y\subseteq S(\M,\tau)$ be a Banach $\M$-bimodule.
Then, $Y$ is Grothendieck if and only if it is positively
Grothendieck.
\end{theorem}
\begin{proof}
It suffices to prove the sufficiency.
Assume by contradiction that $Y$ is positively
Grothendieck but not Grothendieck. 
Then, there exists a weak$^*$-null sequence $(\phi_n)_n\subseteq Y^*$ such that $(\phi_n)$ is not weakly null and $(\phi_n)_n\not\subseteq Y^*_+$.
We divide the proof into several steps.

{\bf Step $1$}. Below,  we   choose Hermitian parts of $(\phi_n)$ which is not weakly null.
Consider two Hermitian functionals $f_n:=\frac{\phi_n+\overline\phi_n}{2}$ and $g_n:=\frac{\phi_n-\overline\phi_n}{2i}$.
Note that $\phi_n=f_n+ig_n$.
Hence, at least
one of $(f_n)$ and $(g_n)$ is weak$^*$ null but not weakly null.
Without loss of generality, we may assume that $(f_n)$  is weak$^*$ null but not weakly null.
Note that $(f_n)$ is bounded in $Y^*$ \cite[D.12]{AK06}.

Since $(f_n)$ is not weakly null, it follows that there exist $y^{**}\in Y^{**}$,
$\varepsilon>0$ and a subsequence, still denoted by $(f_n)$, such
that $|y^{**}(f_n)|\ge\varepsilon$. 
By 
\cite[Theorem~1]{Ando62},
we have that
\[
y^{**}=(y^{**}_1-y^{**}_2)+i(y^{**}_3-y^{**}_4),
\quad y^{**}_j\in(Y^{**})_+, ~1\le j\le 4.
\]
Since $f_n$ are Hermitian, it follows from \cite[p. 252]{DPS} that $y^{**}_j(f_n)\in \mathbb R$.
Then, there exists $j$ such that $|y^{**}_j(f_n)|\ge\frac \varepsilon 4$ for a further
subsequence, denoted by $(f_n)$. 
Passing to a further subsequence, we may assume that $y^{**}_j(f_n)$ have the same sign.
If $y^{**}_j(f_n)<0$ for all $n$, then replace $f_n$ with $-f_n$.
In this case, $y^{**}_j(f_n)\ge\frac \varepsilon 4$ for all $n$.
Letting $x^{**}:=\frac{1}{\max\{1,\norm{y^{**}_j}_{Y^{**}}\}}y^{**}_j$,
we obtain 
$$x^{**}\in B_{Y^{**}}\cap(Y^{**})_+, \quad x^{**}(f_n)>\delta, \quad n\ge1,$$
where $\delta=\frac \varepsilon {8\max\{1,\norm{y^{**}_j}_{Y^{**}}\}}$.


{\bf Step $2$}. Below, we  choose a control element in $Y_+$ obtained from Lemma \ref{control} and  construct bounded linear functionals on $\cM$ with norms tending to zero. 
By \cite[proposition 4.2.2]{DPS}, there exists $(h_n)\subseteq Y^*_+$ such that
\begin{align} \label{hn bounded}
-h_n\le f_n\le h_n,\quad \norm{h_n}_{Y^*}\le 2\norm{f_n}_{Y^*}.
\end{align}
Applying Lemma~\ref{control} to $(h_n)$, there exist a convex block sequence $(\widetilde h_j)\prec (h_n)$ and $u\in Y_+$ satisfying \eqref{control-eq}.
We may assume that $B_j$ are finite, 
mutually disjoint and
$$\widetilde h_j:=\sum_{n\in B_j}\alpha_{j,n}h_n, \quad \alpha_{j,n}\ge0, \quad \sum_{n\in B_j}\alpha_{j,n}=1, \quad \min B_{j+1}> \max B_{j}, \quad j\ge1.$$
For all $j\ge1$, set $\widetilde f_j:=\sum_{n\in B_j}\alpha_{j,n}f_n$.
Note that $-\widetilde  h_n\le \widetilde f_n\le \widetilde h_n$ and $x^{**}(\widetilde f_n)>\delta$ for all $n\ge1$.


For any $n\ge 1$, define
\begin{align} \label{the def of lambdan}
 \lambda_n(a):=\widetilde f_n\left(u^{1/2}\,au^{1/2}\right),\qquad a\in\M.
\end{align}
Note that $u^{1/2}\,au^{1/2}\in Y$ and $j\le \min B_j\to \infty$ as $j\to \infty$.
Since $(f_n)$ is weak$^*$ null in $Y^*$,
it follows 
that 
\[
|\lambda_n(a)|
= \left|\widetilde f_n\left(u^{1/2}\,au^{1/2}\right)\right|
= \left|\sum_{k\in B_n}\alpha_{n,k}f_k\left(u^{1/2}\,au^{1/2}\right)\right|
\le \sup_{k\ge \min B_n}\left|f_k\left(u^{1/2}\,au^{1/2}\right)\right|
\to 0,
\]
which implies that
$(\lambda_n)$ is weak$^*$ null in $\cM^*$, and so $(\lambda_n)$ is bounded \cite[D.12]{AK06}.
By the fact that $\cM$ is Grothendieck \cite[Corollary 7]{P94}, we obtain that $(\lambda_n)$ is weakly null in $\cM^*$.

By 
\cite[Theorem~3.13, p.~67]{R91},
there exist a nonempty finite set of
indices $I_1\subseteq\{1,2,\ldots\}$ and coefficients $c_{1,n}>0$, $n\in I_1$, such that
\[
   \sum_{n\in I_1}c_{1,n}=1,\qquad
   \Lambda_1:=\sum_{n\in I_1}c_{1,n}\lambda_n,\qquad
   \|\Lambda_1\|_{\M^*}<\frac12.
\]
Set $q_1:=\max I_1$.
By induction, we may assume that
$I_j$ has been chosen and $q_j:=\max I_j$. 
Applying \cite[Theorem~3.13, p.~67]{R91} to tail
$(\lambda_n)_{n>q_j}$, we obtain a nonempty finite set $I_{j+1}\subseteq\{q_j+1,q_j+2,\ldots\}$
and coefficients $c_{j+1,n}>0$, $n\in I_{j+1}$, such that
\begin{align} \label{the def of Lambdaj}
   \sum_{n\in I_{j+1}}c_{j+1,n}=1,\qquad
   \Lambda_{j+1}:=\sum_{n\in I_{j+1}}c_{j+1,n}\lambda_n,\qquad
   \|\Lambda_{j+1}\|_{\M^*}<2^{-(j+1)}.
\end{align}
Set $q_{j+1}:=\max I_{j+1}$. 
Note that
\[
   \max I_j<\min I_{j+1},\qquad
   \|\Lambda_j\|_{\M^*}<2^{-j},\quad j\ge1.
\]

{\bf Step $3$}. 
Define
\begin{align} \label{the def of Kj}
   F_j:=\sum_{n\in I_j}c_{j,n}\widetilde f_n,\qquad
   K_j:=\sum_{n\in I_j}c_{j,n}\widetilde h_n.   
\end{align}
Observe that $(K_j)\prec (\widetilde h_n)$.
Since $c_{j,n}\ge0$ and $\sum_{n\in I_j} c_{j,n}=1$, it follows from $-\widetilde  h_n\le \widetilde f_n\le \widetilde h_n$ and $x^{**}(\widetilde f_n)>\delta$ that
\begin{align} \label{x**>delta}
   -K_j\le F_j\le K_j,\qquad
   x^{**}(F_j)=\sum_{n\in I_j} c_{j,n}x^{**}(\widetilde f_n)>\delta.
\end{align}
Moreover, for any $a\in\M$,
\begin{align*} 
   F_j(u^{1/2}au^{1/2})
   =\sum_{n\in I_j} c_{j,n}\widetilde f_n(u^{1/2}au^{1/2})
   \stackrel{\eqref{the def of lambdan}}{=}\sum_{n\in I_j} c_{j,n}\lambda_n(a)
   \stackrel{\eqref{the def of Lambdaj}}{=}\Lambda_j(a).
\end{align*}
Taking the supremum over all $a\in B_{\M}$, we obtain
\begin{align} \label{<2-j}
   \sup_{a\in B_{\M}}\left|F_j(u^{1/2}au^{1/2})\right|
   =\|\Lambda_j\|_{\M^*}\stackrel{\eqref{the def of Lambdaj}}{<}2^{-j}.
\end{align}

Fix $j\ge1$. 
Since $-K_j\le F_j\le K_j$, the functional
$k_j:=(K_j+F_j)/2$ satisfies $0\le k_j\le K_j$. 
Let $Y^{(2)}$ and $N_{K_j}$ be defined as in \eqref{the def of D(2)} and \eqref{the def of Nhj}, respectively.
By the construction in Subsection \ref{6.3},
Let $j_{K_j}:Y^{(2)}\to H_{K_j}$ given by $j_{K_j}(x)=[x]_{K_j}=x+N_{K_j}$ , where $H_{K_j}$ denotes the Hilbert space completion of $Y^{(2)}/N_{K_j}$.
For any $a\in \cM$, set $L_a[x]_{K_j}:=[ax]_{K_j}$ and we denote its extension to $H_{K_j}$ by $\pi_{K_j}(a)$.
Moreover, the map $\pi_{K_j}:\cM \to B(H_{K_j})$ is a unital $*$-representation.

By \eqref{construction of C},
there exists a
positive contraction $C_j\in \pi_{K_j}(\M)'$ with $0\le C_j\le \mathbf 1$ such that
\[
 k_j(y^*x)=\langle C_j[x]_{K_j},[y]_{K_j}\rangle_{H_{K_j}},
 \qquad x,y\in Y^{(2)} .
\]
Set $A_j:=2C_j-\mathbf1$. 
Then, $A_j=A_j^*$ and $\|A_j\|_{B(H_{K_j})}\le1$.
Since $F_j=2k_j-K_j$ and
$K_j(y^*x)\stackrel{\eqref{GNS}}{=}\langle[x]_{K_j},[y]_{K_j}\rangle_{H_{K_j}}$ by $(K_j)_\mathbf 1=K_j$, we obtain
\begin{align} \label{125}
 F_j(y^*x)=\langle A_j[x]_{K_j},[y]_{K_j}\rangle_{H_{K_j}},
 \qquad x,y\in Y^{(2)} .
\end{align}

Note that $u^{1/2}\in Y^{(2)}$.
Let $\alpha_j:=j_{K_j}(u^{1/2})=[u^{1/2}]_{K_j}\in H_{K_j}$.
For any $a\in\M$, we have that
$\pi_{K_j}(a)\alpha_j=L_a[u^{1/2}]_{K_j}=[au^{1/2}]_{K_j}$.
Taking $x=au^{1/2}$ and $y=u^{1/2}$ in \eqref{125} and
using $A_j=A_j^*$, we have
\begin{align} \label{fj=}
 F_j(u^{1/2}au^{1/2})
 =\langle A_j[au^{1/2}]_{K_j},[u^{1/2}]_{K_j}\rangle_{H_{K_j}}
 \stackrel{A_j=A_j^*}{=}\langle\pi_{K_j}(a)\alpha_j,A_j\alpha_j\rangle_{H_{K_j}} .
\end{align}
Applying Lemma~\ref{hilbert} to $\mathcal A=\pi_{K_j}(\cM)'$,
$\xi=\alpha_j$ and $\eta=A_j\alpha_j$, we obtain that there exists an orthogonal
projection $P_j\in \pi_{K_j}(\cM)'$ such that
\begin{eqnarray} \label{normto0}
 \|P_j\alpha_j\|^2_{H_{K_j}}+\|(\mathbf 1-P_j)A_j\alpha_j\|^2_{H_{K_j}}
 &\le&\sup_{a\in B_\M}
      \left|\langle\pi_{K_j}(a)\alpha_j,A_j\alpha_j\rangle_{H_{K_j}}\right|\notag\\
 &\stackrel{\eqref{fj=}}{=}&\sup_{a\in B_\M}\left|F_j(u^{1/2}au^{1/2})\right|
   \stackrel{\eqref{<2-j}}{<}
   2^{-j}.
\end{eqnarray}
Note that $(K_j)$ is a convex block sequence of $(\widetilde h_n)$ by \eqref{the def of Kj}.
By the fact that $(K_j)\prec (\widetilde h_n)$ and $\left((h_n),(\widetilde h_n),u\right)$ satisfy \eqref{control-eq}, we have that 
for any $(\psi_j)\subseteq Y^*$ satisfying $0\le \psi_j\le K_j$ and $\psi_j(u)\to0$,
$(\psi_j)$ is weak$^*$ null in $Y^*$.
Set
\begin{align} \label{the def of Tj}
 T_j:=\frac12\bigl(P_j+A_j(\mathbf 1-P_j)A_j\bigr)\in \pi_{K_j}(\cM)',
 \qquad b_j:=(K_j)_{T_j}
\end{align}
(see the definition of $b_j$ in \eqref{GNS}).

{\bf Step $4$}. Below, we  prove that $0\le b_j\le K_j$ and $(b_j)$ is weakly null in $Y^*$.
For any $\xi\in H_{K_j}$, since $A_j=A_j^*$ and both $P_j$ and
$\mathbf 1-P_j$ are orthogonal projections, we obtain
\begin{align} \label{inner product}
 \langle T_j\xi,\xi\rangle_{H_{K_j}}
 \stackrel{\eqref{the def of Tj}}{=}\frac12\Bigl(\|P_j\xi\|^2_{H_{K_j}}+\|(\mathbf 1-P_j)A_j\xi\|^2_{H_{K_j}}\Bigr).
\end{align}
Noting that $\|P_j\xi\|_{H_{K_j}}\le\|\xi\|_{H_{K_j}}$ and
$\|(\mathbf 1-P_j)A_j\xi\|_{H_{K_j}}\le\|A_j\xi\|_{H_{K_j}}\le\|\xi\|_{H_{K_j}}$ by $\|A_j\|_{B({H_{K_j}})}\le1$,
we have that for any
$\xi\in H_{K_j}$, $0\le\langle T_j\xi,\xi\rangle_{H_{K_j}}\le\|\xi\|^2_{H_{K_j}}$, which implies that $0\le T_j\le \mathbf 1$.
Fix $t\in Y_+$.
By \eqref{GNS} and the definition of the inner product on $H_{K_j}$,
\begin{eqnarray} \label{inner product 2}
 0\le b_j(t)
 &\stackrel{\eqref{the def of Tj},~\eqref{GNS}}{=}&\left\langle T_j[t^{1/2}]_{K_j},[t^{1/2}]_{K_j}\right\rangle_{H_{K_j}}
 \le\left\|[t^{1/2}]_{K_j}\right\|^2_{H_{K_j}}\notag\\
 &\stackrel{K_j(y^*x)=\langle[x]_{K_j},[y]_{K_j}\rangle_{H_{K_j}}}{=}&K_j(t).
\end{eqnarray}
In particular, $0\le b_j\le K_j$.

Moreover,
\begin{eqnarray*}
 2b_j(u)&\stackrel{\eqref{GNS}}{=}&2\langle T_j[u^{1/2}]_{K_j},[u^{1/2}]_{K_j}\rangle_{H_{K_j}}\\
 &\stackrel{\eqref{the def of Tj}}{=}&\langle P_j[u^{1/2}]_{K_j},[u^{1/2}]_{K_j}\rangle_{H_{K_j}}+\langle A_j(\mathbf 1-P_j)A_j[u^{1/2}]_{K_j},[u^{1/2}]_{K_j}\rangle_{H_{K_j}}\\
 &\stackrel{A_j=A_j^*, ~\alpha_j=[u^{1/2}]_{K_j}}{=}&\|P_j\alpha_j\|^2_{H_{K_j}}+
          \|(\mathbf 1-P_j)A_j\alpha_j\|^2_{H_{K_j}}\stackrel{\eqref{normto0}}{\to}0.
\end{eqnarray*}
This together with $0\le b_j\le K_j$ and $(K_j)\prec (\widetilde h_n)$ yields that $(b_j)$ is weak$^*$ null in $Y^*$. 
Since $Y$ is positively Grothendieck, $(b_j)$ is weakly null in $Y^*$.

{\bf Step $5$}.
Given $x\in Y_+$, 
by $A_j=P_jA_j+(\mathbf 1-P_j)A_j$ and $\|A_j\|_{B(H_{K_j})}\le1$, we have that for all $s>0$,
\begin{eqnarray*}
 &&|F_j(x)|\\
 &\stackrel{\eqref{125}}{=}&\left|\left\langle A_j[x^{1/2}]_{K_j},[x^{1/2}]_{K_j}\right\rangle_{H_{K_j}}\right|\\
 &\le&\left(\left\|P_j[x^{1/2}]_{K_j}\right\|_{H_{K_j}}+\left\|(\mathbf 1-P_j)A_j[x^{1/2}]_{K_j}\right\|_{H_{K_j}}\right)\left\|[x^{1/2}]_{K_j}\right\|_{H_{K_j}}\\
 &\stackrel{K_j(y^*x)=\langle[x]_{K_j},[y]_{K_j}\rangle_{H_{K_j}}}{\le} &\sqrt{2\left(\left\|P_j[x^{1/2}]_{K_j}\right\|_{H_{K_j}}^2+\left\|(\mathbf 1-P_j)A_j[x^{1/2}]_{K_j}\right\|_{H_{K_j}}^2\right)}\sqrt{K_j(x)}\\
 &\stackrel{\eqref{inner product}, ~\eqref{inner product 2}}{=}&2\sqrt{b_j(x)K_j(x)}\\
 &\stackrel{2\sqrt {ab}\le a+b}{\le} & sK_j(x)+s^{-1}b_j(x).
\end{eqnarray*}
Then, for any $s>0$,
\[
 -\bigl(s K_j+s^{-1}b_j\bigr)
 \le F_j\le s K_j+s^{-1}b_j.
\]
By the fact that $(b_j)$ is weakly null in $Y^*$ and $x^{**}\in (Y^{**})_+$,
we obtain
\[
 \limsup_j|x^{**}(F_j)|\le \limsup_j x^{**}\left(sK_j+s^{-1}b_j\right)\stackrel{\substack{x^{**}\in B_{Y^{**}}\\x^{**}(b_j)\to0}}{\le}s\sup_j\|K_j\|\stackrel{\eqref{hn bounded},~\eqref{the def of Kj}}{<}\infty.
\]
Letting $s\downarrow0$, we obtain that $|x^{**}(F_j)|\to0$, which contradicts $x^{**}(F_j)\stackrel{\eqref{x**>delta}}{>}\delta$.
This completes the proof.
\end{proof}

We return to $X=E\cap F$.
For fixed $u \in F_+$,
we approximate the values $h_n(u)$ from below using positive functionals on $E + F$ .

\begin{lemma}\label{compatible}
Let $(h_n)\subseteq F^*_+$ be bounded and assume that
$h_n|_X\to0$ in $\sigma(X^*,X)$.
Given $u\in F_+$ and $\varepsilon>0$,
there exist  $N\ge1$ and 
positive functionals
$(w_n)_{n\ge N}\subseteq(E+F)^*$ such that
\begin{equation}\label{common-w}
 \sup_{n\ge N}\|w_n\|_{(E+F)^*}<\infty,\qquad w_n|_F\le h_n,
 \qquad (h_n-w_n|_F)(u)\le\varepsilon.
\end{equation}
\end{lemma}
\begin{proof}
Fix $u\in F_+$ and $\varepsilon>0$.

{\bf Step $1$}. Below, we  prove that
there exist $0<m<\infty$ and $N\ge 1$ such that for any $n\ge N$,
\begin{align}\label{common-cover}
 h_n(u)\le m\|e\|_E+h_n(f)+\varepsilon,
 \qquad e\in E_+,\ f\in F_+,\ u\le e+f.
\end{align}
Assume by contradiction that \eqref{common-cover} fails.
Then, there exist a subsequence of $(h_j)$, still denoted by $(h_j)$,    and $(e_j)\subseteq E_+$, $(f_j)\subseteq F_+$ such that
\begin{align} \label{contradi}
 h_j(u)-h_j(f_j)>2^j\|e_j\|_E+\varepsilon, \qquad u\le e_j+f_j.
\end{align}
Letting $C:=\sup_jh_j(u)<\infty$, we obtain that
$\sum_j\|e_j\|_E<\sum_j ^\infty C2^{-j}<\infty$ and $h_j(f_j)\le C$.

Set $v=u^{1/2}\in F^{(2)}$. 
Since $0\le e_j, f_j\le e_j+f_j$ and $0\le u\le e_j+f_j$, it follows from \cite[Proposition~2.2.24]{DPS} that there exist $c_{1,j},c_{2,j},c_{3,j}\in \cM$ with $\norm{c_{i,j}}_\cM\le1$, $1\le i\le 3$, such that
\begin{align*}
    c_{i,j}=c_{i,j}s(e_j+f_j), &\quad (e_j)^{1/2}=c_{1,j}(e_j+f_j)^{1/2}, ~~1\le i\le 3, \\ (f_j)^{1/2}=c_{2,j}(e_j+f_j)^{1/2}, &\quad v=c_{3,j}(e_j+f_j)^{1/2}.
\end{align*}
Note that $c_{1,j}^*c_{1,j}+c_{2,j}^*c_{2,j}=s\left(e_j+f_j\right)$ \cite[Lemma~3.6(ii)]{NelsonTT}.
Let $x_j:=c_{3,j}c_{1,j}^*(e_j)^{1/2}$ and $y_j:=c_{3,j}c_{2,j}^*(f_j)^{1/2}$.
In this case, 
\begin{eqnarray*}
x_j+y_j&=&
c_{3,j}c_{1,j}^*(e_j)^{1/2}+c_{3,j}c_{2,j}^*(f_j)^{1/2}\\
&=&
c_{3,j}c_{1,j}^*c_{1,j}(e_j+f_j)^{1/2}+c_{3,j}c_{2,j}^*c_{2,j}(e_j+f_j)^{1/2}\\
&=&c_{3,j}\left(c_{1,j}^*c_{1,j}+c_{2,j}^*c_{2,j}\right)(e_j+f_j)^{1/2}\\   
&=&c_{3,j}s\left(e_j+f_j\right)(e_j+f_j)^{1/2}=v.
\end{eqnarray*}

Let $E^{(2)},F^{(2)},X^{(2)}$ and $N_{h_j}$ be defined as in \eqref{the def of D(2)} and \eqref{the def of Nhj}, respectively.
Let $j_{h_j|_X}:X^{(2)}\to H_{h_j|_X}$ given by $j_{h_j|_X}(x)=[x]_{h_j|_X}=x+N_{h_j|_X}$, $x\in X^{(2)}$, where $H_{h_j|_X}$ denotes the Hilbert space completion of $X^{(2)}/N_{h_j|_X}$.
By $\norm{c_{i,j}}_\cM\le 1$, we have that $(c_{3,j}c_{i,j})^*(c_{3,j}c_{i,j})\le \mathbf 1$.
Then, $x_j^*x_j\le e_j$ and $y_j^*y_j\le f_j$.
In particular, $x_j\in E^{(2)}$ and $y_j\in F^{(2)}$.
Since $x_j=v-y_j$ and $v,y_j\in F^{(2)}$, we have $x_j\in X^{(2)}$.
Noting that
$$2u+2y_j^*y_j-x_j^*x_j\stackrel{\substack{x_j=v-y_j\\v^*v=u}}{=}(v+y_j)^*(v+y_j)\ge0,$$
we have that $x_j^*x_j\le 2u+2y_j^*y_j$.
Then, by  $\langle [x]_{h_j|_X}, [y]_{h_j|_X} \rangle_{H_{h_j|_X}}=h_j|_X(y^*x)$ for all $x,y\in X^{(2)}$,
\begin{eqnarray} \label{le 2sqrt C}
 \|j_{h_j|_X}(x_j)\|_{H_{h_j|_X}}=\sqrt{h_j|_X(x_j^*x_j)}&=&\sqrt{h_j(x_j^*x_j)}\notag\\
 &\stackrel{x_j^*x_j\le 2u+2y_j^*y_j}{\le}& \sqrt{h_j\left(2u+2y_j^*y_j\right)}\notag\\
 &\stackrel{\substack{y_j^*y_j\le f_j\\h_j(f_j),h_j(u)\le C}}{\le}& 2\sqrt C
 \end{eqnarray}
 and
 \begin{eqnarray} \label{lower bound}
 2{\operatorname{Re}} h_j(v^*x_j)=h_j(v^*x_j)+h_j(x_j^*v)
 &\stackrel{y_j=v-x_j}{=}&h_j(u)+h_j(x_j^*x_j)-h_j(y_j^*y_j)\notag\\
 &\stackrel{\substack{h_j(x_j^*x_j)\ge0\\y_j^*y_j\le f_j}}{\ge}&h_j(u)-h_j(f_j)
 \stackrel{\eqref{contradi}}{>}\varepsilon.
\end{eqnarray}

By $\|x_jx_j^*\|_E=\|x_j^*x_j\|_E\le\|e_j\|_E$ and $\sum_j\norm{e_j}_E<\infty$, we have that the series $d_0:=\sum_jx_jx_j^*$ converges in $E_+$.
Then, $d_0-x_jx_j^*\in E_+$ for every $j$.
For any $k>0$, let
\[
 p:=e^{k^2d_0-vv^*}[0,\infty),\qquad q:=\mathbf1-p.
\]
Then,
\begin{align} \label{the second}
 pvv^*p\le k^2pd_0p,\qquad
 qx_jx_j^*q\stackrel{d_0-x_jx_j^*\ge0}{\le} qd_0q\le k^{-2}qvv^*q.
\end{align}
This implies that
$(pv)(pv)^*\in E_+$, i.e. $pv\in E^{(2)}$.
By $v\in F^{(2)}$ and $p\in\M$, we have that $pv\in F^{(2)}$, and so $pv\in X^{(2)}$.
A combination of \eqref{the second} and \cite[Corollary 2.2.28]{DPS} yields $|x_j^*q|\le |k^{-1}v^*q|$, which  implies that
$\|qx_j\|_{ F^{(2)}}\le k^{-1}\|v\|_{ F^{(2)}}$.
Consequently, 
\begin{eqnarray*}
 &&|h_j(v^*x_j)|\\
 &\le&  |h_j\left(v^*px_j\right)|+|h_j\left(v^*qx_j\right)|\\
 &=& |h_j((pv)^*(px_j))|+|h_j((qv)^*(qx_j))|\\
 &\stackrel{pv,~x_j\in X^{(2)}}{\le}& \left|\langle j_{h_j|_X}(px_j),j_{h_j|_X}(pv) \rangle_{H_{h_j|_X}}\right|+\norm{h_j}_{F^*}\norm{(qv)^*(qx_j)}_F\\
 &\stackrel{\eqref{mixed}}{\le}& \norm{j_{h_j|_X}(px_j)}_{H_{h_j|_X}}\norm{j_{h_j|_X}(pv)}_{H_{h_j|_X}}+\norm{h_j}_{F^*}\norm{qv}_{F^{(2)}}\norm{qx_j}_{F^{(2)}}  \\
 &\stackrel{\|qx_j\|_{ F^{(2)}}\le k^{-1}\|v\|_{ F^{(2)}}}{\le}& \left(h_j|_X\left((px_j)^*(px_j)\right)\right)^{1/2}\left(h_j|_X\left((pv)^*(pv)\right)\right)^{1/2}+k^{-1}\sup_j\|h_j\|_{F^*}\|v\|_{ F^{(2)}}^2\\
 &\le&  \left(h_j|_X\left(x_j^*x_j\right)\right)^{1/2}\left(h_j|_X\left(v^*pv\right)\right)^{1/2}+k^{-1}\sup_j\|h_j\|_{F^*}\|v\|_{ F^{(2)}}^2\\
 &\stackrel{\eqref{le 2sqrt C}}{\le }& 2\sqrt C\,\left(h_j(v^*pv)\right)^{1/2}
       +k^{-1}\sup_j\|h_j\|_{F^*}\|v\|_{ F^{(2)}}^2.
\end{eqnarray*}
Noting that for fixed $k$, $v^*pv=(pv)^*(pv)\in X_+$ is independent of $j$,
by $h_j|_X\to 0$ in $\sigma(X^*,X)$, we have that $\left(h_j(v^*pv)\right)^{1/2}\to0$ as $j\to \infty$.
Taking $k\to \infty$, we obtain that $|h_j(v^*x_j)|\to 0$ as $j\to \infty$, which contradicts \eqref{lower bound}.
This proves \eqref{common-cover}.

{\bf Step $2$}. We aim to construct positive linear functionals on $ E+F$ satisfying \eqref{common-w}.
For any $n\ge N$, let
\begin{align} \label{the def of pn}
 p_n(z):=\inf\{m\|e\|_E+h_n(f):
            e\in E_+,\ f\in F_+,\ z\le e+f\}, \quad z\in (E+F)_{\rm sa}.
\end{align}
Note that $p_n$ is finite and nonnegative.

We claim that $p_n$ is sublinear.
Indeed, it suffices to prove that $p_n$ is subadditive.
Given $n\ge N$, let
$z_1,z_2\in(E+F)_{\rm sa}$, and let $\delta>0$.
By the definition of $p_n(z_i)$, for each $i=1,2$, there exist
$e_i\in E_+$ and $f_i\in F_+$ such that
\begin{align} \label{subadditive}
 z_i\le e_i+f_i,\qquad
 m\|e_i\|_E+h_n(f_i)<p_n(z_i)+\frac{\delta}{2}.
\end{align}
Then, $e_1+e_2\in E_+$, $f_1+f_2\in F_+$, and
\[
 z_1+z_2\le(e_1+e_2)+(f_1+f_2).
\]
This implies that
\begin{eqnarray*}
 p_n(z_1+z_2)
 &\le& m\|e_1+e_2\|_E+h_n(f_1+f_2)\\
 &\le& \bigl(m\|e_1\|_E+h_n(f_1)\bigr)
      +\bigl(m\|e_2\|_E+h_n(f_2)\bigr)\\
 &\stackrel{\eqref{subadditive}}{<}&p_n(z_1)+p_n(z_2)+\delta.
\end{eqnarray*}
Letting $\delta\downarrow0$, we obtain that
$p_n(z_1+z_2)\le p_n(z_1)+p_n(z_2)$.
This proves the claim.

Fix $n\ge N$. 
Moreover,
\begin{align} \label{-wn}
 h_n(u)-\varepsilon\stackrel{\eqref{common-cover}, ~\eqref{the def of pn}}{\le} p_n(u)\stackrel{\eqref{the def of pn}}{\le} h_n(u).
\end{align}
For any $z\in(E+F)_+$, we have that $p_n(-z)=0$.
Define $\lambda_n:\mathbb Ru\to\mathbb R$ by $\lambda_n(tu)=t p_n(u)$.
Note that $\lambda_n$ is a real linear functional.
Then,
\[
 \begin{aligned}
  \lambda_n(tu)&=t p_n(u)=p_n(tu), && t\ge0,\\
  \lambda_n(tu)&=t p_n(u)\le0=p_n(tu), && t<0.
 \end{aligned}
\]
By the  Hahn--Banach theorem \cite[Theorem~1.1]{B11},
$\lambda_n$ extends to a real linear functional on $(E+F)_{\rm sa}$,
still denoted by $\lambda_n$, satisfying
\begin{align} \label{exten of ln}
 \lambda_n(z)\le p_n(z), \quad z\in(E+F)_{\rm sa},
 \qquad \lambda_n(u)=p_n(u).
\end{align}
Define $w_n:E+F\to\mathbb C$ by
\[
 w_n(x+iy):=\lambda_n(x)+i\lambda_n(y),
 \qquad x,y\in(E+F)_{\rm sa}.
\]
Then, $w_n$ is complex linear, $w_n(u)\stackrel{\eqref{exten of ln}}{=}p_n(u)$, and
\[
 w_n(z)=-\lambda_n(-z)\stackrel{\eqref{exten of ln}}{\ge} -p_n(-z)=0,
 \qquad z\in(E+F)_+.
\]
Furthermore, 
\begin{align*}
\left(h_n-w_n|_F\right)(u)=h_n(u)-w_n(u)
=h_n(u)-p_n(u)\stackrel{\eqref{-wn}}{\le}\varepsilon.
\end{align*}
Note that
\[
 \begin{aligned}
  0\le w_n|_E(e)=\lambda_n(e)&\stackrel{\eqref{exten of ln}}{\le} p_n(e)\stackrel{\eqref{the def of pn}}{\le} m\|e\|_E, && e\in E_+,\\
  0\le w_n|_F(f)=\lambda_n(f)&\stackrel{\eqref{exten of ln}}{\le} p_n(f)\stackrel{\eqref{the def of pn}}{\le} h_n(f), && f\in F_+.
 \end{aligned}
\]
This together with \cite[Corollary 4.2.8]{DPS} yields that
\[
 \|w_n|_E\|_{E^*}\le m,\qquad 0\le w_n|_F\le h_n.
\]
For any $x=y+z\in E+F$ with $y\in E$ and $z\in F$, we have that
\begin{align*}
|w_n(x)|=|w_n(y+z)|&\le |w_n|_E(y)|+|w_n|_F(z)|\\
&\le m\norm{y}_E+\norm{h_n}_{F^*}\norm{z}_F\le \max\{m, \norm{h_n}_{F^*}\}\left(\norm{y}_E+\norm{z}_F\right).
\end{align*}
Taking infimum over all decompositions, we obtain that $\|w_n\|_{(E+F)^*}\le\max\{m,\|h_n\|\}$.
This shows that $\sup_{n\ge N}\|w_n\|_{(E+F)^*}\le\max\{m,\sup_n\|h_n\|\}<\infty$.
The proof is complete.
\end{proof}

Now, we apply above results to show that $E \cap F$ inherits the Grothendieck
property from $E$ and $F$.


\begin{proof}[Proof of Theorem \ref{main}]
{\bf Step $1$}. Below, We find linear functionals 
whose restrictions to $X$ are not weakly null and obtain a convex block sequence of $(h_n)$ as well as a control element in $F_+$
from Lemma \ref{control}.
By Theorem~\ref{criterion}, it suffices to consider positive
weak$^*$ null sequences in $X^*$. 
Assume by contradiction that there exists a positive weak$^*$-null sequence $(\phi_n)$ in $ X^*$ but not weakly null.
We claim that there exist $a_n \in E^*_+$ and $h_n \in F^*_+$
such that
\begin{align} \label{decompose}
    \phi_n(f) = a_n|_X(f) + h_n|_X(f), \quad \forall f \in X.
\end{align}
Indeed, 
given $f=f^*\in X$, let $p:=e^f{(0,\infty)}$ and  $f=f_+-f_-$ with $f_+,f_-\in X_+$.
For any $u\in E_+ $ and $v\in F_+$, we have $0\le f_+=pfp\le p(f+u)p$ and $f_+\le p(f+v)p$.
Then,
\[
    \phi_n(f)\le\phi_n(f_+)\le\|\phi_n\|_{X^*}\|f_+\|_X
    \le\|\phi_n\|_{X^*}\,\|(f+u,f+v)\|_{E \overset{\ell_\infty}{\oplus} F} .
\]
Applying \cite[Corollary~2.1 and~(2.6)]{MD00} to $E=\left(E \overset{\ell_\infty}{\oplus} F\right)_{\mathrm{sa}}$ and $K=E_+ \overset{\ell_\infty}{\oplus} F_+$, we obtain that there exists a linear positive continuous extension $ g_n$ on $\left(E \overset{\ell_\infty}{\oplus} F\right)_{\mathrm{sa}}$ of the mapping $(f,f) \mapsto \phi_n(f)$, $f=f^*\in X$ such that
$$\norm{g_n}=\norm{\phi_n}_{X^*}, \quad g_n(f,f)=\phi_n(f), \quad f=f^*\in X.$$
Consider the complex linear extension of $g_n$ given by
\[
    \psi_n(u+iv):=g_n(u)+i g_n(v),\qquad u,v\in \left(E \overset{\ell_\infty}{\oplus} F\right)_{\mathrm{sa}}.
\]
Observe that $\norm{\psi_n}=\norm{\phi_n}_{X^*}$.
Define
$a_n(\cdot) = \psi_n(\cdot,0)$, $h_n(\cdot) = \psi_n(0,\cdot)$.
Note that $a_n\in E^*_+$, $h_n\in F^*_+$ and
$$\phi_n(f) = a_n|_X(f) + h_n|_X(f), \quad \forall f \in X.$$
This proves the claim.

Moreover,
\[
    \|\phi_n\|_{X^*}
    = \|a_n\|_{E^*} + \|h_n\|_{F^*}.
\]
Indeed, for any $e \in E$, $f \in F$ with
$\max\{\|e\|_E, \|f\|_F\} = 1$,
$\psi_n(e,f) = a_n(e) + h_n(f)$ implies that
$\|\psi_n\| \le \|a_n\|_{E^*} + \|h_n\|_{F^*}$.
For the reverse inequality, given $\varepsilon > 0$,
choose $e_0 \in E$, $f_0 \in F$ with
$\|e_0\|_E, \|f_0\|_F \le 1$ such that
\[
    |a_n(e_0)| > \|a_n\|_{E^*} - \frac{\varepsilon}{2},
    \qquad
    |h_n(f_0)| > \|h_n\|_{F^*} - \frac{\varepsilon}{2}.
\]
Choose two unimodular scalars $\alpha$ and $\beta$ such that
$\alpha a_n(e_0) = |a_n(e_0)|$ and
$\beta h_n(f_0) = |h_n(f_0)|$.
Then, $(\alpha e_0,\beta f_0)$ belongs to the unit ball of
$E \overset{\ell_\infty}{\oplus} F$, and hence
\[
    \|\psi_n\|
    \ge |\psi_n(\alpha e_0,\beta f_0)|
    = |a_n(e_0)| + |h_n(f_0)|
    > \|a_n\|_{E^*} + \|h_n\|_{F^*} - \varepsilon.
\]
Letting $\varepsilon \downarrow 0$, we obtain
\[
    \|\phi_n\|_{X^*}
    = \|\psi_n\|
    = \|a_n\|_{E^*} + \|h_n\|_{F^*}.
\]

Note that
\[
 \sup_n\bigl(\|a_n\|_{E^*}+\|h_n\|_{F^*}\bigr)
 =\sup_n\|\phi_n\|_{X^*}\stackrel{\tiny\mbox{\cite[D.12]{AK06}}}{<}\infty.
\]
By $\phi_n\to 0$ in $\sigma(X^*,X)$ and \eqref{decompose}, we have that  for any $x\in X_+$, $0\le a_n(x),h_n(x)\le\phi_n(x)\to0$.
Then, $(a_n|_X)$ and $(h_n|_X)$ are weak$^*$ null in $X^*$.
Since $(\phi_n)$ fails to be weakly null in $X^*$, it follows from \eqref{decompose} that without loss of generality, we may assume that $(h_n|_X)$ is not weakly null.
Then, there exist a subsequence, denoted by $h_n$, an element $Y_0\in B_{X^{**}}$ and a constant $c\ne0$ such that
$Y_0(h_n|_X)\to c\ne0$.
Set $Y:=\frac{\overline c}{|c|}Y_0$ and $\delta:=\frac{|c|}4$. 
Then,
\[
 Y(h_{n}|_X)\longrightarrow |c|>2\delta.
\]
There exists a subsequence, denoted by $(h_n)$, such that
\begin{align} \label{<delta}
 \bigl|Y(h_n|_X)-|c|\bigr|<\delta,
 \qquad n\ge1.
\end{align}
For $D\in\{E,F\}$, define $Y_D\in D^{**}$ by $Y_D(b)=Y(b|_X), ~b\in D^*$.
Then, by $\|b|_X\|_{X^*}\le\|b\|_{D^*}$, we have that $\|Y_D\|_{D^{**}}\le1$.


Let 
$(H_j)\prec(h_n)$  and $u_F\in F_+$ be given by Lemma \ref{control} such that 
\begin{align} \label{the def of Hj}
H_j=\sum_{n\in I_j}\alpha_{j,n} h_n, \quad \alpha_{j,n}>0, \quad \sum_{n\in I_j}\alpha_{j,n}=1,
\end{align}
where $I_j$ are finite,  
mutually disjoint and $\max I_j<\min I_{j+1}$.
Since $(h_n|_X)$ is weak$^*$ null in $X^*$ and $j\le \min I_j\to \infty$ as $j\to \infty$, we have that for all $x\in X$,
\begin{align*}
    |H_j(x)|=\left|\sum_{n\in I_j}\alpha_{j,n} h_n(x)\right|\le \sup_{n\ge \min I_j}|h_n(x)|\to0,
\end{align*}
which implies that
$H_j|_X\xrightarrow{w^*}0$.
Note that
\begin{align} \label{<<delta}
 \left|Y_F(H_j)-|c|\right|\stackrel{\substack{\eqref{the def of Hj}\\Y_F(\cdot)=Y(\cdot|_X)}}{=}\left|\sum_{n\in I_j}\alpha_{j,n} Y(h_n|_X)-|c|\right|\stackrel{\eqref{<delta}}{<} \delta \sum_{n\in I_j}\alpha_{j,n}=\delta.
\end{align}

{\bf Step $2$}.
We aim to prove that there exist $\eta>0$ and $N\ge1$ such that 
\begin{equation}\label{at-Y}
 G\in\operatorname{co}\{H_j:j\ge N\},\quad 0\le r\le G,\quad r(u_F)\le\eta
 \quad\Longrightarrow\quad |Y_F(r)|<\delta.
\end{equation}
Otherwise, 
let $N_1=1$ and $\eta=1$. 
There exist $T_1\in\operatorname{co}\{H_j:j\ge1\}$ and
$r_1\in(F^*)_+$ such that
\[
 T_1=\sum_{j\in D_1}\alpha_{1,j}'H_j,~
 \alpha_{1,j}'>0,~ \sum_{j\in D_1}\alpha_{1,j}'=1, \quad 0\le r_1\le T_1,\quad r_1(u_F)\le1,\quad |Y_F(r_1)|\ge\delta,
\]
where $D_1$ is a nonempty finite subset of $\{j:j\ge 1\}$.
Now, consider $N_2=\max D_1+1$ and $\eta=1/2$.
There exist $T_2\in\operatorname{co}\{H_j:j\ge \max D_1+1\}$ and $r_2\in(F^*)_+$ such that
\[
 T_2=\sum_{j\in D_2}\alpha_{2,j}'H_j,~
 \alpha_{2,j}'>0,~ \sum_{j\in D_2}\alpha_{2,j}'=1, \quad0\le r_2\le T_2,\quad r_2(u_F)\le\tfrac12,\quad
 |Y_F(r_2)|\ge\delta,
\]
where $D_2$ is a nonempty finite subset of $\{j:j\ge \max D_1+1\}$.
Note that $D_1$ and $D_2$ are disjoint.
Arguing inductively, we obtain  a convex block sequence $(T_k)\prec(H_j)$ and $r_k\in F^*_+$ satisfying
\begin{align} \label{contradicts}
 0\le r_k\le T_k,\qquad r_k(u_F)\le 1/k,\qquad |Y_F(r_k)|\ge\delta.
\end{align}
By the fact that $(T_k)\prec(H_j)$ and $\left((h_n),(H_j),u_F\right)$ satisfy \eqref{control-eq}, we have that
$r_k\xrightarrow{w^*}0$ in $F^*$.
Since $F$ is Grothendieck, we have that $r_k\xrightarrow{w}0$, and so $Y_F(r_k)\to0$,
which contradicts \eqref{contradicts}.
Then, \eqref{at-Y} holds.

{\bf Step $3$}. We intend to construct functionals on $E+F$ approximating $H_j(u_F)$ and attain a convex block sequence of  above functionals
together with a control element in $E_+$ from Lemma \ref{control}.
Applying Lemma~\ref{compatible} to $(h_n)=(H_j)$, $u=u_F$ and $\varepsilon=\eta$, we obtain that there exist $N_1\ge N$ and
$(w_j)_{j\ge N_1}\subseteq(E+F)^*_+$ such that
\begin{align} \label{rj(uF)le eta}
 C:=\sup_{j\ge N_1}\|w_j\|_{(E+F)^*}<\infty,\quad
 r_j:=H_j-w_j|_F\ge0,\quad r_j(u_F)\le\eta, \quad j\ge N_1.
\end{align}
Since $0\le w_j|_X\le H_j|_X$ and $H_j|_X\xrightarrow{w^*}0$, we have $w_j|_X\xrightarrow{w^*}0$.
A combination of \eqref{at-Y} and \eqref{rj(uF)le eta} yields that $|Y_F(r_j)|<\delta$ for all $j\ge N_1$.
Then, for all $j\ge N_1$,
\begin{equation}\label{Y-row}
 \bigl|Y(w_j|_X)-|c|\bigr|
 \stackrel{\substack{Y_F(\cdot)=Y(\cdot|_X)\\\eqref{rj(uF)le eta}}}{=}\bigl|Y_F(H_j)-|c|-Y_F(r_j)\bigr|
 \le \bigl|Y_F(H_j)-|c|\bigr|+|Y_F(r_j)|\stackrel{\eqref{<<delta}}{<}2\delta.
\end{equation}

Applying Lemma~\ref{control} to $(g_n)=(w_j|_E)_{j\ge N_1}$, we obtain that there exist $u_E\in E_+$
and a convex block sequence $(A_m)_{m\ge1}$ of
$(w_j|_E)_{j\ge N_1}$ satisfying \eqref{control-eq}, 
where $(J_m)$ are mutually disjoint finite subsets of $\{N_1,N_1+1,\ldots\}$ and
\[
 A_m=\sum_{j\in J_m}\alpha_{mj} w_j|_E,\qquad
 \sum_{j\in J_m}\alpha_{mj}=1,\qquad
 \max J_m<\min J_{m+1}, \quad m\ge1.
\]
Define
\[
 V_m:=\sum_{j\in J_m}\alpha_{mj}w_j\in(E+F)^*_+,\qquad
 B_m:=\sum_{j\in J_m}\alpha_{mj}H_j\in F^*_+, \quad m\ge1.
\]
Note that $V_m|_E=A_m$ for any $m\ge1$ and $(B_m)_{m\ge1}\prec(H_j)_{j\ge N_1}$.
Since $0\le w_j|_F\le H_j$ for any $j\ge N_1$ (see \eqref{rj(uF)le eta}), we have
\begin{align} \label{Vmle Bm}
 0\le V_m|_F
 =\sum_{j\in J_m}\alpha_{mj}w_j|_F
 \le\sum_{j\in J_m}\alpha_{mj}H_j=B_m.
\end{align}
By $\|w_j\|_{(E+F)^*}\le C$ for any $j\in J_m$ (see \eqref{rj(uF)le eta}), we have that $\|V_m\|_{(E+F)^*}\le C$.

Fix $x\in X$. 
By $w_j|_X\xrightarrow{w^*}0$, we have that $w_j(x)\to0$ as $j\to\infty$.
This together with
$\min J_m\to\infty$ as $m\to\infty$ yields that
\[
 |V_m(x)|
 =\left|\sum_{j\in J_m}\alpha_{mj}w_j(x)\right|
 \le\sup_{j\ge\min J_m}|w_j(x)|\longrightarrow0 ~\text{as}~ m\to \infty.
\]
This implies that $V_m|_X\xrightarrow{w^*}0$ in $X^*$. 
For any $m\ge1$,
\begin{align} \label{<2delta}
 \begin{aligned}
 \bigl|Y(V_m|_X)-|c|\bigr|
 &=\left|\sum_{j\in J_m}\alpha_{mj}
          \bigl(Y(w_j|_X)-|c|\bigr)\right|\\
 &\le\sum_{j\in J_m}\alpha_{mj}\bigl|Y(w_j|_X)-|c|\bigr|
 \stackrel{\eqref{Y-row}}{<}2\delta\sum_{j\in J_m}\alpha_{mj}
 =2\delta.
 \end{aligned}
\end{align}


{\bf Step $4$}. We aim to construct a weakly null sequence in $\cM^*$ and choose a convex block sequence  with norms tending to zero. 
Let $v:=u_F^{1/2}$ and $\lambda_m(a):=V_m(v^*au_E^{1/2})$ for any $a\in\M$.
Note that $v\in F^{(2)}\subseteq (E+F)^{(2)}$ and $u_E^{1/2}\in E^{(2)}\subseteq (E+F)^{(2)}$.
Then, by \cite[Lemma 4.2.4]{DPS}, we have that
$v^*au_E^{1/2}\in E+F$ and 
$$|\lambda_m(a)|=\left|V_m(v^*au_E^{1/2})\right|\le \norm{V_m}_{(E+F)^*}\norm{v^*au_E^{1/2}}_{E+F}\stackrel{\eqref{mixed}, ~\eqref{the def of D(2)}}{\le} C\|u_F\|_F^{1/2}\|u_E\|_E^{1/2}\norm{a}_\cM.$$
Given $a\in\M$, set $z:=au_E^{1/2}$,
$p:=e^{zz^*-vv^*}([0,\infty))$ and $q:=\mathbf1-p$. 
Then,
\[
 0\le(pv)(pv)^*\le pzz^*p\in E_+,\qquad
 0\le(qz)(qz)^*\le qvv^*q\in F_+.
\]
This together with $pv\in F^{(2)}$ and $qz\in E^{(2)}$ yields that
$pv,qz\in X^{(2)}$. In particular, $v^*pv=(pv)^*(pv),z^*qz=(qz)^*(qz)\in X_+$.
Then, by $v^*z=(pv)^*(pz)+(qv)^*(qz)$ and $V_m|_X\xrightarrow{w^*}0$ in $X^*$, we obtain
\begin{eqnarray*}
 |\lambda_m(a)|&=&|V_m\left(v^*z\right)|\\
 &\le& |V_m\left(v^*pz\right)|+|V_m\left(v^*qz\right)|\\
 &=&|V_m\left((pv)^*(pz)\right)|+|V_m\left((qv)^*qz\right)|\\
 &\stackrel{\tiny\mbox{\cite[(4.9)]{DPS}}}{\le}& \sqrt{V_m(z^*pz)V_m(v^*pv)}
       +\sqrt{V_m(z^*qz)V_m(v^*qv)}\\
 &\stackrel{\substack{0\le z^*pz\le z^*z\\0\le v^*qv\le u_F}}{\le}&\sqrt{C\|z^*z\|_E}\,\sqrt{(V_m|_X)(v^*pv)}
       +\sqrt{C\|u_F\|_F}\,\sqrt{(V_m|_X)(z^*qz)}
 \longrightarrow 0.
\end{eqnarray*}
Then, $\lambda_m\xrightarrow{w^*}0$ in $\M^*$.
Since $\M$ is Grothendieck \cite[Corollary~7]{P94}, we have that
$\lambda_m\xrightarrow{w}0$ in $\cM^*$.

Arguing as \eqref{the def of Lambdaj}, we have that there exists a convex block sequence $(\Lambda_k)_{k\ge1}\prec (\lambda_m)_{m\ge1}$ satisfying
\begin{align} \label{the def of Lambda}
 \Lambda_k=\sum_{m\in S_k}\beta_{km} \lambda_m,\quad
 \sum_{m\in S_k}\beta_{km}=1,\quad
 \max S_k<\min S_{k+1}, \quad \|\Lambda_k\|_{\M^*}<2^{-k}, \quad k\ge1.
\end{align}
Set
\begin{align} \label{the def of Wk}
 W_k:=\sum_{m\in S_k}\beta_{km}V_m\in (E+F)^*_+,\qquad L_k:=\sum_{m\in S_k}\beta_{km}B_m\in F^*_+, \quad k\ge1.
\end{align}
Then, $(W_k|_E)\prec(V_m|_E)$, $(L_k)\prec(B_m)\prec(H_j)_{j\ge N_1}$, and
$0\le W_k|_F\le L_k$ by \eqref{Vmle Bm}.

{\bf Step $5$}.
Fix $k\ge1$.
Let $j_{W_k}:(E+F)^{(2)}\to H_{W_k}$ given by $j_{W_k}(x)=[x]_{W_k}=x+N_{W_k}$, where $H_{W_k}$ denotes the Hilbert space completion of $(E+F)^{(2)}/N_{W_k}$.
For any $a\in \cM$, set $L_a[x]_{W_k}:=[ax]_{W_k}$ and we denote its extension to $H_{W_k}$ by $\pi_{W_k}(a)$.
Moreover, the map $\pi_{W_k}:\cM \to B(H_{W_k})$ is a unital $*$-representation.
Applying Lemma~\ref{hilbert} to $\xi=j_{W_k}(u_E^{1/2})$, $\eta=j_{W_k}(u_F^{1/2})$ and $\cA=\pi_{W_k}(\M)'$,
we obtain that there exists a projection $P_k\in\pi_{W_k}(\M)'$ such that
\begin{align} \label{ref 6.1}
 &\quad \quad \Biggl\|P_k\left(j_{W_k}\left(u_E^{1/2}\right)\right)\Biggr\|_{H_{W_k}}^2+\Biggl\|\left(\mathbf 1-P_k\right)\left(j_{W_k}\left(u_F^{1/2}\right)\right)\Biggr\|^2_{H_{W_k}}\notag\\
 &\le\sup_{a\in B_{\pi_{W_k}(\M)''}}\left|\left\langle a\left(j_{W_k}\left(u_E^{1/2}\right)\right),j_{W_k}\left(u_F^{1/2}\right)\right\rangle_{H_{W_k}}\right|\notag\\
 &=\sup_{a\in B_{\cM}}\left|\left\langle\pi_{W_k}(a)\left(j_{W_k}\left(u_E^{1/2}\right)\right),j_{W_k}\left(u_F^{1/2}\right)\right\rangle_{H_{W_k}}\right|.
\end{align}
Set $b_k:=({W_k})_{P_k}$ as in \eqref{GNS}. 
Since $A\in \pi_{W_k}(\M)'_h\mapsto ({W_k})_A$ is linear and positive, we have that
$b_k\ge0$ and $W_k-b_k=({W_k})_{\mathbf 1-P_k}\ge0$.
Observe that $b_k=(W_k)_{P_k}\le (W_k)_\mathbf 1=W_k$.
Note that by \eqref{GNS}, we have that 
\begin{align} \label{P norm}
b_k(u_E)=\Biggl\|P_k\left(j_{W_k}\left(u_E^{1/2}\right)\right)\Biggr\|^2_{H_{W_k}},\quad (W_k-b_k)(u_F)=\Biggl\|\left(\mathbf 1-P_k\right)j_{W_k}\left(u_F^{1/2}\right)\Biggr\|^2_{H_{W_k}}. 
\end{align}
For any $a\in B_\M$,
 \begin{eqnarray} \label{=Lambda k(a)}
  \left\langle\pi_{W_k}(a)\left(j_{W_k}\left(u_E^{1/2}\right)\right),j_{W_k}\left(u_F^{1/2}\right)\right\rangle_{H_{W_k}}
  &=&\left\langle j_{W_k}\left(a u_E^{1/2}\right),
                j_{W_k}\left(u_F^{1/2}\right) \right\rangle_{H_{W_k}} \notag\\
  &\stackrel{\eqref{GNS}}{=}&W_k\left(u_F^{1/2}a u_E^{1/2}\right)\notag\\
  &=&\sum_{m\in S_k}\beta_{km}V_m\left(u_F^{1/2}a u_E^{1/2}\right)\notag\\
    &=&\sum_{m\in S_k}\beta_{km}\lambda_m(a)\stackrel{\eqref{the def of Lambda}}{=}\Lambda_k(a).
 \end{eqnarray}
A combination of \eqref{ref 6.1}, \eqref{P norm} and \eqref{=Lambda k(a)} yields that
\[
 b_k(u_E)+(W_k-b_k)(u_F)\le\|\Lambda_k\|_{\M^*}\stackrel{\eqref{the def of Lambda}}{<}2^{-k}.
\]
This implies that $b_k(u_E)\to 0$ and $(W_k-b_k)(u_F)\to0$ as $k\to \infty$.

Now, we have that $(W_k|_E)\prec (V_m|_E)=(A_m)$, $(L_k)\prec (H_j)_{j\ge N_1}$ and $ \left( (w_j|_E)_{j\ge N_1}, (A_m),u_E\right)$, $\left((h_n),(H_j),u_F\right)$ both satisfy \eqref{control-eq},
and
\[
 0\le b_k|_E\le W_k|_E,\quad b_k|_E(u_E)\to0,\quad
 0\le (W_k-b_k)|_F\le L_k,\quad (W_k-b_k)|_F(u_F)\to0. 
\]
By \eqref{control-eq} in $E$ and $F$, respectively,
$b_k|_E\xrightarrow{w^*}0$ in $E^*$ and $(W_k-b_k)|_F\xrightarrow{w^*}0$ in $F^*$.
Since $E$ and $F$ are Grothendieck, it follows that $b_k|_E\xrightarrow{w}0$ and
$(W_k-b_k)|_F\xrightarrow{w}0$. 
Since $\sum_{m\in S_k}\beta_{km}=1$ for all $k\ge1$, we have that
\begin{align} \label{>-2delta}
 \bigl|Y(W_k|_X)-|c|\bigr|
 \stackrel{\eqref{the def of Wk}}{=}\Bigl|\sum_{m\in S_k}\beta_{km}\bigl(Y(V_m|_X)-|c|\bigr)\Bigr|\stackrel{\eqref{<2delta}}{<}2\delta.
\end{align}
By $|c|=4\delta$, $W_k|_X=b_k|_X+(W_k-b_k)|_X$, $b_k|_E\xrightarrow{w}0$ and $(W_k-b_k)|_F\xrightarrow{w}0$, we obtain
\[
 2\delta=|c|-2\delta\stackrel{\eqref{>-2delta}}{<}|Y(W_k|_X)|
 \stackrel{Y_D(\cdot)=Y(\cdot|_X), ~D=E,F}{=}\bigl|Y_E(b_k|_E)+Y_F((W_k-b_k)|_F)\bigr|\longrightarrow0,
\]
a contradiction.
This completes the proof.
\end{proof}

{\bf Statement on the use of artificial intelligence}

This work grew out of discussions between the authors and ChatGPT 6.0 Astra. During these discussions, we found that the \(\ell_\infty\)-sum of Grothendieck spaces considered in this paper (see Corollary \ref{cor:uc-product}) also has the Grothendieck property. This observation led us to explore the Grothendieck property of  the intersection of two Grothendieck spaces. The authors chose the research questions, guided the discussions, and take full responsibility for the results presented in this paper.

{\bf Data availability}
Our work follows a theoretical and mathematical approach, and as such does not generate any datasets.
  
 {\bf Conflict of interest}
 The authors have no conflict of interest to declare that are relevant to the content of this article.

\end{document}